# WEIL–PETERSSON HOMEOMORPHISMS, MINIMAL LAGRANGIAN DIFFEOMORPHISMS, AND MAXIMAL SURFACES IN ANTI-DE SITTER SPACE

FARID DIAF, ALEX MORIANI, RYM SMAÏ, GRAHAM ANDREW SMITH, AND ENRICO TREBESCHI

ABSTRACT. In this paper, we study the class of Weil–Petersson circle homeomorphisms from the point of view of three-dimensional anti-de Sitter space $\mathbf{AdS}^{2,1}$. We show that a homeomorphism $\varphi : \mathbf{RP}^1 \to \mathbf{RP}^1$ is Weil–Petersson if and only if its graph, viewed as a curve in the boundary at infinity of $\mathbf{AdS}^{2,1}$, is the asymptotic boundary of a complete maximal spacelike surface in $\mathbf{AdS}^{2,1}$ with finite renormalized area. As an application, we obtain the following AdS-independent result in Teichmüller theory: a homeomorphism is Weil–Petersson if and only if its minimal lagrangian extension to $\mathbf{H}^2$ has square-integrable Beltrami differential. We also provide two further new technical characterizations, which we believe to be of independent interest, and which are essential for the proofs of our main results.

## Contents

## 1. Introduction

*Weil–Petersson Teichmüller space* is a remarkable space of homeomorphisms of the circle which has in recent years attracted ever more interest from mathematicians and theoretical physicists. Introduced independently under different names by Cui in [Cui00] and Guo in [Guo00], its significance for *universal Teichmüller theory* was established by Takhtajan–Teo in [TT06]. This space is endowed with a rich geometric structure and, in particular, carries an essentially unique, complete, homogeneous Kähler metric, arising as the completion of the Kähler metric defined over the space of smooth diffeomorphisms of $\mathbf{RP}^1$ by Bowick–Rajeev in [BR87]. Weil–Petersson Teichmüller space has since found applications across several areas of contemporary mathematics, including geometric function theory, operator theory, hyperbolic geometry, computer vision, stochastic processes such as Schramm–Loewner evolution, and so on (see, for example, [SM06, She18, Wan19, Bis22, BBVPW25, LTW26]).

In his recent, groundbreaking work [Bis25], Bishop provided no less than two dozen different characterizations of Weil–Petersson Teichmüller space, highlighting, in particular, the role played in this theory by the geometry of 3-dimensional hyperbolic space $\mathbf{H}^3$. The purpose of the present paper is to introduce a new perspective based on the geometry of (2, 1)-dimensional anti-de Sitter space $\mathbf{AdS}^{2,1}$, widely viewed as the lorentzian cousin of $\mathbf{H}^3$, having richer geometric structure. Indeed, since the seminal work [Mes07] of Mess, a deep connection has been revealed between (2, 1)-dimensional anti-de Sitter geometry and classical Teichmüller theory (see, for example, [KS07, BS10, BS12, BMS13, BS18, Sep19, BS20, Tre25].)

In the following sections, we present two new characterizations of the Weil–Petersson Teichmüller space in terms of anti–de Sitter geometry. The first is given in terms of the Beltrami differentials of minimal lagrangian diffeomorphisms of the hyperbolic plane, and the second, of which the first is a consequence, is formulated in terms of complete maximal spacelike surfaces in $\mathbf{AdS}^{2,1}$. We also provide two further, more technical, characterizations in terms of a lorentzian analogue of Jones' beta numbers (see [Jon90]), which we believe to be of independent interest.

1.1. **Minimal lagrangian diffeomorphisms of $\mathbf{H}^2$.** A classical problem in Teichmüller theory is to find *conformally natural extensions* to the hyperbolic plane of homeomorphisms defined over the circle. Notable examples include the Douady–Earle extension [DE86], the earthquake extension [Thu86], the harmonic extension [Mar17], and so on (we refer the reader to the survey [Hu22] for a deeper discussion). Relatively little is known however about the relationship between such extensions and Weil–Petersson Teichmüller space. Our first main result addresses this in the case of the *minimal lagrangian extension*.

Throughout the sequel, we will identify the circle with 1-dimensional real projective space $\mathbf{RP}^1$, and we will restrict attention to homeomorphisms $\varphi : \mathbf{RP}^1 \to \mathbf{RP}^1$ which preserve orientation. We will use the following characterization of the Weil–Petersson class, provided by Shen in [She18, Theorem 1.1].

**Definition 1.1.1.** A homeomorphism $\varphi : \mathbf{RP}^1 \to \mathbf{RP}^1$ is *Weil–Petersson* whenever it has $L^1$ distributional derivative and satisfies

$$\log(\varphi') \in H^{1/2}(\mathbf{RP}^1) \, .$$

*Remark* 1.1.2. For any interval $I \subseteq \mathbf{R}$, a continuous function $f : I \to \mathbf{R}$ has $L^1$ distributional derivative (with respect to Lebesgue measure) if and only if it is absolutely continuous.

It is known that a quasisymmetric homeomorphism $\varphi : \mathbf{RP}^1 \to \mathbf{RP}^1$ is Weil–Petersson if and only if it admits a quasiconformal extension to the hyperbolic plane with square integrable Beltrami differential (see [Cui00, TT06]). It was likewise proven in [Cui00] that the Beltrami differential of the Douady–Earle extension of any such homeomorphism is also square integrable. Our first result (Theorem A) provides an analogous characterization of Weil-Petersson homeomorphisms in terms of minimal lagrangian diffeomorphisms.

The concept of minimal lagrangian diffeomorphism was introduced independently by Labourie in [Lab92] and Schoen in [Sch93] as a modification of the class of harmonic diffeomorphisms into one preserved by the operation of inversion. Let $\mathbf{H}^2$ denote 2-dimensional hyperbolic space. A diffeomorphism $\Phi : \mathbf{H}^2 \to \mathbf{H}^2$ is said to be *minimal lagrangian* whenever its graph is a minimal lagrangian surface[1] in the cartesian product $\mathbf{H}^2 \times \mathbf{H}^2$. Using the natural identification of $\mathbf{RP}^1$ with the ideal boundary of $\mathbf{H}^2$, we say that a minimal lagrangian diffeomorphism $\Phi : \mathbf{H}^2 \to \mathbf{H}^2$ *extends* the homeomorphism $\varphi : \mathbf{RP}^1 \to \mathbf{RP}^1$ whenever $\Phi$ extends to a homeomorphism of $\mathbf{H}^2 \cup \mathbf{RP}^1$ to itself whose restriction to $\mathbf{RP}^1$ coincides with $\varphi$.

In [Lab92] and [Sch93], Labourie and Schoen independently proved the existence of a unique minimal lagrangian extension for every homeomorphism $\varphi : \mathbf{RP}^1 \to \mathbf{RP}^1$ which is *equivariant* under a discrete and faithful action of a compact surface group. This result was extended to all *quasisymmetric* homeomorphisms by Bonsante–Schlenker in [BS10] and to *general* homeomorphisms by the fifth author in [Tre24b] (see also Appendix A).

There is growing evidence to suggest that the geometry of the minimal lagrangian extension is closely tied to that of the initial homeomorphism. Indeed, in [BS10], Bonsante–Schlenker showed that the minimal lagrangian extension is quasiconformal if and only if the initial homeomorphism is quasisymmetric. By a compactness argument, it is straightforward to show (see Theorem 5.2.2 of [Sep24]) that the minimal lagrangian extension is *asymptotically conformal* if and only if the initial homeomorphism is *symmetric*. Recall that the *Beltrami differential* of a diffeomorphism $\Phi : \mathbf{H}^2 \to \mathbf{H}^2$ is given by

$$\mu_\Phi := \frac{\overline{\partial}\Phi}{\partial\Phi} \, .$$

For any homeomorphism $\varphi : \mathbf{RP}^1 \to \mathbf{RP}^1$, we denote by $\mu_\varphi$ the Beltrami differential of its minimal lagrangian extension. Our first result provides a characterization of Weil–Petersson homeomorphisms in terms of the Beltrami differential of their minimal lagrangian extensions.

**Theorem A.** *A homeomorphism $\varphi : \mathbf{RP}^1 \to \mathbf{RP}^1$ is Weil–Petersson if and only if the Beltrami differential $\mu_\varphi$ of its minimal lagrangian extension is square-integrable over $\mathbf{H}^2$, that is, if and only if*

$$\int_{\mathbf{H}^2} \|\mu_\varphi\|^2 \, d\, \mathrm{Area}_{\mathbf{H}^2} < \infty \, ,$$

*where* $d\, \mathrm{Area}_{\mathbf{H}^2}$ *denotes the canonical area form of $\mathbf{H}^2$.*

The proof of Theorem A is explained in Section 1.6.

*Remark* 1.1.3. This answers in the affirmative a question raised by Seppi in [Sep24].

---

[1]The canonical symplectic form on $\mathbf{H}^2 \times \mathbf{H}^2$ is $\pi_1^*\omega \ominus \pi_2^*\omega$, where $\omega$ is the area form of the hyperbolic metric. Thus, the graph of $\Phi$ is lagrangian with respect to this form if and only if $\Phi$ is area-preserving.

*Remark* 1.1.4. Note that, in contrast to [She18], Theorem A does not require the minimal lagrangian extension to be quasiconformal. More precisely, we show that when the boundary map is a homeomorphism, quasiconformality in fact follows from the square integrability of the Beltrami differential.

1.2. **Complete maximal spacelike surfaces in $\mathbf{AdS}^{2,1}$.** Theorem A is proven using the geometry of anti-de Sitter space $\mathbf{AdS}^{2,1}$. In particular, this result has a natural reformulation in terms of complete maximal spacelike surfaces in $\mathbf{AdS}^{2,1}$, which we now describe. Throughout the sequel, we view $\mathbf{AdS}^{2,1}$ as the projectivized quadric

$$\mathbf{AdS}^{2,1} := \{(x_1, x_2, x_3, x_4) \in \mathbf{R}^4 \mid x_1^2 + x_2^2 - x_3^2 - x_4^2 = -1\}/\{\pm 1\} \ .$$

This is a complete lorentzian 3-manifold of constant sectional curvature equal to $-1$. The ideal boundary of $\mathbf{AdS}^{2,1}$, known as $(1,1)$-dimensional Einstein space $\mathbf{Ein}^{1,1}$, naturally identifies with the cartesian product $\mathbf{RP}^1 \times \mathbf{RP}^1$. Under this correspondence, the graph of any homeomorphism $\varphi : \mathbf{RP}^1 \to \mathbf{RP}^1$ is a simple closed *acausal* curve $\Lambda_\varphi$ in $\mathbf{Ein}^{1,1}$, and, conversely, every such curve in $\mathbf{Ein}^{1,1}$ arises in this way. This correspondence is explained in greater detail in Section 2.3.

In [BS10], Bonsante–Schlenker showed that the graph in $\mathbf{Ein}^{1,1}$ of every *quasisymmetric* homeomorphism of $\mathbf{RP}^1$ bounds a unique complete *maximal* spacelike surface $\Sigma$ in $\mathbf{AdS}^{2,1}$, that is, a complete spacelike surface of zero mean curvature, and in [LTW24], Labourie–Toulisse–Wolf extended this result to *general* homeomorphisms. For any homeomorphism $\varphi : \mathbf{RP}^1 \to \mathbf{RP}^1$, we denote by $\Sigma_\varphi$ the unique complete maximal spacelike surface bounded by $\Lambda_\varphi$ and, by a slight abuse of terminology, we say that $\varphi$ *bounds* $\Sigma_\varphi$.

The *renormalized area* of any maximal spacelike surface $\Sigma$ in $\mathbf{AdS}^{2,1}$ is defined by

$$\mathcal{A}(\Sigma) := \int_\Sigma \|\mathrm{I\!I}_\Sigma\|^2 \mathrm{d\,Area}_\Sigma \ ,$$

where $\mathrm{I\!I}_\Sigma$ and $\mathrm{d\,Area}_\Sigma$ denote, respectively, its second fundamental form and its area form. This is the natural analogue of the notion of renormalized area defined for minimal surfaces in Poincaré–Einstein spaces by Graham–Witten in [GW99], and is justified as the constant term in a suitable asymptotic expansion of the areas of large balls in $\Sigma$ (we refer the reader to [AM10] for a more detailed discussion). We can now state our second main result.

**Theorem B.** *A homeomorphism $\varphi : \mathbf{RP}^1 \to \mathbf{RP}^1$ is Weil–Petersson if and only if the unique complete maximal spacelike surface in $\mathbf{AdS}^{2,1}$ that it bounds has finite renormalized area.*

The proof of Theorem B is explained in Section 1.6.

*Remark* 1.2.1. In [SST23], in collaboration with A. Seppi and J. Toulisse, the fourth author showed that a complete maximal spacelike surface $\Sigma$ in $\mathbf{AdS}^{2,1}$ has finite renormalized area whenever its boundary curve is $C^{3,\alpha}$. Theorem B completes this result by providing a full characterization of complete maximal spacelike surfaces in $\mathbf{AdS}^{2,1}$ with finite renormalized area in terms of their asymptotic boundaries.

*Remark* 1.2.2. Note that a complete spacelike surface in $\mathbf{AdS}^{2,1}$ is necessarily properly embedded (see, for example, Lemma 3.11 of [SST23]). In general, a properly embedded spacelike surface in $\mathbf{AdS}^{2,1}$ need not be complete. However, in [Tre24a], the fifth author shows that every properly embedded *maximal* spacelike surface in $\mathbf{AdS}^{2,1}$ is complete. The statement of Theorem B thus remains valid with the word "complete" replaced by "properly embedded".

1.3. **Homeomorphisms, diffeomorphisms, and surfaces.** Theorems A and B together reflect the deep connection between 3-dimensional anti-de Sitter geometry and Teichmüller theory. This relationship, first observed by Mess in [Mes07], arises from the structure of *the Gauss map of spacelike surfaces in* $\mathbf{AdS}^{2,1}$, which yields surfaces in $\mathbf{H}^2 \times \mathbf{H}^2$, which are viewed in turn as graphs of functions defined over subsets of $\mathbf{H}^2$. In this manner, Mess showed that certain, well-chosen spacelike surfaces in $\mathbf{AdS}^{2,1}$ yield natural extensions to $\mathbf{H}^2$ of homeomorphisms of $\mathbf{RP}^1$. He thereby obtained a simple proof of Thurston's earthquake theorem by showing that, for any homeomorphism $\varphi : \mathbf{RP}^1 \to \mathbf{RP}^1$, the desired earthquake maps are induced by the boundary components of the convex hull in $\mathbf{AdS}^{2,1}$ of its graph.

It was later observed by Bonsante–Schlenker in [BS10] that the same procedure identifies suitable maximal spacelike surfaces in $\mathbf{AdS}^{2,1}$ with minimal lagrangian diffeomorphisms of $\mathbf{H}^2$. This yields the following picture, completed by the fifth author in [Tre24b] (see also Appendix A), and illustrated in Figure 1: a homeomorphism of $\mathbf{RP}^1$ determines a maximal spacelike surface in $\mathbf{AdS}^{2,1}$, which determines a minimal lagrangian diffeomorphism of $\mathbf{H}^2$, which in turn determines the initial homeomorphism.

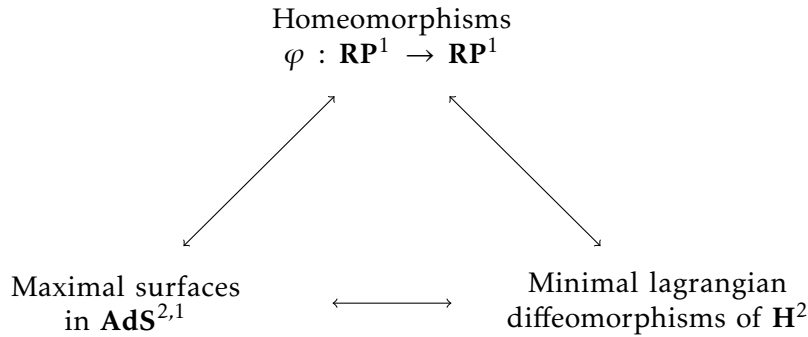

FIGURE 1. Homeomorphisms of $\mathbf{RP}^1$ identify with maximal spacelike surfaces in $\mathbf{AdS}^{2,1}$, which in turn correspond to minimal lagrangian diffeomorphisms of $\mathbf{H}^2$.

It turns out that the geometric properties of the homeomorphism, the minimal lagrangian diffeomorphism, and the maximal spacelike surface are closely tied to one another. Indeed, in [BS10] Bonsante–Schlenker showed that the homeomorphism is quasisymmetric if and only if the minimal lagrangian diffeomorphism is quasiconformal, if and only if the maximal spacelike surface has extrinsic curvature $K := \det(A)$ satisfying $\|K\|_{L^\infty} < 1$. In a similar manner, it is straightforward to show (see Theorem 5.2.2 of [Sep24]) that the homeomorphism is symmetric if and only if the minimal lagrangian diffeomorphism is asymptotically conformal, if and only if the maximal spacelike surface is asymptotically hyperbolic, and so on. Theorems A and B thus fit perfectly into this schema by showing that the homeomorphism is Weil–Petersson if and only if the minimal lagrangian diffeomorphism has $L^2$ Beltrami coefficient, if and only if the maximal spacelike surface has finite renormalized area. This is summarized in Table 1.

| Homeos. of $\mathbf{RP}^1$ | Maximal surfaces in $\mathbf{AdS}^{2,1}$ | Min. lag. diffeos. of $\mathbf{H}^2$ |
|---|---|---|
| Quasisymmetric | $\|K\|_{L^\infty} < 1$ | $\|\mu\|_{L^\infty} < 1$ |
| Symmetric | $\lim_{p \to \infty} \|K(p)\| = 0$ | $\lim_{p \to \infty} \|\mu(p)\| = 0$ |
| Weil–Petersson | $\|K\|_{L^2} < \infty$ | $\|\mu\|_{L^2} < \infty$ |

TABLE 1. Homeomorphisms, maximal surfaces, and minimal lagrangian diffeomorphisms. Properties of the maximal surfaces are expressed in terms of the extrinsic curvature $K := \det(A)$. Properties of the minimal lagrangian diffeomorphism are expressed in terms of the Beltrami differential $\mu$.

1.4. **Conformal welding and minimal surfaces in $\mathbf{H}^3$.** Theorem B may be viewed as a lorentzian analogue of the characterization given by Bishop in Definition 19 of [Bis25], given there in terms of Jordan curves in the Riemann sphere $\mathbf{CP}^1$. These two perspectives are related by the operation of *conformal welding*, which associates homeomorphisms of $\mathbf{RP}^1$ to oriented Jordan curves in $\mathbf{CP}^1$ as follows. Let $\Gamma$ be an oriented Jordan curve in $\mathbf{CP}^1$, and let $\Omega$ and $\Omega^*$ denote the components of its complement lying respectively on its left and on its right. Let $\mathbf{D}$ denote the Poincaré disk, denote $\mathbf{D}^* := \mathbf{CP}^1 \setminus \overline{\mathbf{D}}$, and let $f : \mathbf{D} \to \Omega$ and $f^* : \mathbf{D}^* \to \Omega^*$ be uniformizing maps. These maps respectively extend to homeomorphisms $f : \overline{\mathbf{D}} \to \overline{\Omega}$ and $f^* : \overline{\mathbf{D}}^* \to \overline{\Omega}^*$, and the *welding homeomorphism* of $\Gamma$ is defined by

$$\varphi_\Gamma := (f^*)^{-1} \circ f : \mathbf{RP}^1 \to \mathbf{RP}^1 \, . \tag{1.1}$$

This homeomorphism is well-defined up to post- and pre-composition by elements of $\mathbf{PSL}(2, \mathbf{R})$. The conformal welding problem asks for the reconstruction of the Jordan curve $\Gamma$ from a given homeomorphism $\varphi$ of $\mathbf{RP}^1$. This problem remains open in full generality, and it is not even known which homeomorphisms arise as conformal weldings (see, for example, [Bis]). However, it is known that conformal welding maps the space of quasicircles, quotiented by the action of $\mathbf{PSL}(2, \mathbf{C})$, homeomorphically onto the space of quasisymmetries, quotiented by the action of $\mathbf{PSL}(2, \mathbf{R}) \times \mathbf{PSL}(2, \mathbf{R})$.

A Jordan curve is *Weil–Petersson* whenever its conformal welding is a Weil–Petersson homeomorphism. In [And83] Anderson showed that that every Jordan curve in the Riemann sphere, viewed as the ideal boundary of $\mathbf{H}^3$, bounds a minimal surface in $\mathbf{H}^3$. It is natural to study the relationship between the regularity of the curve and the geometry of the minimal surface. In [Gra00], Graham–Witten introduced the concept of renormalized area for properly embedded minimal surfaces in Poincaré–Einstein spaces with sufficiently regular ideal boundary. Building on this work, Alexakis–Mazzeo proved in [AM10] that a $C^{3,\alpha}$ Jordan curve in $\mathbf{CP}^1$ bounds a minimal surface in $\mathbf{H}^3$ with finite renormalized area. The work [Bis25] of Bishop then completes this picture, by showing that a Jordan curve in $\mathbf{CP}^1$ is Weil–Petersson if and only if it bounds a minimal surface of finite geometry and of finite renormalized area.

However, properly embedded minimal surfaces in $\mathbf{H}^3$ can be problematic for many applications, principally due to their non-uniqueness (see, for example, [HW15, HLS26]). In particular, Bishop's characterization does not yield a bijection between two explicitly defined classes of objects. In contrast, uniqueness, and even stability, of complete maximal spacelike surfaces is always guaranteed in anti-de Sitter space. For reasons such as this, together with the identification of Mess discussed above, such surfaces are expected to present a far more useful class for potential applicatons.

Finally, we draw attention to the recent work [Wan25] of Wang, where conformal welding is interpreted as a correspondence between Jordan curves in $\mathbf{CP}^1$ and acausal closed curves in $\mathbf{AdS}^{2,1}$: every Jordan curve in $\mathbf{CP}^1$ corresponds to the graph in $\mathbf{Ein}^{1,1}$ of its welding. From the perspective of [Wan25], Bishop's work and Theorem B tell us that finiteness of the renormalized area is preserved by conformal welding. We refer the reader to [Wan25] for a conjectural picture of a general correspondence between geometric quantities in $\mathbf{H}^3$ and $\mathbf{AdS}^{2,1}$.

1.5. **Other related work.** Finally, we reviewing related results on similar classes of surfaces in Lorentzian space forms, describing the relationship between their geometry and the regularity of their ideal boundaries.

Closely related to the work [BS10] of Bonsante–Schlenker, in [BMS13], Bonsante–Mondello–Schlenker introduced *landslide diffeomorphisms*, which are a natural generalization of minimal Lagrangian diffeomorphisms, arising from constant mean curvature (CMC) surfaces in $\mathbf{AdS}^{2,1}$ via the Gauss map. A quantitative study of the relationship between the geometry of these CMC surfaces and the optimality of their associated landslides was recently carried out by the fifth author in [Tre25]. We expect the techniques developed in this paper to also yield to a characterization of CMC surfaces that bound Weil–Petersson circle homeomorphisms.

We also highlight related work on the infinitesimal geometry of Weil–Petersson Teichmüller space. In [Dia25], the first author establishes a correspondence between smooth spacelike surfaces in three-dimensional co-Minkowski space $\mathbf{HP}^3$ and divergence-free vector fields on the hyperbolic plane $\mathbf{H}^2$. We recall that *co-Minkowski space*, also known as *half-pipe space*, is the dual space to Minkowski space, having as projective model $\mathbf{H}^2 \times \mathbf{R}$. In this framework, continuous vector fields on the circle may be viewed as continuous real-valued functions on $\mathbf{S}^1$, and thus identify with closed curves in the ideal boundary $\mathbf{S}^1 \times \mathbf{R}$ of $\mathbf{HP}^3$. In [BF20], Barbot–Fillastre showed that every such curve bounds a unique *mean surface* in $\mathbf{HP}^3$, which is the analogue of a maximal surface in this context.

Finally, in [Dia25], the first author characterizes vector fields on the circle belonging to the Sobolev space $H^{3/2}(\mathbf{S}^1)$. This space is identified with the tangent space at the identity of the Weil–Petersson Teichmüller space. It is shown that a vector field on the circle lies in $H^{3/2}(\mathbf{S}^1)$ if and only if the renormalized area of the corresponding mean surface in $\mathbf{HP}^3$ is finite. This suggests that co-Minkowski geometry reflects the infinitesimal geometry of the Weil–Petersson Teichmüller space. This analogy between hyperbolic, anti–de Sitter, and co-Minkowski geometries can be viewed as an instance of the phenomenon of geometric transition, as described by Danciger in [Dan11].

1.6. **Proofs of main theorems and structure of paper.** Theorems A and B are proven through a chain of equivalences. In the process, we provide two further new, but more technical, characterizations of Weil–Petersson homeomorphisms which we also believe to be of independent interest. The five characterizations of the Weil–Petersson class used in this paper are summarized in Table 2, and the chain of equivalences used to prove Theorems A and B is summarized in Figure 2.

| Characterization | Definition/Result | Brief description |
| --- | --- | --- |
| 1 | Definition 1.1.1 | $\log(\varphi') \in H^{3/2}(\mathbf{S}^1)$ |
| 2 | Theorem 3.2.1 | Finite beta sum |
| 3 | Theorem 4.2.1 | Finite epsilon sum |
| 4 | Theorem B | Finite renormalized area |
| 5 | Theorem A | $L^2$ Beltrami differential |

TABLE 2. Characterizations of the Weil–Petersson class and their brief descriptions. Characterization (1) is given in [She18]. Characterizations (2), (3), (4), and (5) are new.

Our first new technical characterization is in terms of finiteness of what we call the beta sum. This is defined in Section 3.1, and may be viewed as the analogue in (1,1)-dimensional Einstein space of the characterization given by Bishop in Definition 11 of [Bis25]. The concept of beta numbers was first introduced by Jones in [Jon90] as part of his travelling salesman theorem. In the present context, *beta numbers* (see Definition 3.1.1) describe the *infinitesimal* structure of the homeomorphism. They are defined for each dyadic interval of $\mathbf{RP}^1$, and measure the deviation of the homeomorphism $\varphi$ from its $L^\infty$ best linear estimator over every such interval. As we will see presently (see Lemma 3.1.8), they provide a logarithmic measure of the rate of change between successive dyadic intervals of the mean gradient of $\varphi$ over these intervals.

The second is in terms of finiteness of what we call the epsilon sum. This is defined in Section 4.1, and may be viewed as the analogue in (1,1)-dimensional Einstein space of characterization given by Bishop in Definition 14 of [Bis25]. We may view *epsilon numbers* (see Definition 4.1.1) as a modification of beta numbers that take into account the *local*, as opposed to *infinitesimal*, structure of the homeomorphism. In this manner, they provide a bridge between our beta numbers of the homeomorphism $\varphi$ and the asymptotic structure of the complete maximal spacelike surface bounded by its graph.

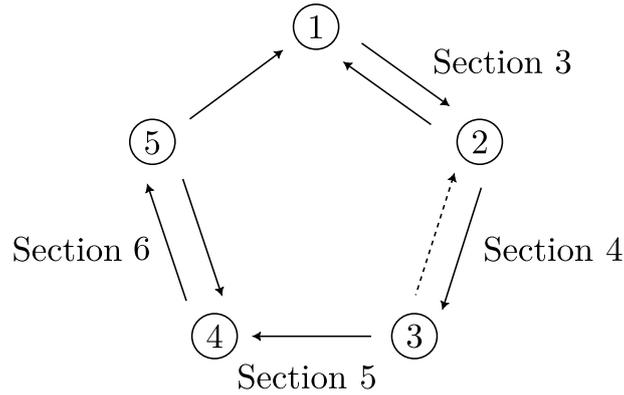

FIGURE 2. Characterizations of the Weil–Petersson class and the implications between them. The implication indicated by the dashed arrow requires the additional hypothesis of quasisymmetry, which is not required to prove the other implications. In particular, the graph can be followed in the clockwise direction, proving the equivalence of all five characterizations, without the need for this hypothesis.

The paper is structured as follows.

(1) Section 2 collects preliminary results on the geometry of $(2,1)$-dimensional anti-de Sitter space, 1-dimensional real projective space, and $(1,1)$-dimensional Einstein space.
(2) Section 3 introduces the beta numbers and the beta sums. In Theorem 3.2.1, we show that the Weil–Petersson property is equivalent to finiteness these beta sums.
(3) Section 4 introduces the epsilon numbers and the epsilon sums. In Theorem 4.2.1, we show that the finiteness of the beta sums is equivalent to quasisymmetry of the homeomorphism together with finiteness of the epsilon sums. Although the extra condition of quasisymmetry is required to directly recover finiteness of the beta sums, it is not actually necessary to complete the chain of equivalences, hence the dashed arrow in Figure 2.
(4) Section 5 addresses the renormalized area of maximal spacelike surfaces. In Theorem 5.1.1, we show that finiteness of any epsilon sum implies finite renormalized area of the complete maximal spacelike surface bounded by the graph of the homeomorphism. We note again that this step does not require us to suppose that the homeomorphism is quasisymmetric.
(5) Section 6 addresses the Beltrami differentials of minimal lagrangian diffeomorphisms. In Theorem 6.1.1 we show that the complete maximal spacelike surface bounded by the graph of the homeomorphism has finite renormalized area if and only if its minimal lagrangian extension has $L^2$ Beltrami differential. In addition, in Lemma 6.4.4, we prove the useful fact that if the minimal lagrangian extension $\Phi_\varphi$ of a homeomorphism $\varphi : \mathbf{RP}^1 \to \mathbf{RP}^1$ has $L^2$ Beltrami differential, then it is quasiconformal. Finally, bearing in mind the result [She18] of Shen, which states that $\varphi$ is Weil–Petersson whenever $\Phi_\varphi$ is quasiconformal and has $L^2$ Beltrami differential, this completes the chain of equivalences shown in Figure 2, and thus proves Theorems A and B.
(6) For the reader's convenience, in Appendix A, we provide a proof of the fact, proven by the fifth author in [Tre24b], that every homeomorphism $\varphi : \mathbf{RP}^1 \to \mathbf{RP}^1$ extends to a unique minimal lagrangian diffeomorphism of $\mathbf{H}^2$.

1.7. **Acknowledgements.** The authors are grateful to François Labourie, Andrea Seppi and Jérémy Toulisse for helpful conversations and the continued interest they have shown in this work.

A. Moriani thanks James Farre for sparking his interest in Bishop's work.

G. Smith's contribution was carried out in part whilst visiting the Dipartimento di Matematica "Giuseppe Peano" at Università di Torino, and he is grateful for the welcoming atmosphere and excellent working conditions provided during that stay.

E. Trebeschi thanks Francesco Bonsante, Timothé Lemistre and Simone Nati Poltri for several related discussion.

Part of this work was also carried out during meetings at the Laboratoire J. A. Dieudonné de l'Université Côte d'Azur, and the authors are grateful for the excellent working conditions provided during this stay.

The authors acknowledge funding by the European Research Council under the ERC-Advanced grant 101095722 AnSur (Geometric Analysis and Surface Groups) and the ERC-consolidator grant 101124349 GENERATE (GeomEtry and aNalysis for $(G, X)$–structurEs and their defoRmATion spacEs). Views and opinions expressed are however those of the authors only and do not necessarily reflect those of the European Union or the European Research Council Executive Agency. Neither the European Union nor the granting authority can be held responsible for them.

## 2. Preliminaries

The main spaces used in this paper will be $(2,1)$-dimensional anti-de Sitter space $\mathbf{AdS}^{2,1}$, 1-dimensional real projective space $\mathbf{RP}^1$, and $(1,1)$-dimensional Einstein space $\mathbf{Ein}^{1,1}$. It will be convenient to work with a number of different models and parametrizations of these spaces which we review in the present section.

2.1. **Anti-de Sitter space.**

2.1.1. *The projectivized quadric model.* For all $(p,q)$, let $\mathbf{R}^{p,q}$ denote the vector space $\mathbf{R}^{p,q}$ endowed with the quadratic form
$$\langle \mathbf{x}, \mathbf{y} \rangle_{p,q} := x_1 y_1 + \cdots + x_p y_p - x_{p+1} y_{p+1} - \cdots - x_{p+q} y_{p+q} \,.$$
We ommit the subscript $p,q$, and write $\langle \cdot, \cdot \rangle$, whenever it is unambiguous. Recall that a vector $\mathbf{x} \in \mathbf{R}^{p,q}$ is said to be *spacelike*, *timelike*, or *lightlike* according to whether its norm-squared is positive, negative, or null respectively.

We define $(2,1)$-dimensional *anti-de Sitter space* to be the projectivized quadric of vectors in $\mathbf{R}^{2,2}$ of norm-squared equal to $-1$, that is
$$\mathbf{AdS}^{2,1} := \{\mathbf{x} \in \mathbf{R}^{2,2} \mid \langle \mathbf{x}, \mathbf{x}\rangle = -1\}/\{\pm 1\} \,.$$
This is a lorentzian manifold of constant sectional curvature equal to $-1$. The quadric in $\mathbf{R}^{2,2}$ given by
$$\widehat{\mathbf{AdS}}^{2,1} := \{\mathbf{x} \in \mathbf{R}^{2,2} \mid \langle \mathbf{x}, \mathbf{x} \rangle = -1\}$$
is its double cover.

Recall that geodesics in any quadric are given by intersections of that quadric with planes passing through the origin. The causal nature of the geodesic corresponds to the signature of the plane. Hence, a geodesic in $\widehat{\mathbf{AdS}}^{2,1}$ is spacelike if it lies in a plane of signature $(1,1)$, timelike if it lies in a plane of signature $(0,2)$, and lightlight if it lies in a plane of signature $(0,1)$. Note that these are the only possible signatures for planes that meet $\widehat{\mathbf{AdS}}^{2,1}$ non-trivially.

Similarly, totally geodesic surfaces in any quadric are given by intersections with 3-dimensional vector subspaces, and the causal nature of the surface corresponds to the signature of the subspace: a surface is spacelike if it lies in a subspace of signature $(2,1)$, timelike if it lies in a subspace of signature $(1,2)$, and so on.

Consider the spacelike totally geodesic surface $H$ in $\mathbf{AdS}^{2,1}$ given by

$$(2.1) \qquad H := \{[\mathbf{x}] \in \mathbf{AdS}^{2,1} \mid x_4 = 0\} = \{(x_1, x_2, x_3) \in \mathbf{R}^3 \mid x_1^2 + x_2^2 - x_3^2 = -1\}/\{\pm 1\}.$$

This is just 2-dimensional hyperbolic space $\mathbf{H}^2$, and it is straightforward to see that every complete, spacelike totally geodesic surface in $\mathbf{AdS}^{2,1}$ is in fact isometric to $\mathbf{H}^2$. Similarly, every complete, timelike totally geodesic surface in $\mathbf{AdS}^{2,1}$ is isometric to $(1,1)$-dimensional anti-de Sitter space $\mathbf{AdS}^{1,1}$, which is itself isometric to $(1,1)$-dimensional de Sitter space $\mathrm{dS}^{1,1}$ with the sign of its metric reversed.

2.1.2. *The matrix model.* The study of $\mathbf{AdS}^{2,1}$ is facilitated by its natural identification with the Lie group $\mathbf{PSL}(2,\mathbf{R})$. We call this model the *matrix model* of $\mathbf{AdS}^{2,1}$. It is given as follows. Let $\mathcal{M}(2,\mathbf{R})$ denote the algebra of $2 \times 2$ real matrices, let $\mathbf{J} \in \mathcal{M}(2,\mathbf{R})$ denote the standard complex structure of $\mathbf{R}^2$, that is

$$(2.2) \qquad \mathbf{J} := \begin{pmatrix} 0 & -1 \\ 1 & 0 \end{pmatrix},$$

and denote

$$\langle \mathbf{M}, \mathbf{N} \rangle_{\mathrm{mat}} := \frac{1}{2} \mathrm{Tr}(\mathbf{M} \mathbf{J} \mathbf{N}^t \mathbf{J}).$$

This bilinear form is non-degenerate, symmetric, and of signature $(2,2)$. Furthermore, for all $\mathbf{M} \in \mathcal{M}(2,\mathbf{R})$,

$$\det(M) = -\langle \mathbf{M}, \mathbf{M} \rangle_{\mathrm{mat}}.$$

This yields the desired identification

$$\mathbf{AdS}^{2,1} = \{\mathbf{M} \in \mathcal{M}(2,\mathbf{R}) \mid \det(M) = 1\}/\{\pm 1\} = \mathbf{PSL}(2,\mathbf{R}).$$

More explicitly, we define the map $\mathcal{M} : \mathbf{R}^{2,2} \to \mathcal{M}(2,\mathbf{R})$ by

$$\mathcal{M}(\mathbf{x}) := \begin{pmatrix} x_3 - x_1 & x_2 - x_4 \\ x_2 + x_4 & x_3 + x_1 \end{pmatrix}.$$

For all $\mathbf{x}, \mathbf{y} \in \mathbf{R}^{2,2}$,

$$\langle \mathbf{x}, \mathbf{y} \rangle_{2,2} = \langle \mathcal{M}(\mathbf{x}), \mathcal{M}(\mathbf{y}) \rangle_{\mathrm{mat}}.$$

Thus $\mathcal{M}$ defines an isometry from $\mathbf{AdS}^{2,1}$ to $\mathbf{PSL}(2,\mathbf{R})$.

2.1.3. *Isometries.* Using the projective quadric model, we see that the isometry group $\mathsf{Isom}(\mathbf{AdS}^{2,1})$ of $\mathbf{AdS}^{2,1}$ naturally identifies with $\mathbf{PO}(2,2)$. This group consists of two connected components, the first made up of elements that preserve the orientation and time orientation, the second made of elements that reverse these orientations. We denote its identity component by $\mathsf{Isom}_0(\mathbf{AdS}^{2,1})$.

The matrix model provides a natural identification of $\mathsf{Isom}_0(\mathbf{AdS}^{2,1})$ with $\mathbf{PSL}(2,\mathbf{R}) \times \mathbf{PSL}(2,\mathbf{R})$. Indeed, the action of $\mathbf{PSL}(2,\mathbf{R}) \times \mathbf{PSL}(2,\mathbf{R})$ on $\mathcal{M}(2,\mathbf{R})$ is given by

$$(2.3) \qquad \alpha(\mathbf{M}, \mathbf{N}) \cdot \mathbf{A} := \mathbf{M} \cdot \mathbf{A} \cdot \mathbf{N}^{-1}.$$

This action is faithful, preserves the determinant, and descends to an action on the projectivization of this space. It therefore preserves $\mathbf{PSL}(2,\mathbf{R})$ and acts isometrically on this space, yields the desired identification of $\mathbf{PSL}(2,\mathbf{R}) \times \mathbf{PSL}(2,\mathbf{R})$ with $\mathsf{Isom}_0(\mathbf{AdS}^{2,1})$.

2.1.4. *Kleinian charts.* Kleinian charts are most useful for studying the geodesic and asymptotic structure of $\mathbf{AdS}^{2,1}$. Throughout the sequel, with $H$ as in (2.1), above, we will work with the *kleinian chart* of $\mathbf{AdS}^{2,1}$ defined by

$$\mathcal{K} : \mathbf{AdS}^{2,1} \setminus H \to \mathbf{R}^{2,1}; [\mathbf{x}] \mapsto \frac{1}{x_4}(x_1, x_2, x_3).$$

The image of this chart is the open set

$$\Omega := \left\{ \mathbf{y} \in \mathbf{R}^{2,1} \mid \langle \mathbf{y}, \mathbf{y} \rangle < 1 \right\}.$$

The kleinian chart maps geodesics in $\mathbf{AdS}^{2,1}$ into straight lines in $\mathbf{R}^{2,1}$. Furthermore, its image is the set bounded by the quadric

$$\widehat{\mathbf{dS}}^{1,1} := \partial\Omega = \left\{ \mathbf{y} \in \mathbf{R}^{2,1} \mid \langle \mathbf{y}, \mathbf{y} \rangle = 1 \right\},$$

which itself identifies with the double cover of $(1,1)$-dimensional de Sitter space. Spacelike geodesics are sent to lines which intersect $\widehat{\mathbf{dS}}^{1,1}$ in two distinct points. Timelike geodesics are sent to lines that neither intersect $\widehat{\mathbf{dS}}^{1,1}$, nor are tangent to this quadric. Lightlike geodesics are sent to lines which either are tangent to $\widehat{\mathbf{dS}}^{1,1}$ or are asymptotic to this quadric. This is illustrated in Figure 3.

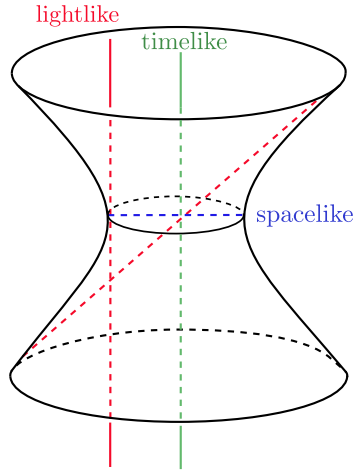

Figure 3. Geodesics of $\mathbf{AdS}^{2,1}$ in the kleinian chart.

## 2.2. Real projective space.

2.2.1. *The quotient model.* We denote by $\mathbf{RP}^1$ the one-dimensional real projective space, that is, the space of lines in $\mathbf{R}^2$ passing by the origin. The automorphism group of $\mathbf{RP}^1$ is the projective special linear group $\mathbf{PSL}(2,\mathbf{R})$. In particular, $\mathbf{RP}^1$ does not carry a canonical riemannian metric. Throughout the rest of this paper, we will make an explicit choice of metric over $\mathbf{RP}^1$. The highbrow approach would be to select a complex structure over $\mathbf{R}^2$ and note that its stabilizer in $\mathbf{PSL}(2,\mathbf{R})$ is a copy of $\mathsf{SO}(2,\mathbf{R})$.

We adopt the following lowbrow approach. Here and in the rest of the paper, we furnish $\mathbf{R}^2$ with its canonical metric and we denote by $e_1, e_2$ its canonical orthonormal basis. We define $\mathcal{P}: \mathbf{R} \to \mathbf{RP}^1$ by

$$\mathcal{P}(x) := [\cos(x) : \sin(x)].$$

This function is a universal covering map of $\mathbf{RP}^1$, and thus projects to a diffeomorphism from $\mathbf{R}/\pi\mathbf{Z}$ to $\mathbf{RP}^1$, which we also denote by $\mathcal{P}$. In this manner, we identify

$$\mathbf{RP}^1 \simeq \mathbf{R}/\pi\mathbf{Z}.$$

We call this the *quotient model* of $\mathbf{RP}^1$. Using this model, we furnish $\mathbf{RP}^1$ with the riemannian metric that it inherits from $\mathbf{R}$. We emphasize the non-canonicity of this choice because, although our results are projectively invariant, the constructions used to prove them are not.

Finally, given an interval $I \subseteq \mathbf{RP}^1$, we say that a function $f : I \to \mathbf{RP}^1$ is *linear* whenever it has the form

$$f(x) = ax + b,$$

in this model.

2.2.2. *The affine chart.* *Affine charts* of $\mathbf{RP}^1$ are most useful for studying elements of $\mathbf{PSL}(2,\mathbf{R})$. Throughout the sequel, we work in the affine chart defined by

$$\mathcal{A}([x_1 : x_2]) = \frac{x_2}{x_1} \, .$$

Its domain is the complement of $[0 : 1]$ in $\mathbf{RP}^1$ and its image is the whole of $\mathbf{R}$. We thus view $[0 : 1]$ as the point at infinity. Note that, composing with the function $\mathcal{P}$ defined above yields, over the fundamental domain $]-\pi/2, \pi/2[$,

$$(\mathcal{A} \circ \mathcal{P})(x) = \tan(x) \, .$$

In the affine chart, elements of $\mathbf{PSL}(2,\mathbf{R})$ have a particularly simple structure. Indeed, they act by fractional linear maps, that is, functions of the form

$$f(t) = \frac{at+b}{ct+d} \, ,$$

where $a$, $b$, $c$, and $d$ are real numbers satisfying

$$ad - bc = 1 \, .$$

In what follows, it will be useful to view this family instead as the union of the set of linear functions $f : \mathbf{R} \to \mathbf{R}$ of the form

$$f(t) = At + B \, ,$$

where $A$ and $B$ are real numbers and $A$ is positive, with the set of rational functions $f : \mathbf{R} \to \mathbf{R}$ of the form

$$(2.4) \qquad f(t) = \frac{P}{Q-t} - R \, ,$$

where $P$, $Q$ and $R$ are real numbers and $P$ is positive. The graph of any such rational function is a hyperbola with horizontal and vertical asymptotes given respectively by

$$H := \{y = -R\} \qquad \text{and} \qquad V := \{x = Q\} \, .$$

We call the point of intersection of these two asymptotes $(Q, -R)$ the *center* of the hyperbola. This structure is illustrated in Figure 7, and we will make considerable use of it throughout the sequel.

2.3. **Einstein space.**

2.3.1. *The projectivized quadric model.* We define $(1,1)$-dimensional *Einstein space* to be the projectivized cone of null vectors in $\mathbf{R}^{2,2}$, that is

$$\mathbf{Ein}^{1,1} := \{\mathbf{x} \in \mathbf{R}^{2,2} \setminus \{0\} \mid \langle \mathbf{x}, \mathbf{x} \rangle = 0\}/\mathbf{R}^* \, .$$

The metric of $\mathbf{R}^{2,2}$ projects to a conformal class of non-degenerate metrics over $\mathbf{Ein}^{1,1}$ of signature $(1,1)$ making this space into a conformal lorentzian manifold. Viewed as a subset of $\mathbf{RP}^3$, $\mathbf{Ein}^{1,1}$ coincides with the topological boundary of $\mathbf{AdS}^{2,1}$, and is thus identified with the ideal boundary $\partial_\infty \mathbf{AdS}^{2,1}$ of the latter.

2.3.2. *The matrix-angle model.* The matrix-angle model is most useful for studying the global structure of $\mathbf{Ein}^{1,1}$ and the action of its automorphism group $\mathbf{PSL}(2,\mathbf{R}) \times \mathbf{PSL}(2,\mathbf{R})$. In order to describe this model, we first review the geometry of $\mathbf{Ein}^{1,1}$ in the matrix model given in Section 2.1.2, above. In this model, Einstein space coincides with the projectivization of the space of rank 1 matrices in $\mathcal{M}$, that is

$$\mathbf{Ein}^{1,1} = \{\mathbf{M} \in \mathcal{M} \mid \mathrm{rank}(\mathbf{M}) = 1\}/\mathbf{R}^* \, .$$

Projective classes of rank 1 matrices over $\mathbf{R}^2$ are determined by their kernel and their image. Using the canonical basis of $\mathbf{R}^2$ given in Section 2.2.1, we thus define $\mathcal{N} : \mathbf{R}^2 \to \mathcal{M}$ by

$$\mathcal{N}(x_1, x_2) := \begin{pmatrix} \cos(x_1) \\ \sin(x_1) \end{pmatrix} \otimes \big(\cos(x_2), \sin(x_2)\big) = \begin{pmatrix} \cos(x_1)\cos(x_2) & \cos(x_1)\sin(x_2) \\ \sin(x_1)\cos(x_2) & \sin(x_1)\sin(x_2) \end{pmatrix} \, .$$

We denote the image of $\mathcal{N}$ by $\widetilde{\mathbf{Ein}}^{1,1}$. It is a smooth, embedded surface in $\mathcal{M}(2,\mathbf{R})$ projectivizing to a double cover of Einstein space. The map $\mathcal{N}$ is doubly-periodic with period lattice generated by $(\pi,\pi)$ and $(-\pi,\pi)$. A fundamental domain is $]0,\pi[\times]0,2\pi[$, illustrated in Figure 4. The left and right edges are identified by the period vector $(\pi,\pi)$, so that matching these two edges involves a vertical translation by $\pi$. The upper and lower edges are identified by the period vector $(0,2\pi)$, so that no horizontal translation is involved in matching these edges.

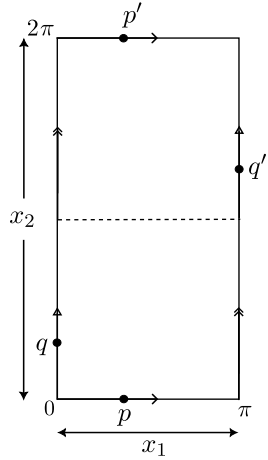

FIGURE 4. A fundamental domain of the double cover of $\mathbf{Ein}^{1,1}$.

Pulling back $\langle\cdot,\cdot\rangle_{\mathrm{mat}}$ through $\mathcal{N}$ yields

(2.5) $$\mathcal{N}^*\langle\cdot,\cdot\rangle_{\mathrm{mat}} = \frac{1}{2}\bigl(d\alpha_1 d\alpha_2 + d\alpha_2 d\alpha_1\bigr).$$

This is a non-degenerate metric of signature $(1,1)$ with causal structure as in Figure [???]. At all points, horizontal and vertical vectors are null, hence lightlike; vectors in the upper-right and lower-left quadrants are positive, hence spacelike; and vectors in the upper-left and lower-right quadrants are negative, hence timelike. In particular, the lightlike geodesics of this metric are the horizontal and vertical lines.

Projectivizing $\mathcal{N}$ yields a map, which we also denote by $\mathcal{N}$, from $\mathbf{R}^2$ into $\mathbf{Ein}^{1,1}$. This map is doubly-periodic with period lattice generated by $(\pi,0)$ and $(0,\pi)$. A fundamental domain is $]0,\pi[\times]0,\pi[$, illustrated in Figure 5. Bearing in mind the identification of $\mathbf{RP}^1$ with $\mathbf{R}/\pi\mathbf{Z}$ given in Section 2.2.1, $\mathcal{N}$ yields a diffeomorphism from $\mathbf{RP}^1 \times \mathbf{RP}^1$ to $\mathbf{Ein}^{1,1}$ which we also denote by $\mathcal{N}$. We call this diffeomorphism the *matrix angle chart* of $\mathbf{Ein}^{1,1}$.

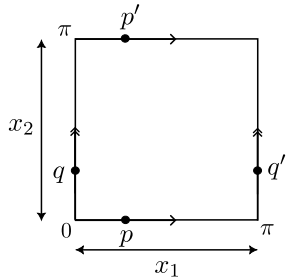

FIGURE 5. A fundamental domain of $\mathcal{N}: \mathbf{R}^2 \to \mathbf{Ein}^{1,1}$.

For all $(L_1,L_2) \in \mathbf{RP}^1 \times \mathbf{RP}^1$, $\mathcal{N}(L_1,L_2)$ is the projectivization of the matrix with image $L_1$ and kernel $L_2^\perp$. It follows from (2.3) that the natural action of $\mathbf{PSL}(2,\mathbf{R}) \times \mathbf{PSL}(2,\mathbf{R})$ on $\mathbf{Ein}^{1,1}$ is given in the matrix model by

$$\tilde{\beta}(\mathbf{M},\mathbf{N}) := (\mathbf{M}\cdot L_1, \mathbf{N}^t \cdot L_2).$$

In particular, this action preserves the set of all curves $\Gamma \in \mathbf{Ein}^{1,1}$ of the form
$$\Gamma = \{(x, f(x)) \mid x \in \mathbf{RP}^1\},$$
for some $f \in \mathrm{PSL}(2, \mathbf{R})$. The set of all such curves, which we call *acausal circles*, is fundamental to the study of the geometry of $\mathbf{Ein}^{1,1}$.

2.3.3. *The Penrose chart.* The Penrose chart is most useful for studying acausal circles in $\mathbf{Ein}^{1,1}$. We first work in the projectivized quadric model of Section 2.3.1. Consider the null subset
$$N := \{[\mathbf{x}] \in \mathbf{Ein}^{1,1} \mid x_3 - x_1 = 0\}.$$
This set is the light cone in $\mathbf{Ein}^{1,1}$ of the point $[1:0:-1:0]$, and consists of two null geodesics. The *Penrose chart* is defined over the complement of $N$ in $\mathbf{Ein}^{1,1}$ by
$$\mathcal{P}^{\mathrm{pr}} : \mathbf{Ein}^{1,1} \setminus N \to \mathbf{R}^2; [\mathbf{x}] \mapsto \frac{1}{x_3 - x_1}(x_2, -x_4).$$
This chart provides a conformal parametrization of the complement of $N$ in $\mathbf{Ein}^{1,1}$ by $(1,1)$-dimensional Minkowski space $\mathbf{R}^{1,1}$, as we will see presently.

The Penrose chart has a simple expression in the matrix angle model. Indeed, with $\mathcal{M}$ and $\mathcal{N}$ as in Sections 2.1.2 and 2.3.2 respectively, we obtain
$$\mathcal{P}^{\mathrm{mat}}(x_1, x_2) := (\mathcal{P}^{\mathrm{pr}} \circ \mathcal{M}^{-1} \circ \mathcal{N})(x_1, x_2) = \frac{1}{2}\big(\tan(x_1) + \tan(x_2), -\tan(x_1) + \tan(x_2)\big).$$
Upon composing with the rotation-scaling
$$\mathbf{R} := \begin{pmatrix} 1 & -1 \\ 1 & 1 \end{pmatrix},$$
we obtain the *rotated Penrose chart*
$$\mathcal{P}^{\mathrm{rot}}(x_1, x_2) := \big(\tan(x_1), \tan(x_2)\big).$$

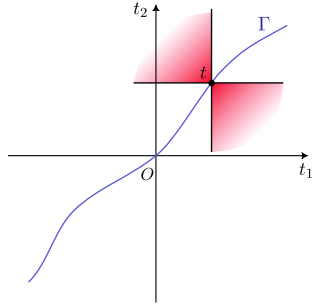

FIGURE 6. The causal structure of the rotated Penrose chart. In particular, an acausal curve is the graph of a monotone increasing function.

Bearing in mind Section 2.2.2, in the matrix angle model $\mathbf{RP}^1 \times \mathbf{RP}^1$ of $\mathbf{Ein}^{1,1}$, the rotated Penrose chart is simply obtained upon taking the affine charts separately in each of the two components. In particular, by (2.5), the metric of $\mathbf{Ein}^{1,1}$ in this chart is conformally equivalent to the metric $(dt_1 dt_2 + dt_2 dt_1)$, which is just the standard metric of $(1,1)$-dimensional Minkowski space rotated through an angle of $\pi/4$. The causal structure of $\mathbf{Ein}^{1,1}$ in the rotated Penrose chart is thus the same as that of the matrix angle model described in the preceding section and illustrated in Figure 6. Thus, as before, horizontal and vertical vectors are null, hence lightlike; vectors in the upper-right and lower-left quadrants are positive, hence spacelike; and vectors in the upper-left and lower-right quadrants are negative, hence timelike. Furthermore, in this chart, the set of acausal circles is the union of the set of straight lines with positive gradient and the set of euclidean hyperbolae with positive gradient and horizontal and vertical asymptotes, as in Figure 7.

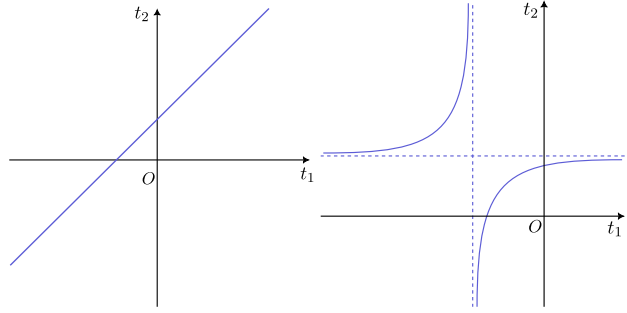

Figure 7. Two acausal circles in the rotated Penrose chart.

## 2.4. Dyadic decompositions.

2.4.1. *Dyadic decompositions.* Throughout the sequel, we will make considerable use of dyadic decompositions of $\mathbf{RP}^1$. These are defined as follows. We identify $\mathbf{RP}^1 := \mathbf{R}/\pi\mathbf{Z}$ as in Section 2.2.1. We choose a base point $x_0 \in \mathbf{RP}^1$ and, for every integer $m \geqslant 0$, we partition $\mathbf{RP}^1$ into $2^m$ arcs of equal length

$$I_{m,k} := x_0 + \left[\frac{k\pi}{2^m}, \frac{(k+1)\pi}{2^m}\right] \bmod \pi, \quad k = 0,\ldots, 2^m - 1.$$

For all such $m$, we call the set

$$\mathcal{D}_{x_0,m} := \{I_{m,k} \mid k = 0,\ldots, 2^m - 1\}$$

the $x_0$-*based dyadic partition at depth m* of $\mathbf{RP}^1$. The $x_0$-*based dyadic decomposition* of $\mathbf{RP}^1$ is then defined to be the union of all dyadic partitions over all depths, that is

$$\mathcal{D}_{x_0} := \bigcup_{m \geqslant 0} \mathcal{D}_{x_0,m}.$$

It will also be useful to denote, for all $k \geqslant 0$,

$$\mathcal{D}_{x_0,m \geqslant k} := \bigcup_{m \geqslant k} \mathcal{D}_{x_0,m}.$$

Note that any two dyadic decompositions differ by a rotation. In what follows, the base point will usually be taken to be equal to zero, and will not be mentioned explicitly.

*Remark* 2.4.1. We underline that this construction depends on a choice of riemannian metric over $\mathbf{RP}^1$, itself depending on a choice of basis of $\mathbf{R}^2$. The families of all dyadic partitions and all dyadic decompositions are thus rotationally invariant, but not projectively invariant. Although this may seem a surprising choice to make, as the theory of Weil–Petersson homeomorphisms is itself projective, it parallels the techniques developed by Bishop in [Bis25].

We call any element of a dyadic decomposition a *dyadic interval*. We denote the length of any interval $I$ by $\ell(I)$. We say that two dyadic intervals $I$ and $J$ are *successive* whenever $I$ is contained in $J$ and $\ell(J) = 2\ell(I)$.

Given any dyadic interval $I$, and any real number $r \geqslant 1$, we define $rI$ to be the interval centred on $I$ of length equal to $r$ times $\ell(I)$. In the sequel, we will typically work with triples of dyadic intervals instead of the dyadic intervals themselves, for two related reasons. Firstly, a dyadic point of the form $k \cdot \pi \cdot 2^{-m}$ only lies in the interior of dyadic intervals of length at least $\pi \cdot 2^{-(m-1)}$, so that dyadic intervals in themselves are insufficient for studying infinitesimal properties near any such point. Secondly, given two dyadic decompositions $\mathcal{D}$ and $\mathcal{D}'$, there is in general no relationship between elements of $\mathcal{D}$ and $\mathcal{D}'$. By contrast, for all $I \in \mathcal{D}$ at depth $m$, there is a dyadic interval $J$ in $\mathcal{D}'$ at depth $(m-1)$ such that $3I \subseteq 3J$. In this manner, triples allows us to directly compare sums over $\mathcal{D}$ with sums over $\mathcal{D}'$.

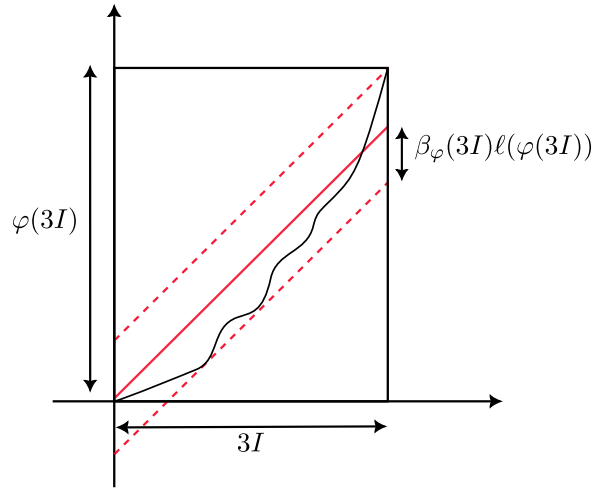

FIGURE 8. Definition of $\beta_\varphi(3I)$ and the $L^\infty$ best linear estimator of $\varphi$ over $3I$.

## 3. BETA SUMS

Our first new characterization of the Weil–Petersson class will be in terms of finiteness of what we will call the beta sum. This is the analogue on the present framework of the quantity studied in Section 4 of [Bis25]. Using techinical results proven by Bishop in [Bis22], we will show that a homeomorphism $\varphi : \mathbf{RP}^1 \to \mathbf{RP}^1$ is Weil–Petersson if and only if it has finite beta sum.

### 3.1. Beta numbers.
We will describe the beta sum characterization using the quotient model $\mathbf{RP}^1 = \mathbf{R}/\pi\mathbf{Z}$ given in Section 2.2.1. Let $\mathcal{D}$ be a dyadic decomposition of $\mathbf{RP}^1$, as in Section 2.4.1.

**Definition 3.1.1.** Given a homeomorphism $\varphi : \mathbf{RP}^1 \to \mathbf{RP}^1$, and a dyadic interval $I$, we define

$$\beta_\varphi(I) := \frac{1}{\ell(\varphi(I))} \inf_a \|\varphi - a\|_{L^\infty(I)}, \tag{3.1}$$

where the infimum is taken over all linear functions $a : I \to \mathbf{RP}^1$.

*Remark* 3.1.2. By taking $a$ to be the constant function sending $I$ to the mid-point of $\varphi(I)$, we see that

$$\beta_\varphi(I) \leq \frac{1}{2}.$$

In fact, the reader may verify that, by strict monotonicity of $\varphi$, this inequality is strict for all $I$.

The quantity $\beta_\varphi(I)$ should be understood as a scale-invariant measure of the distance of $\phi$ to some $L^\infty$ best linear estimator over $I$, analogous to the $L^2$ best linear estimator typically used in linear regression (see Figure 8). Although not strictly necessary for what follows, it is useful to note that this $L^\infty$ best linear estimator always exists and is unique.

**Lemma 3.1.3.** *Let $I$ be a compact subinterval of $\mathbf{R}$. Every continuous function $\varphi : I \to \mathbf{R}$ has a unique $L^\infty$ best linear estimator.*

*Remark* 3.1.4. Uniqueness is not wholly trivial on account of the non-strict convexity of the $L^\infty$ norm, and even fails without the hypothesis of continuity. Indeed, the reader may verify that the Heaviside step function does not have a unique $L^\infty$ best linear estimator over the interval $[-1, 1]$.

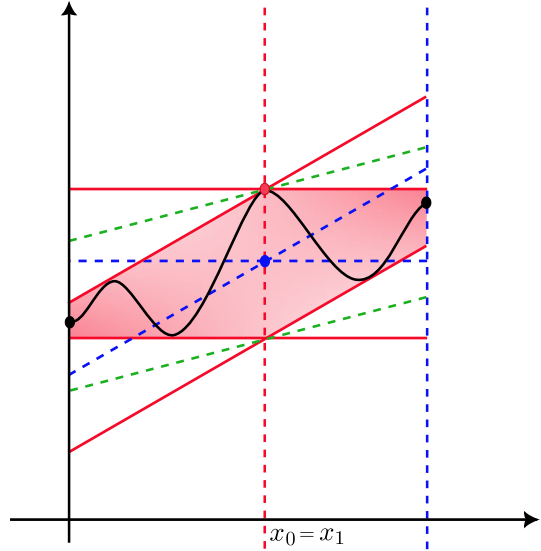

Figure 9. Uniqueness of the $L^\infty$ best linear estimator.

*Proof.* Consider the function
$$E(a,b) := \|f - (\gamma x + \delta)\|_{L^\infty(I)}\,.$$
Since this function is continuous and tends to infinity as $(\gamma,\delta)$ diverges, it attains a minimum value $E_0$, say, at some point $(\gamma_1,\delta_1)$, say. This proves existence.

We now prove uniqueness. Suppose the contrary, so that $E$ also attains its minimum at some other point $(\gamma_2,\delta_2)$, say. Note that we may suppose that $E_0 > 0$, for otherwise $f$ is linear and uniqueness follows trivially. By definition of $E$, for each $i$,
$$(\gamma_i x + \delta_i) - E_0 \leqslant f \leqslant (\gamma_i x + \delta_i) + E_0\,.$$
Consequently,

(3.2) $$\max(\gamma_1 x + \delta_1, \gamma_2 x + \delta_2) - E_0 \leqslant f \leqslant \min(\gamma_1 x + \delta_1, \gamma_2 x + \delta_2) + E_0\,.$$

Hence, denoting $\gamma := (\gamma_1 + \gamma_2)/2$ and $\delta := (\delta_1 + \delta_2)/2$,

(3.3) $$(\gamma x + \delta) - E_0 \leqslant f \leqslant (\gamma x + \delta) + E_0\,.$$

We now claim that one of the inequalities in (3.3) is strict. Indeed, otherwise, by continuity, there exist points $x_1, x_2 \in I$ such that
$$f(x_1) = (\gamma x_1 + \delta) - E_0 \qquad \text{and} \qquad f(x_2) = (\gamma x_2 + \delta) + E_0\,.$$
However, by (3.2) again, these equalities can only hold at the unique point of intersection of the graphs of $(\gamma_1 x + \delta_1)$ and $(\gamma_2 x + \delta_2)$, as in Figure 9. It follows that $x_1 = x_2$, and so $E_0 = 0$. This is absurd, and strictness of one of the inequalities in (3.3) follows. Without loss of generality, we may suppose that the first inequality in (3.3) is strict. Then, by continuity again, for some $\epsilon > 0$,
$$(\gamma x + (\delta + \epsilon)) - (E_0 - \epsilon) = (\gamma x + \delta) - (E_0 - 2\epsilon) \leqslant f \leqslant (\gamma x + \delta) + E_0 = (\gamma x + (\delta + \epsilon)) + (E_0 - \epsilon)\,,$$
so that
$$E(\gamma, \delta + \epsilon) \leqslant E_0 - \epsilon < E_0\,.$$
This is absurd, and uniqueness follows. $\square$

We now establish various elementary properties of $\beta_\varphi(I)$ that will be required in the sequel. In what follows, for any interval $I$, we denote the $L^\infty$ best linear estimator of $\varphi$ over $I$ by $a_I$ and we denote its gradient by $\gamma(I)$.

**Lemma 3.1.5.** *For every interval $I$,*

$$(3.4) \qquad \left(1 - 2\beta_\varphi(I)\right) \cdot \frac{\ell(\varphi(I))}{\ell(I)} \leqslant \gamma(I) \leqslant \left(1 + 2\beta_\varphi(I)\right) \cdot \frac{\ell(\varphi(I))}{\ell(I)} \,.$$

*In particular, for $\beta_\varphi(I) < 1/2$, $\gamma(I) > 0$.*

*Remark* 3.1.6. Thus, when $\beta_\varphi(I)$ is small, $\gamma(I)$ is equivalent to the mean gradient of $\varphi$ over the interval $I$.

*Remark* 3.1.7. In particular, by Remark 3.1.2, for all $I$, $\gamma(I) > 0$.

*Proof.* Let $x_-$ and $x_+$ denote respectively the lower and upper end of $I$. By definition

$$|\varphi(x_\pm) - a_I(x_\pm)| \leqslant \beta_\varphi(I) \cdot \ell(\varphi(I)) \,.$$

Thus,

$$\begin{aligned}
\left| \gamma(I) - \frac{\ell(\varphi(I))}{\ell(I)} \right| &= \left| \frac{1}{\ell(I)} \big(a_I(x_+) - a_I(x_-)\big) - \frac{1}{\ell(I)} \big(\varphi(x_+) - \varphi(x_-)\big) \right| \\
&\leqslant \frac{1}{\ell(I)} \big(|a_I(x_+) - \varphi(x_+)| + |a_I(x_-) - \varphi(x_-)|\big) \\
&\leqslant \frac{2\beta_\varphi(I) \cdot \ell(\varphi(I))}{\ell(I)} \,,
\end{aligned}$$

and the result follows. $\square$

**Lemma 3.1.8.** *Choose $\lambda > 1$, and let $I \subseteq I'$ be intervals such that*

$$\ell(I') \leqslant \lambda \ell(I) \,.$$

*If $\beta_\varphi(I') \leqslant 1/4$, then*

$$(3.5) \qquad 1 - 8\lambda \cdot \beta_\varphi(I') \leqslant \frac{\gamma(I)}{\gamma(I')} \leqslant 1 + 8\lambda \cdot \beta_\varphi(I') \,.$$

*Remark* 3.1.9. Thus, when $\beta_\varphi(I)$ is small, this quantity may be viewed as the rate of change of the mean gradient of $\varphi$ on a logarithmic scale. This simple observation will prove fundamental in what follows.

*Proof.* Let $x_-$ and $x_+$ denote respectively the lower and upper ends of $I$. By definition,

$$|\varphi(x_\pm) - a_{I'}(x_\pm)| \leqslant \|\varphi - a_{I'}\|_{L^\infty(I')} = \beta_\varphi(I') \cdot \ell(\varphi(I')) \,.$$

Likewise, bearing in mind that $a_I$ is the $L^\infty$ best estimator of $\varphi$ over $I$,

$$|\varphi(x_\pm) - a_I(x_\pm)| \leqslant \|\varphi - a_I\|_{L^\infty(I)} \leqslant \|\varphi - a_{I'}\|_{L^\infty(I)} \leqslant \|\varphi - a_{I'}\|_{L^\infty(I')} = \beta_\varphi(I') \cdot \ell(\varphi(I')) \,.$$

Thus,

$$\begin{aligned}
|\gamma_I - \gamma_{I'}| &= \left| \frac{1}{\ell(I)} \big(a_I(x_+) - a_I(x_-)\big) - \frac{1}{\ell(I)} \big(a_{I'}(x_+) - a_{I'}(x_-)\big) \right| \\
&\leqslant \frac{1}{\ell(I)} \big(|a_I(x_+) - \varphi(x_+)| + |a_{I'}(x_+) - \varphi(x_+)| \\
&\qquad\qquad + |a_I(x_-) - \varphi(x_-)| + |a_{I'}(x_-) - \varphi(x_-)|\big) \\
&\leqslant \frac{4\beta_\varphi(I') \cdot \ell(\varphi(I'))}{\ell(I)} \\
&\leq \frac{4\lambda \cdot \beta_\varphi(I') \cdot \ell(\varphi(I'))}{\ell(I')} \,.
\end{aligned}$$

Dividing by $\gamma(I')$ then yields

$$\left| \frac{\gamma(I)}{\gamma(I')} - 1 \right| \leqslant \frac{4\lambda \cdot \beta_\varphi(I') \cdot \ell(\varphi(I'))}{\ell(I') \cdot \gamma(I')} \,,$$

so that, by (3.4),
$$\left|\frac{\gamma(I)}{\gamma(I')} - 1\right| \leq \frac{4\lambda \cdot \beta_\varphi(I')}{(1 - 2\beta_\varphi(I'))}.$$

Thus, since $\beta_\varphi(I') \leq 1/4$,
$$\left|\frac{\gamma(I)}{\gamma(I')} - 1\right| \leq 8\lambda \cdot \beta_\varphi(I'),$$

as desired. □

**Lemma 3.1.10.** *Choose $\lambda > 1$, and let $I \subseteq I'$ be intervals such that*
$$\ell(I') \leq \lambda \ell(I).$$

*If $\beta_\varphi(I') \leq 1/16\lambda$, then*
$$\ell(\varphi(I')) \leq 4\lambda \cdot \ell(\varphi(I)).$$

*Remark* 3.1.11. In other words, when $\beta_\varphi(I')$ is small, the function $\varphi$ is almost linear over scales comparable to the length of $I'$. In particular, if $I$ and $I'$ have equivalent lengths, then so too do $\varphi(I)$ and $\varphi(I')$.

*Proof.* Indeed, by (3.4),
$$\frac{\ell(\varphi(I'))}{\ell(I')} \leq \frac{\gamma(I')}{1 - 2\beta_\varphi(I')} \quad \text{and} \quad \frac{\ell(\varphi(I))}{\ell(I)} \geq \frac{\gamma(I)}{1 + 2\beta_\varphi(I)}.$$

Thus, bearing in mind (3.5) and the hypotheses on $\beta_\varphi(I')$,
$$\frac{\ell(\varphi(I'))}{\ell(\varphi(I))} \leq \frac{\gamma(I') \cdot (1 + 2\beta_\varphi(I)) \cdot \ell(I')}{\gamma(I) \cdot (1 - 2\beta_\varphi(I')) \cdot \ell(I)} \leq \frac{\lambda \cdot (1 + 2\beta_\varphi(I))}{(1 - 2\beta_\varphi(I')) \cdot (1 - 8\lambda\beta_\varphi(I'))} \leq \frac{8\lambda \cdot (1 + 2\beta_\varphi(I))}{3}.$$

Finally, since $a_I$ is the $L^\infty$ best linear estimator of $\varphi$ over $I$,
$$\beta_\varphi(I) \cdot \ell(\varphi(I)) = \|\varphi - a_I\|_{L^\infty(I)} \leq \|\varphi - a_{I'}\|_{L^\infty(I)} \leq \|\varphi - a_{I'}\|_{L^\infty(I')} = \beta_\varphi(I') \cdot \ell(\varphi(I')),$$

so that
$$\beta_\varphi(I) \leq \frac{\ell(\varphi(I'))}{\ell(\varphi(I))} \cdot \beta_\varphi(I') \leq \frac{\ell(\varphi(I'))}{16\lambda \cdot \ell(\varphi(I))}.$$

Combining this with the preceding identity yields
$$\frac{\ell(\varphi(I'))}{\ell(\varphi(I))} \leq \frac{8\lambda}{3} + \frac{\ell(\varphi(I'))}{3\ell(\varphi(I))},$$

so that
$$\frac{\ell(\varphi(I'))}{\ell(\varphi(I))} \leq 4\lambda,$$

and the result follows. □

### 3.2. The beta sum characterization of the Weil–Petersson class.
The main result of this section is the following characterization of Weil–Petersson homeomorphisms in terms of beta numbers.

**Theorem 3.2.1.** *Choose $\lambda \geq 1$, and let $\varphi \colon \mathbf{RP}^1 \to \mathbf{RP}^1$ be a homeomorphism. The following are equivalent.*
  (1) *$\varphi$ is Weil–Petersson.*
  (2) *For every dyadic decomposition $\mathcal{D}$ of $\mathbf{RP}^1$,*
$$\sum_{I \in \mathcal{D}} \beta_\varphi(\lambda I)^2 < \infty. \tag{3.6}$$

Theorem 3.2.1 follows immediately from Lemmas 3.4.1 and 3.4.1, below.

Our proof of Theorem 3.2.1 uses a straightforward adaptation of the ideas presented in Section 6 of [Bis22]. We first reduce to the case where $\lambda = 1$.

**Lemma 3.2.2.** *Choose $\lambda \geqslant 1$, and let $\varphi : \mathbf{RP}^1 \to \mathbf{RP}^1$ be a homeomorphism. Then*

$$\sum_{I \in \mathcal{D}} \beta_\varphi(\lambda I)^2 < \infty \tag{3.7}$$

*for every dyadic decomposition $\mathcal{D}$ of $\mathbf{RP}^1$ if and only if*

$$\sum_{I \in \mathcal{D}} \beta_\varphi(I)^2 < \infty \tag{3.8}$$

*for every dyadic decomposition $\mathcal{D}$ of $\mathbf{RP}^1$.*

*Proof.* We first suppose that (3.7) holds for every dyadic decomposition $\mathcal{D}$ of $\mathbf{RP}^1$, and we show that (3.8) also holds for every such dyadic decomposition. Let $\mathcal{D}$ be one such dyadic decomposition. It suffices to show that for all $I$ in $\mathcal{D}$,

$$\beta_\varphi(I) \lesssim \beta_\varphi(\lambda I). \tag{3.9}$$

However, let $I \in \mathcal{D}$ be a dyadic interval. By finiteness in (3.7), we may suppose that $\beta_\varphi(\lambda I) < 1/16\lambda$. Let $a_I$ and $a_{\lambda I}$ denote respectively the $L^\infty$ best linear estimators of $\varphi$ over $I$ and $\lambda I$. By definition of $\beta_\varphi(I)$ and $\beta_\varphi(\lambda I)$,

$$\beta_\varphi(I) \cdot \ell(\varphi(I)) = \|\varphi - a_I\|_{L^\infty(I)} \leqslant \|\varphi - a_{\lambda I}\|_{L^\infty(I)} \leqslant \|\varphi - a_{\lambda I}\|_{L^\infty(\lambda I)} = \beta_\varphi(\lambda I) \cdot \ell(\varphi(\lambda I)),$$

so that, bearing in mind (3.1.10),

$$\beta_\varphi(I) \leqslant \frac{\ell(\varphi(\lambda I))}{\ell(\varphi(I))} \cdot \beta_\varphi(\lambda I) \leqslant 4\lambda \cdot \beta_\varphi(\lambda I),$$

and 3.9 follows, as desired.

We now suppose that (3.8) holds for every dyadic decomposition $\mathcal{D}$ of $\mathbf{RP}^1$, and we show that (3.7) also holds for every such dyadic decomposition. Let $\mathcal{D}$ be one such dyadic decomposition, and let $\mathcal{D}'$ and $\mathcal{D}''$ denote the two dyadic decompositions obtained by rotating $\mathcal{D}$ by $\pm\pi/3$. Let $k$ denote the least integer greater than $\ln(3(\lambda-1))/\ln(2)$. We define the function $\alpha : \mathcal{D}_{\geqslant k} \to \mathcal{D}' \cup \mathcal{D}''$ as follows. For all $m$, and for every dyadic interval $I \in \mathcal{D}_{m+k}$, there exists a dyadic interval $J \in \mathcal{D}'_m \cup \mathcal{D}''_m$ such that $\lambda I \subseteq J$. For all such $I$, we thus define $\alpha(I)$ to be one such $J$. As in the preceding paragraph, we now show that, for all $I$,

$$\beta_\varphi(\lambda I) \lesssim \beta_\varphi(\alpha(I)).$$

Note now that, for any dyadic interval $J \in \mathcal{D}' \cup \mathcal{D}''$, the preimage $\alpha^{-1}(\{J\})$ has cardinality at most $2^k$. Consequently,

$$\sum_{I \in \mathcal{D}_{\geqslant k}} \beta_\varphi(\lambda I)^2 \lesssim \sum_{I \in \mathcal{D}_{\geqslant k}} \beta_\varphi(\alpha(I))^2 \lesssim 2^k \sum_{J \in \mathcal{D}'} \beta_\varphi(J)^2 + 2^k \sum_{J \in \mathcal{D}''} \beta_\varphi(J)^2 < \infty.$$

Since $\mathcal{D}$ is arbitrary, (3.7) follows for every dyadic decomposition of $\mathbf{RP}^1$, and this completes the proof. □

**3.3. Weil–Petersson implies finite beta sum.** In this section, we will make considerable use of the results of Section 6 of [Bis22]. For the sake of compatibility with the notation of that paper, we identify $\mathbf{S}^1$ with the quotient $\mathbf{R}/2\pi\mathbf{Z}$, and we identify $\mathbf{RP}^1$ with $\mathbf{S}^1$ via the map

$$\mathbf{R}/\pi\mathbf{Z} \to \mathbf{R}/2\pi\mathbf{Z}; [x] \mapsto [2x].$$

Let $\varphi : \mathbf{S}^1 \to \mathbf{S}^1$ be a Weil–Petersson homeomorphism, let $\mathcal{D}$ be a dyadic decomposition of $\mathbf{S}^1$, let $\mathbf{D}$ denote the unit disk in $\mathbf{C}$, and let $u : \overline{\mathbf{D}} \to \mathbf{R}$ denote the harmonic extension of $\ln(\varphi')$. Recall from [Bis22] that

$$\|u\|_{H^1}^2 \simeq \|\ln(\varphi')\|_{H^{1/2}}^2.$$

In particular, since $\varphi$ is Weil–Petersson,

$$\int_{\mathbf{D}} |\nabla u|^2 \, dxdy \lesssim \|\ln(\varphi')\|_{H^{1/2}}^2 < \infty.$$

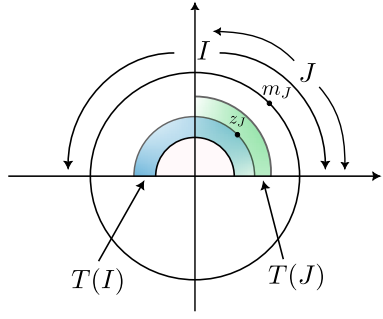

FIGURE 10. The "Carleson square" partition of the disk. More precisely, the partition used in this paper consists of the inner halves of all Carleson squares.

We now require the following definitions. First, for every dyadic interval $I \in \mathcal{D}$, we denote

$$T(I) := \left\{ re^{i\theta} \,\middle|\, e^{i\theta} \in I \text{ and } \frac{\ell(I)}{4\pi} \leqslant 1 - r \leqslant \frac{\ell(I)}{2\pi} \right\},$$

as in Figure 10. Note that the family $\{T(I) \mid I \in \mathcal{D}\}$ constitutes a partition of $\mathbf{D}$.

Next, for every dyadic interval $I \in \mathcal{D}$, and for all $\alpha \in\, ]0, 1[$, we define the quantity $\lambda_u(I, \alpha)$ such that

$$\lambda_u(I, \alpha)^2 := \sum_{J \subseteq 3I} \int_{T(J)} |\nabla u|^2 \frac{\ell(J)^\alpha}{\ell(I)^\alpha} \mathrm{dxdy}\,.$$

This is the quantity denoted by $\epsilon(u, Q_I, 1 - \alpha)$ in [Bis22]. In Lemma 6.1 of [Bis22], Bishop shows that, for all $\alpha \in [0, 1[$,

(3.10) $$\sum_{I \in \mathcal{D}} \lambda_u(I, \alpha)^2 \lesssim \int_{\mathbf{D}} |\nabla u|^2 \,\mathrm{dxdy} < \infty\,.$$

Finally, for every dyadic interval $I$, let $m_I \in \mathbf{S}^1$ denote its mid-point, and let $z_I$ denote the midpoint of the innermost side of $T(I)$, that is

$$z_I := \frac{\ell(I)}{2\pi} \cdot m_I\,,$$

as in Figure 10 again.

**Lemma 3.3.1.** *For every dyadic interval $I \in \mathcal{D}$, and for all $\alpha \in\, ]0, 1[$,*

(3.11) $$\int_I |\exp(u(s)) - \exp(u(z_I))|\, ds \lesssim \exp(u(z_I)) \cdot \ell(I) \cdot \lambda_u(I, \alpha)\,.$$

*Proof.* Denote $\tilde{u} := u - u(z_I)$. Trivially, $\nabla \tilde{u} = \nabla u$, so that, for all $\alpha$,

$$\lambda_{\tilde{u}}(I, \alpha) = \lambda_u(I, \alpha)\,.$$

Since $\tilde{u}(z(I)) = 0$, by Corollary 6.5 of [Bis22],[2]

$$\frac{1}{\exp(u(z_I)) \cdot \ell(I)} \int_I |\exp(u(s)) - \exp(u(z_I))|\, ds = \frac{1}{\ell(I)} \int_I |\exp(\tilde{u}(s)) - 1|\, ds$$

$$\lesssim \lambda_{\tilde{u}}(I, \alpha)$$

$$= \lambda_u(I, \alpha)\,,$$

and the result follows. □

---

[2]Note that, although the statement of Corollary 6.5 of [Bis22] claims only to estimate the integral of $(\exp(\tilde{u}(s)) - 1)$, as opposed to its absolute value, the estimate for the absolute value nonetheless follows on account of the Cauchy–Schwarz inequality used in the first step of the proof.

We are now ready to prove the first part of Theorem 3.2.1.

**Lemma 3.3.2.** *Let $\varphi : \mathbf{RP}^1 \to \mathbf{RP}^1$ be a homeomorphism. If $\varphi$ is Weil–Petersson then, for every dyadic decomposition $\mathcal{D}$ of $\mathbf{S}^1$,*

$$\sum_{I \in \mathcal{D}} \beta_\varphi(I)^2 < \infty . \tag{3.12}$$

*Proof.* We continue to use the above notation, and we fix $\alpha \in \,]0,1[$. Let $I \in \mathcal{D}$ be dyadic, denote $\gamma := \exp(u(z_I))$, and let $a : I \to \mathbf{R}$ denote the linear function defined such that

$$a' = \gamma \quad \text{and} \quad a(m_I) = \varphi(m_I) .$$

Since $\varphi' = \exp(u)$, for all $x \in I$,

$$|\varphi(x) - a(x)| \leqslant \int_I |\varphi'(y) - a'(y)| \, \mathrm{d}y .$$

Thus, by (3.11),

$$\|\varphi - a\|_{L^\infty(I)} \leqslant \int_I |\varphi'(y) - a'(y)| \, \mathrm{d}y = \int_I |\varphi'(y) - \gamma| \, \mathrm{d}y \lesssim \gamma \cdot \ell(I) \cdot \lambda_u(I, \alpha) .$$

Similarly,

$$\left|\ell(\varphi(I)) - \gamma \cdot \ell(I)\right| = \left|\int_I \varphi'(y) - \gamma \, \mathrm{d}y\right| \leqslant \int_I |\varphi'(y) - \gamma| \, \mathrm{d}y \lesssim \gamma \cdot \ell(I) \cdot \lambda_u(I, \alpha)$$

By (3.10), $\lambda_u(I, \alpha)$ tends to 0 as $\ell(I)$ tends to 0, so that, for $\ell(I)$ sufficiently small,

$$\gamma \cdot \ell(I) \lesssim \ell(\varphi(I)) .$$

Consequently

$$\beta_\varphi(I) \leqslant \frac{1}{\ell(\varphi(I))} \|\varphi - a\|_{L^\infty(I)} \lesssim \frac{\gamma \cdot \ell(I) \cdot \lambda_u(I, \alpha)}{\ell(\varphi(I))} \lesssim \lambda_u(I, \alpha) .$$

Bearing in mind (3.10) again, it follows that, for some subset $\mathcal{D}'$ of $\mathcal{D}$ containing all sufficiently small dyadic intervals,

$$\sum_{I \in \mathcal{D}'} \beta_\varphi(I)^2 \lesssim \sum_{I \in \mathcal{D}'} \lambda_u(I, \alpha)^2 \lesssim \int_{\mathbf{D}} |\nabla u|^2 \, \mathrm{dxdy} < \infty .$$

Since the dyadic decomposition $\mathcal{D}$ is arbitrary, the result follows. □

### 3.4. Finite beta sum implies Weil–Petersson.

The second part of the proof of Theorem 3.2.1 is a straightforward consequence of [Bis22].

**Lemma 3.4.1.** *Let $\varphi : \mathbf{RP}^1 \to \mathbf{RP}^1$ be a homeomorphism. If, for every dyadic decomposition $\mathcal{D}$ of $\mathbf{RP}^1$,*

$$\sum_{I \in \mathcal{D}} \beta_\varphi(I)^2 < \infty ,$$

*then $\varphi$ is Weil–Petersson.*

*Proof.* As in the preceeding section, we denote the mid-point of any subinterval $J$ of $\mathbf{RP}^1$ by $m_J$. Let $\mathcal{D}$ be a dyadic decomposition of $\mathbf{RP}^1$. For all dyadic $I \in \mathcal{D}$, we define

$$\mathrm{qs}_\varphi(I) := \frac{\left|\varphi(m_I) - m_{\varphi(I)}\right|}{\ell(\varphi(I))} ,$$

We claim that, for all dyadic $I \in \mathcal{D}$,

$$\mathrm{qs}_\varphi(I) \leqslant 2\beta_\varphi(I) .$$

Indeed, let $I$ be such a dyadic interval, and let $a : I \to \mathbf{RP}^1$ denote the $L^\infty$ best linear estimator of $\varphi$ over $I$. Denoting by $I_\pm$ the two extremities of $I$,

$$|\varphi(I_\pm) - a(I_\pm)| \leqslant \beta_\varphi(I) \cdot \ell(\varphi(I)) ,$$

so that
$$\left|m_{\varphi(I)} - a(m_I)\right| \leqslant \beta_\varphi(I) \cdot \ell(\varphi(I)) .$$
However,
$$|\varphi(m_I) - a(m_I)| \leqslant \beta_\varphi(I) \cdot \ell(\varphi(I)) ,$$
so that
$$\left|\varphi(m_I) - m_{\varphi(I)}\right| \leqslant 2\beta_\varphi(I) \cdot \ell(\varphi(I)) ,$$
from which the assertion follows. It follows that
$$\sum_{I \in \mathcal{D}} \mathrm{qs}_\varphi(I)^2 \leqslant 4 \sum_{I \in \mathcal{D}} \beta_\varphi(I)^2 < \infty .$$
Since $\mathcal{D}$ is arbitrary, the result now follows by [Bis22, Theorem 1.6]. □

## 4. Epsilon sums

Our second new characterization of the Weil–Petersson class will be in terms of the finiteness of what we call the epsilon sum. This is the analogue in the present context of Definition 14 of [Bis25]. We will show that finiteness of the beta sum implies finiteness of the epsilon sum. Conversely, we will show that finiteness of the epsilon sum implies finiteness of the beta sum under the additional hypothesis that the homeomorphism $\varphi : \mathbf{RP}^1 \to \mathbf{RP}^1$ is quasisymmetric.

### 4.1. Epsilon numbers.
In order to describe the epsilon sum characterization, we will continue to work with the quotient model $\mathbf{RP}^1 = \mathbf{R}/\pi\mathbf{Z}$, as in Section 2.2.1.

**Definition 4.1.1.** Let $\varphi : \mathbf{RP}^1 \to \mathbf{RP}^1$ be a homeomorphism with lift $\hat\varphi : \mathbf{R} \to \mathbf{R}$. Let $\mathcal{D}$ be a dyadic decomposition of $\mathbf{RP}^1$. Given a dyadic interval $I \in \mathcal{D}$, we denote by $\varepsilon_\varphi(I)$ the infimum of all $\varepsilon > 0$ for which there exists a point $x$ in $I$, with lift $\hat{x}$, and projective linear maps $f_\pm \in \mathrm{PSL}(2,\mathbf{R})$, with lifts $\hat{f}_\pm$, such that

$$\hat{f}_- \leqslant \hat{\varphi} \leqslant \hat{f}_+ , \tag{4.1}$$

and

$$\left|\hat{f}_\pm(\hat{x}) - \hat\varphi(\hat{x})\right| \leqslant \varepsilon \cdot \ell(\varphi(3I)) ,$$
$$\left|\hat{f}'_\pm(\hat{x}) - \frac{\ell(\varphi(3I))}{\ell(3I)}\right| \leqslant \varepsilon \cdot \frac{\ell(\varphi(3I))}{\ell(3I)} , \text{ and} \tag{4.2}$$
$$\left|\hat{f}''_\pm(\hat{x})\right| \leqslant \varepsilon \cdot \frac{\ell(\varphi(3I))}{\ell(3I)^2} .$$

By convention, we set $\varepsilon_\varphi(I) := 1$ whenever no such triplet $(x, f_\pm)$ exists.

*Remark* 4.1.2. Morally, Condition (4.1) means that, in $\mathbf{Ein}^{1,1}$, the graph of $f_-$ lies below that of $\varphi$, which in turn lies below that of $f_+$, and this is how it should be understood. Since $\mathbf{Ein}^{1,1}$ is not simply-connected, we work in the universal cover in order to provide a correct formal definition that we may work with. By mild abuse of terminology, we will say in the sequel that the graph of $f_-$ *lies below* that of $\varphi$ and the graph of $\varphi$ *lies below* that of $f_+$.

The quantity $\varepsilon_\varphi(I)$ is illustrated in Figure 11. It should be understood as a scale-invariant measure of the extent to which the graph of $\varphi$ differs from that of a fractional linear map. When $\varepsilon_\varphi(I)$ is small, the graph of $\varphi$ is pinched between those of the fractional linear maps $f_-$ and $f_+$, as in Figure 12. In the limit as $\varepsilon_\varphi(I)$ tends to zero, these two functions coincide, and $\varphi$ is also fractional linear.

The coefficients on the right hand side of (4.2) arise from the appropriate scales of the objects considered. Indeed, the scale in the codomain is measured by the length of $\varphi(3I)$, and the scale in the domain is measured by the length of $3I$. Each derivative introduces a factor on the scale of the domain, and this explains the above choices of coefficients.

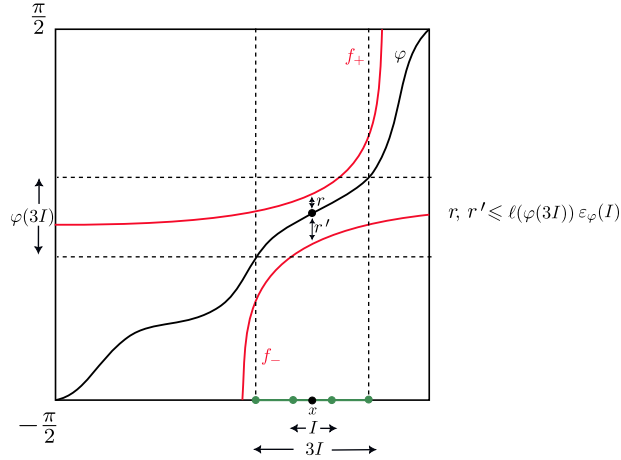

Figure 11. The definition of $\varepsilon_\varphi(I)$.

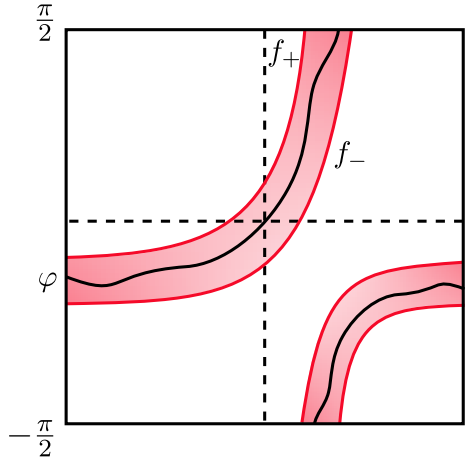

Figure 12. The definition of $\varepsilon_\varphi(I)$.

Although not strictly necessary for what follows, it is useful to note that, provided $\varepsilon_\varphi(I) < 1$, this quantity is always realized by some triplet $(x, f_\pm)$.

**Lemma 4.1.3.** *Let $\varphi : \mathbf{RP}^1 \to \mathbf{RP}^1$ be a homeomorphism, and let $\mathcal{D}$ be a dyadic decomposition of $\mathbf{RP}^1$. For every dyadic interval $I \in \mathcal{D}$ such that $\varepsilon_\varphi(I) < 1$, there exists a triplet $(x, f_\pm)$ realizing $\varepsilon_\varphi(I)$.*

*Proof.* For each $x \in I$, define $\varepsilon_\varphi(I, x)$ to be the infimum of all $\varepsilon > 0$ such that there exist fractional linear maps $f_\pm$, with lifts $\hat{f}_\pm$, satisfying (4.2) at $x$. For any $A, B > 0$, the set of fractional linear maps $f$ such that
$$|f(x) - \varphi(x)| \leq B, \qquad |f'(x) - A| \leq B, \qquad \text{and} \qquad |f''(x)| \leq B$$
is compact. Consequently, for all $x$, there exists a pair $(f_{\pm,x})$ realizing $\varepsilon_\varphi(I, x)$. Note now that $\varepsilon_\varphi(I, x)$ is lower semicontinuous in $x$. Thus, by compactness of $I$, $\varepsilon_\varphi(I, x)$ attains its minimum at some point of $I$, and this completes the proof. □

4.2. **The epsilon sum characterization of the Weil–Petersson class.** The main result of this section is the following characterization of the Weil–Petersson class in terms of epsilon numbers.

**Theorem 4.2.1.** *Let $\varphi : \mathbf{RP}^1 \to \mathbf{RP}^1$ be a homeomorphism. The following are equivalent.*
  (1) *$\varphi$ is Weil–Petersson.*

(2) $\varphi$ is quasisymmetric and there exists a dyadic decomposition $\mathcal{D}$ of $\mathbf{RP}^1$ such that

(4.3) $$\sum_{I \in \mathcal{D}} \varepsilon_\varphi(I)^2 < \infty \,.$$

(3) $\varphi$ is quasisymmetric and, for every dyadic decomposition $\mathcal{D}$ of $\mathbf{RP}^1$,

(4.4) $$\sum_{I \in \mathcal{D}} \varepsilon_\varphi(I)^2 < \infty \,.$$

Theorem 4.2.1 follows immediately from Corollaries 4.3.2 and 4.4.2, below.

*Remark* 4.2.2. The hypothesis of quasisymmetry can be weakened to $\mu$-Hölder, for some $\mu > 0$. This hypothesis is required to account for the distortions introduced by Penrose charts. Interestingly, it is not actually necessary. Indeed, in Section 5, below, we show by a different approach that finiteness of the epsilon sum alone is sufficient to prove the Weil–Petersson property. Lemma 4.3.1 and Corollary 4.3.2, below, are nonetheless of interest for being more direct at little extra cost.

Throughout this section, we will make considerable use of rotated Penrose charts, where projective linear maps take on a particularly simple form. Recall from Section 2.3.2 that $\mathbf{Ein}^{1,1}$ identifies with $\mathbf{RP}^1 \times \mathbf{RP}^1$. Furthermore, $\mathbf{RP}^1$ identifies with $\mathbf{R}/\pi\mathbf{Z}$, and the rotated Penrose chart about $(0,0)$ is given by

$$\mathcal{P}^{\mathrm{rot}} : (-\pi, \pi) \times (-\pi, \pi) \longrightarrow \mathbf{R} \times \mathbf{R}$$
$$(x_1, x_2) \longmapsto (t_1, t_2) := (\tan(x_1), \tan(x_2)) \,.$$

Note that the tangent function is linear up to third order. Consequently, the deviation of a Penrose chart from a linear map becomes arbitrarily small as we approach the origin.

### 4.3. The epsilon sum dominates the beta sum.
The first part of the proof of Theorem 4.3 is a consequence of the following estimate.

**Lemma 4.3.1.** *Let $\varphi : \mathbf{RP}^1 \to \mathbf{RP}^1$ be a homeomorphism and let $\mathcal{D}$ be a dyadic decomposition of $\mathbf{RP}^1$. For any dyadic interval $I \in \mathcal{D}$,*

(4.5) $$\beta_\varphi(I) \lesssim \varepsilon_\varphi(I) + \ell(I)^2 + \ell(\varphi(I))^2 \,.$$

**Corollary 4.3.2.** *Let $\varphi : \mathbf{RP}^1 \to \mathbf{RP}^1$ be a quasisymmetric homeomorphism. If*

$$\sum_{I \in \mathcal{D}} \varepsilon_\varphi(I)^2 < \infty$$

*for every dyadic decomposition $\mathcal{D}$ of $\mathbf{RP}^1$, then*

$$\sum_{I \in \mathcal{D}} \beta_\varphi(I)^2 < \infty$$

*for every dyadic decomposition $\mathcal{D}$ of $\mathbf{RP}^1$. In particular, when this holds, $\varphi$ is Weil–Petersson.*

*Proof.* Indeed, since $\varphi$ is quasisymmetric, it is $\mu$-Hölder for some $\mu > 0$. Let $\mathcal{D}$ be a dyadic decomposition of $\mathbf{RP}^1$. By (4.5) and the algebraic-geometric mean inequality,

$$\sum_{I \in \mathcal{D}} \beta_\varphi(I)^2 \lesssim \sum_{I \in \mathcal{D}} \varepsilon_\varphi(I)^2 + \sum_{I \in \mathcal{D}} \ell(I)^4 + \sum_{I \in \mathcal{D}} \ell(\varphi(I))^4 \,.$$

By hypothesis, the first term on the right hand side is finite. For all $m$, the subset $\mathcal{D}_m$ contains $2^m$ elements, each of length $\pi \cdot 2^{-m}$. Hence

$$\sum_{I \in \mathcal{D}} \ell(I)^4 = \sum_{m=0}^\infty \sum_{I \in \mathcal{D}_m} \pi^4 \cdot 2^{-4m} = \pi^4 \sum_{m=0}^\infty 2^{-3m} < \infty \,.$$

Finally, since $\varphi$ is a homeomorphism, for all $m$,

$$\sum_{I \in \mathcal{D}_m} \ell(\varphi(I)) = \pi \,.$$

Hence, bearing in mind that $\varphi$ is $\mu$-Hölder,

$$\sum_{I \in \mathcal{D}} \ell(\varphi(I))^4 = \sum_{m=0}^{\infty} \sum_{I \in \mathcal{D}_m} \ell(\varphi(I))^4$$
$$\leqslant \sum_{m=0}^{\infty} \left( \sup_{I \in \mathcal{D}_m} \ell(\varphi(I))^3 \right) \cdot \sum_{I \in \mathcal{D}_m} \ell(\varphi(I))$$
$$\lesssim \sum_{m=0}^{\infty} 2^{-3m\mu}$$
$$< \infty \,.$$

Upon combining these estimates, we see that

$$\sum_{I \in \mathcal{D}} \beta_\varphi(I)^2 < +\infty \,,$$

and since $\mathcal{D}$ is arbitrary, the result follows. $\square$

*Proof of Lemma 4.3.1.* Let $\mathcal{D}$ be a dyadic decomposition of $\mathbf{RP}^1$, and let $I \in \mathcal{D}$ be a dyadic interval. Since it suffices to prove the estimate for small $\ell(I)$, we may suppose that

$$\ell(I) < \delta \,,$$

for some $\delta > 0$ to be determined presently. Furthermore, since the result holds trivially for $\varepsilon_\varphi(I) \geqslant \delta$, we may likewise suppose that

$$\varepsilon_\varphi(I) < \delta \,.$$

Let $(x_0, f_\pm)$ realize $\varepsilon_\varphi(I)$, as in Lemma 4.1.3. Upon rotating $\mathbf{RP}^1 \times \mathbf{RP}^1$ if necessary, we may suppose that $\theta_0 = \varphi(\theta_0) = 0$. We now work in the rotated Penrose chart about $(0, 0)$. Let $\tilde{f}_\pm$ denote the image of $f_\pm$ in this chart, that is

(4.6) $$\tilde{f}_\pm = \tan \circ f_\pm \circ \arctan \,.$$

This is a fractional linear map satisfying

$$\tilde{f}_\pm(0) = \mathrm{O}\big(\varepsilon_\varphi(I) \cdot \ell(\varphi(I))\big) \,,$$
$$\tilde{f}'_\pm(0) = \frac{\ell(\varphi(I))}{\ell(I)} + \mathrm{O}\left( \frac{\varepsilon_\varphi(I) \cdot \ell(\varphi(I))}{\ell(I)} \right) \,, \text{ and}$$
$$\tilde{f}''_\pm(0) = \mathrm{O}\left( \frac{\varepsilon_\varphi(I) \cdot \ell(\varphi(I))}{\ell(I)^2} \right) \,.$$

Upon perturbing slightly if necessary, we may suppose that neither $\tilde{f}_\pm$ sends the point at infinity to itself, and we may then write

$$\tilde{f}_\pm(t) = \frac{P}{Q-t} - R = \frac{P}{Q} \cdot \frac{1}{1 - t/Q} - R \,,$$

where

$$\frac{P}{Q} - R = \mathrm{O}\big(\varepsilon_\varphi(I) \cdot \ell(\varphi(I))\big) \,,$$
$$\frac{P}{Q^2} = \frac{\ell(\varphi(I))}{\ell(I)} + \mathrm{O}\left( \frac{\varepsilon_\varphi(I) \cdot \ell(\varphi(I))}{\ell(I)} \right) \,, \text{ and}$$
$$\frac{P}{Q^3} = \mathrm{O}\left( \frac{\varepsilon_\varphi(I) \cdot \ell(\varphi(I))}{\ell(I)^2} \right) \,.$$

The quotient of the third identity by the second yields
$$\frac{\ell(I)}{Q} = \mathrm{O}(\varepsilon_\varphi(I))\,,$$
so that the quotient $t/Q$ is small for all $t \in [-2\ell(I), 2\ell(I)]$. Consequently, for all $t$ in this interval,

$$
\begin{aligned}
\tilde{f}_\pm(t) &= \frac{P}{Q} \cdot \frac{1}{(1 - t/Q)} - R \\
&= \frac{P}{Q} \cdot \left(1 + \frac{t}{Q} + \frac{t^2}{Q^2} \cdot \frac{1}{(1 - t/Q)}\right) - R \\
&= \frac{P}{Q^2} t + \left(\frac{P}{Q} - R\right) + \frac{P}{Q^3} \frac{t^2}{(1 - t/Q)} \\
&= \frac{\ell(\varphi(I))}{\ell(I)} t + \ell(\varphi(I)) \cdot \mathrm{O}\!\left(\varepsilon_\varphi(I)\right).
\end{aligned}
\tag{4.7}
$$

In particular, for all such $t$,

$$\tilde{f}_\pm(t) = \mathrm{O}\!\left(\ell(\varphi(I))\right). \tag{4.8}$$

It remains to transfer the estimate (4.7) back to the quotient model of $\mathbf{RP}^1$. To this end, choose $x \in I \subseteq [-\ell(I), \ell(I)]$ and denote $t := \tan(x)$. Note that, for $\delta$ sufficiently small, $t \in [-2\ell(I), 2\ell(I)]$, and the above estimates hold. Furthermore, by uniform continuity, $\ell(\varphi(I))$ tends to zero as $\ell(I)$ tends to zero, so that, by (4.8),

$$f_\pm(x) = \arctan\!\left(\tilde{f}_\pm(t)\right) = \tilde{f}_\pm(t) + \mathrm{O}\!\left(\tilde{f}_\pm(t)^3\right) = \tilde{f}_\pm(t) + \mathrm{O}\!\left(\ell(\varphi(I))^3\right).$$

Thus, by (4.7), for all $x \in I$,

$$
\begin{aligned}
f_\pm(x) &= \frac{\ell(\varphi(I))}{\ell(I)} \cdot t + \ell(\varphi(I)) \cdot \mathrm{O}\!\left(\varepsilon_\varphi(I) + \ell(\varphi(I))^2\right) \\
&= \frac{\ell(\varphi(I))}{\ell(I)} \cdot x + \ell(\varphi(I)) \cdot \mathrm{O}\!\left(\varepsilon_\varphi(I) + \ell(I)^2 + \ell(\varphi(I))^2\right).
\end{aligned}
$$

Thus, denoting
$$a(x) := \frac{\ell(\varphi(I))}{\ell(I)} \cdot x\,,$$
we obtain, over $I$,
$$|\varphi(x) - a(x)| \lesssim \ell(\varphi(I)) \cdot \left(\varepsilon_\varphi(I) + \ell(I)^2 + \ell(\varphi(I))^2\right),$$
so that
$$\beta_\varphi(I) \lesssim \varepsilon_\varphi(I) + \ell(I)^2 + \ell(\varphi(I))^2\,,$$
as desired. $\square$

### 4.4. The beta sum dominates the epsilon sum.
The second part of the proof of Theorem 4.2.1 is a consequence of the following estimate.

**Lemma 4.4.1.** *Let $\varphi : \mathbf{RP}^1 \to \mathbf{RP}^1$ be a homeomorphism, and let $\mathcal{D}$ be a dyadic decomposition of $\mathbf{RP}^1$. If*

$$\sum_{I \in \mathcal{D}} \beta_\varphi(3I)^2 < \infty\,, \tag{4.9}$$

*then, for all $\delta > 0$ and any dyadic interval $I \in \mathcal{D}$,*

$$\varepsilon_\varphi(I) \lesssim \ell(I)^{3/4} + \sum_{J \supseteq I} \beta_\varphi(J) \left(\frac{\ell(I)}{\ell(J)}\right)^{1-\delta}. \tag{4.10}$$

**Corollary 4.4.2.** *Let $\varphi : \mathbf{RP}^1 \to \mathbf{RP}^1$ be a homeomorphism, and let $\mathcal{D}$ be a dyadic decomposition of $\mathbf{RP}^1$. If*

$$\sum_{I \in \mathcal{D}} \beta_\varphi(3I)^2 < \infty,$$

*then,*

$$\sum_{I \in \mathcal{D}} \varepsilon_\varphi(I)^2 < \infty.$$

*Proof.* Indeed, by (4.10) and the algebraic-geometric mean inequality, for all $\delta > 0$,

$$\sum_{I \in \mathcal{D}} \varepsilon_\varphi(I)^2 \lesssim \sum_{I \in \mathcal{D}} \ell(I)^{3/2} + \sum_{I \in \mathcal{D}} \left( \sum_{J \supseteq I} \beta_\varphi(3J) \left( \frac{\ell(I)}{\ell(J)} \right)^{1-\delta} \right)^2.$$

For all $m$, the number of elements of $\mathcal{D}_m$ is equal to $2^m$. Hence,

$$\sum_{I \in \mathcal{D}} \ell(I)^{3/2} = \sum_{m=0}^{\infty} \sum_{I \in \mathcal{D}_m} \ell(I)^{3/2} \lesssim \sum_{m=0}^{\infty} 2^{-1/2} < \infty.$$

By the Cauchy–Schwarz inequality,

$$\sum_{I \in \mathcal{D}} \left( \sum_{J \supseteq I} \beta_\varphi(3J) \left( \frac{\ell(I)}{\ell(J)} \right)^{1-\delta} \right)^2 \leq \sum_{I \in \mathcal{D}} \left( \sum_{J \supseteq I} \beta_\varphi(3J)^2 \left( \frac{\ell(I)}{\ell(J)} \right)^{3/2-\delta} \right) \cdot \left( \sum_{J \supseteq I} \left( \frac{\ell(I)}{\ell(J)} \right)^{1/2-\delta} \right).$$

Since the second factor on the right-hand side is a geometric series, this yields

$$\sum_{I \in \mathcal{D}} \left( \sum_{J \supseteq I} \beta_\varphi(3J) \left( \frac{\ell(I)}{\ell(J)} \right)^{1-\delta} \right)^2 \lesssim \sum_{I \in \mathcal{D}} \sum_{J \supseteq I} \beta_\varphi(3J)^2 \left( \frac{\ell(I)}{\ell(J)} \right)^{3/2-\delta}$$

$$= \sum_{J \in \mathcal{D}} \beta_\varphi(3J)^2 \cdot \sum_{I \subseteq J} \left( \frac{\ell(I)}{\ell(J)} \right)^{3/2-\delta}$$

$$= \sum_{J \in \mathcal{D}} \beta_\varphi(3J)^2 \cdot \sum_{m=0}^{\infty} \sum_{\substack{I \subseteq J \\ \ell(I) = 2^{-m}\ell(J)}} \left( \frac{\ell(I)}{\ell(J)} \right)^{3/2-\delta}.$$

For all $J$, and for all $m$, $J$ contains $2^m$ dyadic intervals of length $2^{-m}\ell(J)$. Consequently, for $\delta < 1/2$,

$$\sum_{I \in \mathcal{D}} \left( \sum_{J \supseteq I} \beta_\varphi(3J) \left( \frac{\ell(I)}{\ell(J)} \right)^{(1-\delta)} \right)^2 \lesssim \sum_{J \in \mathcal{D}} \beta_\varphi(3J)^2 \cdot \sum_{m=0}^{\infty} 2^{-(1/2-\delta)} \lesssim \sum_{J \in \mathcal{D}} \beta_\varphi(3J)^2 < \infty.$$

The result follows upon combining these estimates. □

**4.5. Interlude.** The proof of Lemma 4.4.1 is the most technically challenging part of this paper. It is thus worthwhile before proceeding to reflect on the relatively straightforward intuition behind it. Consider first the case where $\varphi$ is twice differentiable, and suppose we wished to show that

$$\varphi(t) \leq \frac{1}{1-t},$$

over some interval $I$ containing the origin. A straightforward approach would be to show that $\varphi(0) = 1$, $\varphi'(0) = 1$, and that the second derivative of $\varphi$ is less than that of $1/(1-t)$ over this interval. The assertion would then follow upon integrating twice.

The argument used here is essentially a discretization of this technique, allowing for the fact that $\varphi$ is not necessarily differentiable. Indeed, we recall that, by Lemma 3.1.5, above, the quantity $\beta_\varphi(I)$ measures the variation of the mean gradient on a logarithmic scale, and may thus be viewed as a discrete analogue of the second derivative. Lemmas 4.5.1 and 4.5.2 then constitute respectively the discrete analogues of the first and the second integral. The remainder of the argument consists of comparing the resulting estimates with those arising from a suitably well-chosen fractional linear map, thereby yielding an upper bound for $\varepsilon_\varphi(I)$. Finally, the term $\ell(I)^{3/4}$ in (4.10) serves to exclude behaviours that may occur over the large scale with respect to the interval $I$.

4.5.1. *The first and second discrete integrals.* We now prove Lemma 4.4.1. Note first that, by (4.9),
$$\lim_{\ell(I) \to 0} \beta_\varphi(3I) = 0 \,.$$
As in Section 3.1, for all $I \in \mathcal{D}$, we denote by $a_{3I}$ the $L^\infty$ best linear estimator of $\varphi$ over $3I$, and we denote by $\gamma(3I)$ its gradient.

**Lemma 4.5.1.** *Choose $0 < \eta \leq 1$, and let*
$$I_0 \subseteq I_1 \subseteq \ldots \subseteq I_m$$
*be a sequence of successive dyadic intervals, such that, for all $k$,*
$$\beta_\varphi(3I_k) \leq \frac{2^\eta - 1}{32} \,.$$
*For all $k$,*

(4.11) $$\left|\gamma(3I_k) - \gamma(3I_0)\right| \leq 32 \cdot \gamma(3I_0) \cdot 2^{\eta(k-1)} \cdot \sum_{i=0}^{k} \beta_\varphi(3I_i) \,.$$

*Proof.* Choose $0 \leq k \leq m$. For brevity, for all $i$, we denote
$$\gamma_i := \gamma(3I_i) \qquad \text{and} \qquad \beta_i := \beta_\varphi(3I_i) \,.$$
For all $i$, since $\beta_i \leq 1/32$, (3.5) yields

(4.12) $$\frac{\gamma_{i+1}}{\gamma_i} \leq \frac{1}{1 - 16\beta_{i+1}} \leq 1 + 32\beta_{i+1} \,,$$

and

(4.13) $$\frac{\gamma_{i+1}}{\gamma_i} \geq \frac{1}{1 + 16\beta_{i+1}} \geq 1 - 32\beta_{i+1} \,.$$

Iterating (4.12) yields,
$$\gamma_k = \gamma_0 \cdot \prod_{i=0}^{k-1} \frac{\gamma_{i+1}}{\gamma_i}$$
$$\leq \gamma_0 \cdot \prod_{i=0}^{k-1} \left(1 + 32\beta_{i+1}\right)$$
$$= \gamma_0 + \gamma_0 \cdot \sum_{i=0}^{k-1} \left(32\beta_{i+1} \cdot \prod_{j=i+1}^{k-1} \left(1 + 32\beta_{j+1}\right)\right) \,.$$

Since, by hypothesis, $\beta_i \leq (2^\eta - 1)/32$ for all $i$, this yields,
$$\gamma_k - \gamma_0 \leq 32 \cdot \gamma_0 \cdot 2^{\eta(k-1)} \cdot \sum_{i=0}^{k} \beta_i \,,$$

as desired. The lower bound follows in a similar manner by iterating (4.13) instead of (4.12), and this completes the proof. □

In the sequel, given a dyadic interval $I$ and a point $x$ in $I$, we will denote by $I^x$ the translate of $I$ centered on $x$. In particular, for all such $I$ and $x$, the interval $I^x$ is contained in $3I$.

**Lemma 4.5.2.** *Choose $0 < \eta \leq 1$, and let*

$$I_0 \subseteq I_1 \subseteq \ldots \subseteq I_m$$

*be a sequence of successive dyadic intervals, such that, for all $k$,*

$$\beta_\varphi(3I_k) \leq \frac{2^\eta - 1}{32} \, .$$

*For all $x$ in $I_0$, for all $k$, and for all $y$ in $I_k^x$,*

$$(4.14) \qquad \left|\varphi(y) - \varphi(x) - \gamma(3I_0) \cdot (y-x)\right| \leq 76 \cdot \gamma(3I_0) \cdot \ell(I_k) \cdot 2^{\eta(k-1)} \cdot \sum_{i=0}^{k} \beta_\varphi(3I_i) \, .$$

*Proof.* Choose $0 \leq k \leq m$. For brevity, we denote, for all $i$,

$$a_i := a_{3I_i}, \qquad \gamma_i := \gamma(3I_i), \qquad \text{and} \qquad \beta_i := \beta_\varphi(3I_i) \, .$$

By definition of $a_k$, $\beta_k$, and $\gamma_k$, for all $y \in I_k^x$,

$$\begin{aligned}
\varphi(y) &\leq a_k(y) + \beta_k \cdot \ell(\varphi(3I_k)) \\
&= a_k(x) + \gamma_k \cdot (y-x) + \beta_k \cdot \ell(\varphi(3I_k)) \\
&\leq \varphi(x) + \gamma_k \cdot (y-x) + 2\beta_k \cdot \ell(\varphi(3I_k)) \, ,
\end{aligned}$$

so that

$$\varphi(y) - \varphi(x) - \gamma_0 \cdot (y-x) \leq (\gamma_k - \gamma_0) \cdot (y-x) + 2\beta_k \cdot \ell(\varphi(3I_k)) \, .$$

Note that $|y-x| < \ell(I_k)$. Applying (3.4) therefore yields

$$\varphi(y) - \varphi(x) - \gamma_0 \cdot (y-x) \leq (\gamma_k - \gamma_0) \cdot \ell(I_k) + \frac{2\beta_k \cdot \gamma_k \cdot \ell(3I_k)}{(1 - 2\beta_k)} \, ,$$

and, since $\beta_k < 1/4$,

$$\begin{aligned}
\varphi(y) - \varphi(x) - \gamma_0 \cdot (y-x) &\leq \left|\gamma_k - \gamma_0\right| \cdot \ell(I_k) + 12 \cdot \beta_k \cdot \gamma_k \cdot \ell(I_k) \\
&\leq \left|\gamma_k - \gamma_0\right| \cdot \ell(I_k) \cdot (1 + 12\beta_k) + 12 \cdot \beta_k \cdot \gamma_0 \cdot \ell(I_k) \, .
\end{aligned}$$

Since $\beta_i < (2^\eta - 1)/32$ for all $0 \leq i \leq k$, by (4.11),

$$\gamma_k \leq \gamma_0 + 32\gamma_0 \cdot 2^{\eta(k-1)} \cdot \sum_{i=0}^{k} \beta_i \, .$$

Hence

$$\varphi(y) - \varphi(x) - \gamma_0 \cdot (y-x) \leq 32 \cdot \gamma_0 \cdot \ell(I_k) \cdot (1 + 12\beta_k) \cdot 2^{\eta(k-1)} \cdot \sum_{i=0}^{k} \beta_i + 12\beta_k \cdot \gamma_0 \cdot \ell(I_k) \, .$$

Since, in particular, $\beta_k < 1/12$, this yields

$$\varphi(y) - \varphi(x) - \gamma_0 \cdot (y-x) \leq 76 \cdot \gamma_0 \cdot \ell(I_k) \cdot 2^{\eta(k-1)} \cdot \sum_{i=0}^{k} \beta_i \, ,$$

as desired. The lower bound is obtained in a similar manner, and this completes the proof. $\square$

4.6. **The local quadratic majorant I.** We now construct a local quadratic majorant of $\varphi$. We henceforth fix $0 < \eta \leq 1$, and let $\delta > 0$ be such that, for every dyadic interval $I$ of length less than $24\delta$,
$$\beta_\varphi(3I) < \frac{(2^\eta - 1)}{32}.$$

Let $I$ be a dyadic interval of length less than $\delta^4$, let $x$ be a point of $I$, and let $C$ be a positive constant to be determined presently. We define the real constants $P$, $Q$ and $R$ by

(4.15) $$\frac{1}{Q} := \ell(I)^{-\frac{1}{4}} + \frac{C}{\ell(I)} \cdot \sum_{\substack{I \subset J \\ \ell(J) < 24\delta}} \beta_\varphi(3J) \cdot \left(\frac{\ell(I)}{\ell(J)}\right)^{1-\eta},$$

(4.16) $$\frac{P}{Q^2} := \gamma(3I), \text{ and}$$

(4.17) $$\frac{P}{Q} - R := 2\beta_\varphi(3I) \cdot \ell(\varphi(3I)),$$

and we define the quadratic function $p$ by

(4.18) $$p(y) := \varphi(x) + \left(\frac{P}{Q} - R\right) + \frac{P}{Q^2} \cdot (y - x) + \frac{P}{16Q^3} \cdot (y - x)^2.$$

**Lemma 4.6.1.** *For all $y \in [x - 5Q, x + 3Q]$,*
$$\varphi(y) \leq p(y).$$

*Proof.* Indeed, let
$$I =: I_0 \subseteq I_1 \subseteq I_2 \subseteq \cdots \subseteq I_m$$
be a maximal sequence of successive dyadic intervals such that $\ell(I_m) \leq 24\delta$. In particular, by maximality
$$\ell(I_m) \geq 12\delta.$$
By hypothesis,
$$12\delta \geq 12\ell(I_0)^{1/4},$$
and, by (4.15),
$$12\ell(I_0)^{1/4} \geq 12Q,$$
so that
$$\ell(I_m) \geq 12Q.$$
In particular,
$$[x - 5Q, x + 3Q] \subseteq I_m.$$

Choose $1 \leq k \leq m$ and $y \in I_k^x \setminus I_{k-1}^x$. By definition of $P$, $Q$, $R$, and $p$,
$$p(y) - \varphi(x) - \gamma(3I_0) \cdot (y - x) = 2\beta_\varphi(3I_0) \cdot \ell(\varphi(3I_0)) + \frac{\gamma(3I_0)}{16Q} \cdot (y - x)^2.$$

Note that the first term on the right-hand side is non-negative. Since, furthermore,
$$|y - x| \geq \ell(I_k)/4,$$
this yields
$$p(y) - \varphi(x) - \gamma(3I_0) \cdot (y - x) \geq \frac{\gamma(3I_0) \cdot \ell(I_k)^2}{256Q}.$$

Consequently
$$p(y) - \varphi(y) = \big(p(y) - \varphi(x) - \gamma(3I_0) \cdot (y - x)\big) - \big(\varphi(y) - \varphi(x) - \gamma(3I_0) \cdot (x - y)\big)$$
$$\geq -\big(\varphi(y) - \varphi(x) - \gamma(3I_0) \cdot (x - y)\big) + \frac{\gamma(3I_0) \cdot \ell(I_k)^2}{256Q}.$$

Setting $C_1 := 76 \cdot 2^{-\eta}$, (4.14) now yields

$$p(y) - \varphi(y) \geq \frac{\gamma(3I_0) \cdot \ell(I_k)^2}{256Q} - C_1 \cdot \gamma(3I_0) \cdot \ell(I_k) \cdot 2^{\eta k} \cdot \sum_{i=0}^{k} \beta_\varphi(3I_i).$$

Note that, by Lemma 3.1.5, we may suppose that $\gamma(I_0) \geq 0$. Thus, since $y \in I_k^x \setminus I_{k-1}^x$ is arbitrary, $p \geq \varphi$ over $I_k^x \setminus I_{k-1}^x$ provided that

$$\frac{1}{Q} \geq \frac{256 \cdot C_1}{\ell(I_k)} \cdot 2^{\eta k} \cdot \sum_{i=0}^{k} \beta_\varphi(3I_i).$$

Since $\ell(I_k) = 2^k \ell(I_0)$, this in turn holds provided that

(4.19) $$\frac{1}{Q} \geq \frac{256 \cdot C_1}{\ell(I_0)} \cdot 2^{(\eta-1)k} \cdot \sum_{i=0}^{k} \beta_\varphi(3I_i).$$

This inequality holds for all $1 \leq k \leq m$ provided that

$$\frac{1}{Q} \geq \frac{256 \cdot C_1}{\ell(I_0)} \cdot \sum_{k=0}^{m} \left( 2^{(\eta-1)k} \cdot \sum_{i=0}^{k} \beta_\varphi(3I_i) \right).$$

However,

$$\sum_{k=0}^{m} \left( 2^{(\eta-1)k} \cdot \sum_{i=1}^{k} \beta_\varphi(3I_i) \right) = \sum_{i=0}^{m} \left( \beta_\varphi(3I_i) \cdot \sum_{k=i}^{m} 2^{(\eta-1)k} \right)$$
$$\leq \frac{1}{(1 - 2^{(\eta-1)})} \cdot \sum_{i=0}^{m} 2^{(\eta-1)i} \cdot \beta_\varphi(3I_i)$$
$$= \frac{1}{(1 - 2^{(\eta-1)})} \cdot \sum_{i=0}^{m} \beta_\varphi(3I_i) \cdot \left( \frac{\ell(I_0)}{\ell(I_i)} \right)^{1-\eta}.$$

Thus, by definition of $Q$, upon setting

$$C := \frac{256 \cdot C_1}{(1 - 2^{(\eta-1)})},$$

(4.19) holds for all $1 \leq k \leq m$, so that $p \geq \varphi$ over $I_m^x \setminus I_0^x$. Finally, by definition of $P$, $Q$, $R$, and $p$, for all $y \in I_0^x$,

$$p(y) \geq \varphi(x) + \gamma(3I_0) \cdot (y - x) + 2\beta_\varphi(3I_0) \cdot \ell(\varphi(3I_0)).$$

Since $\gamma(3I_0)$ is the gradient of the $L^\infty$ best estimator of $\varphi$ over $3I_0$, this yields, for all such $y$,

$$p(y) \geq \varphi(y).$$

It follows that, with the above chosen value of $C$, $p \geq \varphi$ over the whole of $I_m$, and hence over the whole of $[x - 5Q, x + 3Q]$, as desired. □

**4.7. The local quadratic majorant II.** We now work in a rotated Penrose chart, as in Section 2.3.3. Upon rotating $\mathbf{RP}^1 \times \mathbf{RP}^1$, we may suppose first that

$$x = 0 \quad \text{and} \quad p(0) = \varphi(0) + \left( \frac{P}{Q} - R \right) = 0.$$

We denote by $\tilde{\varphi}$ and $\tilde{p}$ the respective images of $\tilde{\varphi}$ and $\tilde{p}$ in the rotated Penrose chart about $(0, 0)$, that is

$$\tilde{\varphi} := \tan \circ \varphi \circ \arctan \quad \text{and} \quad \tilde{p} := \tan \circ p \circ \arctan.$$

It is precisely at this stage that we require the Hölder property of $\varphi$.

**Lemma 4.7.1.** *The map $\varphi$ is $\mu$-Hölder for all $\mu < 1$.*

*Proof.* Indeed, by [GR15], $\varphi \in H^s(\mathbf{RP}^1)$ for all $s < 3/2$, and the result now follows by the Sobolev embedding theorem. □

**Lemma 4.7.2.** *Over the interval* $[-4Q, 2Q]$,

$$\tilde{\varphi}(t) \leqslant \tilde{p}(t) \leqslant \frac{P}{Q^2} \cdot t + \frac{P}{8Q^3} \cdot t^2 . \tag{4.20}$$

*Proof.* Consider first the quadratic polynomial

$$q(x) := \frac{A}{B} \cdot x + \frac{A}{16B^2} \cdot x^2 ,$$

for suitable constants $A, B > 0$. Let $\tilde{q}$ denote its image in the rotated Penrose chart about $(0, 0)$, that is

$$\tilde{q} := \tan \circ q \circ \arctan .$$

We claim that there exists $\delta > 0$ such that, for $0 < A, B < \delta$, and for all $t \in [-4B, 2B]$,

$$\tilde{q}(t) \leqslant \frac{A}{B} \cdot t + \frac{A}{8B^2} \cdot t^2 . \tag{4.21}$$

Indeed, the Taylor series of tan and arctan at 0 are

$$\tan(x) = x\bigl(1 + O(x^2)\bigr) \quad \text{and} \quad \arctan(t) = t\bigl(1 + O(t^2)\bigr) .$$

Thus, for $t \in [-4B, 2B]$,

$$q(\arctan(t)) = \frac{A}{B} \cdot t\bigl(1 + O(t^2)\bigr) + \frac{A}{16B^2} \cdot t^2\bigl(1 + O(t^2)\bigr)$$

$$= \frac{A}{B} \cdot t + \frac{A}{16B^2} \cdot t^2 \cdot \bigl(1 + O(B^2)\bigr) ,$$

It follows that, for $\delta$ sufficiently small, and for all $t$ in this interval,

$$q(\arctan(t)) \leqslant \frac{A}{B} \cdot t + \frac{3A}{32B^2} \cdot t^2 .$$

Composing with the tangent now yields, for all $t \in [-4B, 2B]$,

$$\tilde{q}(t) = \tan(q(\arctan(t)))$$

$$= \frac{A}{B} \cdot t + \frac{3A}{32B^2} \cdot t^2 + O\left(\frac{A^3}{B^3} \cdot t^3 + \frac{A^3}{B^6} \cdot t^6\right)$$

$$= \frac{A}{B} \cdot t + \frac{3A}{32B^2} \cdot t^2 \cdot \bigl(1 + O(A^2)\bigr) .$$

Thus, for $\delta$ sufficiently small, and for all $t$ in this interval,

$$\tilde{q}(t) \leqslant \frac{A}{B} \cdot t + \frac{A}{8B^2} \cdot t^2 ,$$

as asserted.

By (4.2),

$$Q \leqslant \ell(I)^{1/4} .$$

Thus, bearing in mind (3.4) and (4.16),

$$\frac{P}{Q} = \frac{P}{Q^2} \cdot Q = \gamma(3I) \cdot Q \lesssim \frac{\ell(\varphi(3I))}{\ell(3I)} \cdot Q \lesssim \ell(\varphi(3I)) \cdot \ell(I)^{-3/4} .$$

However, by Lemma 4.7.1, $\varphi$ is $\mu$–Hölder for all $\mu < 1$, so that

$$\frac{P}{Q} \lesssim \ell(I)^{1/8} .$$

It follows that both $P/Q$ and $Q$ tend to zero as $\ell(I)$ tends to zero. The result now follows upon applying the argument of the preceding paragraph with $A := P/Q$ and $B := Q$. □

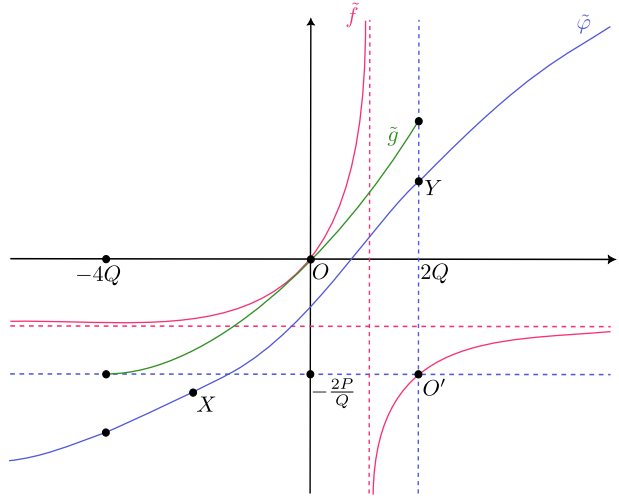

Figure 13. The quadratic and fractional linear majorants viewed in the rotated Penrose chart.

### 4.8. The fractional linear majorant.
It now only remains to construct a fractional linear majorant of $\varphi$ for $\ell(I)$ sufficiently small.

**Lemma 4.8.1.** *For $\ell(I)$ sufficiently small, the graph of $\varphi$ lies below the graph of the unique fractional linear map $f$ satisfying*

$$(4.22) \qquad f(x) = \varphi(x) + \left(\frac{P}{Q} - R\right), \quad f'(x) = \frac{P}{Q^2} \quad \text{and} \quad f''(x) = \frac{2P}{Q^3}.$$

*Proof.* Note first that, up to rotating $\mathbf{RP}^1 \times \mathbf{RP}^1$, we may suppose that

$$(4.23) \qquad x = 0 \quad \text{and} \quad p(x) = \varphi(x) + \left(\frac{P}{Q} - R\right) = 0.$$

Now let $f$ be as in (4.22), and denote by $\tilde{f}$ its image in the rotated Penrose chart, that is

$$\tilde{f} = \tan \circ f \circ \arctan.$$

Since $\tan(0) = \tan''(0) = 0$ and $\tan'(0) = 1$,

$$\tilde{f}(0) = f(0) = 0, \quad \tilde{f}'(0) = f'(0) = \frac{P}{Q^2}, \quad \text{and} \quad \tilde{f}''(0) = f''(0) = \frac{2P}{Q^3}.$$

Since $\tilde{f}$ is a fractional linear map, this yields

$$\tilde{f}(t) = \frac{P}{Q-t} - \frac{P}{Q}.$$

In particular, the graph of $\tilde{f}$ passes through the point $(0,0)$, has a horizontal asymptote at $\{y = -P/Q\}$, a vertical asymptote at $\{x = Q\}$, and is centered on the point $(Q, -P/Q)$, as in Figure 13.

Consider now the quadratic function

$$\tilde{g}(t) := \frac{P}{Q^2} \cdot t + \frac{P}{8Q^3} \cdot t^2.$$

Note that

$$\tilde{f}(t) - \tilde{g}(t) = \frac{P t^2 \cdot (t + 7Q)}{8Q^3(Q-t)},$$

so that, over the interval $[-7Q, Q[$,

$$\tilde{g} \leqslant \tilde{f}.$$

However, by Lemma 4.7.2, over the interval $[-4Q, 2Q]$,

$$\tilde{\varphi} \leqslant \tilde{g},$$

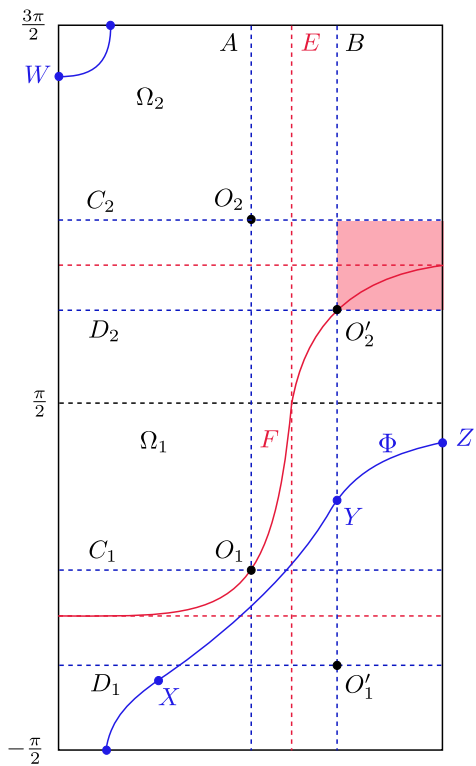

FIGURE 14. The fractional linear majorant viewed in the double cover of the matrix angle model.

so that, over the interval $[-4Q, Q[$

$$\tilde{\varphi} \leq \tilde{f} \ .$$

It remains only to show that the graph of $\varphi$ does not meet the graph of $f$ at any other point. Indeed, let $\hat{\varphi} : \mathbf{R} \to \mathbf{R}$ be a lift of $\varphi$, let $\hat{x} \in \mathbf{R}$ be a lift of $x$, and let $\hat{f}_+$ denote the unique lift of $f$ satisfying

$$\hat{f}_+(\hat{x}) = \hat{\varphi}(\hat{x}) + \left(\frac{P}{Q} - R\right).$$

If the graph of $\varphi$ does not meet the graph of $f$ at any other point, then it follows, in particular, that it does not *cross* the graph of $f$ at any point. Since the same must also hold for the lifts, we have that

$$\hat{\varphi} \leq \hat{f}_+ \ ,$$

over the whole of $\mathbf{R}$, and the result follows.

In order to show that the graph of $\varphi$ does not meet that of $f$ at any other point, it is now helpful to work in the double cover of the matrix angle model, described in Section 2.3.2, and whose fundamental domain is shown in Figure 14. Let $\Omega_1$ and $\Omega_2$ denote the two fundamental domains of this double cover. Denote $O := (0, 0)$, and let $O'$ denote the image of this point under reflection through the centre of the graph of $f$. Let $O_1$ and $O_2$ (resp. $O'_1$ and $O'_2$) denote the respective lifts of $O$ (resp. $O'$) in $\Omega_1$ and $\Omega_2$. Let $A$ (resp. $B$) denote the vertical line passing through $O_1$ and $O_2$ (resp. $O'_1$ and $O'_2$), and let $C_1$ and $C_2$ (resp. $D_1$ and $D_2$) denote the respective horizontal lines passing through $O_1$ and $O_2$ (resp. $O'_1$ and $O'_2$). Let $F$ denote the unique lift of the graph of $f$ passing through $O_1$, and let $\Phi$ denote the unique lift of the graph of $\varphi$ passing near this point. This construction is illustrated in Figure 14.

It will suffice to prove that $\Phi$ does not meet $F$. To show this, we will compare the geometries of the graphs of $\tilde{f}$ and $\tilde{\varphi}$ in the rotated Penrose chart, shown in Figure 13, with the geometries of $\Phi$ and $F$ in the double cover of the matrix angle model, shown in Figure 14.

We first claim that $\Phi$ contains a point $X$, say, lying below $D_1$ and to the left of $A$ in $\Omega_1$. Indeed, in the rotated Penrose chart,

$$\tilde{g}(-4Q) = -\frac{2P}{Q}.$$

Since $\tilde{\varphi}(0)$ is finite, and since $\tilde{\varphi} \leqslant \tilde{g}$ over $[-4Q, 0]$, it follows by continuity that there exists $t \in [-4Q, 0]$ such that

$$\tilde{\varphi}(t) \leqslant -\frac{2P}{Q}.$$

Let $X$ denote the point of $\Phi$ that projects onto $(t, \tilde{\varphi}(t))$. This point lies below $D_1$ and to the left of $A$, as desired.

We now claim that $\Phi$ passes through $B$ at some point $Y$, say, of $\Omega_1$. Indeed, in the rotated Penrose chart, since $\tilde{\varphi}(0)$ is finite, since $\tilde{\varphi} \leqslant \tilde{g}$ over $[0, 2Q]$, and since $\tilde{\varphi}$ is increasing over $[0, 2Q]$, this function is finite over the whole of this interval. Let $Y$ denote the point of $\Phi$ that projects onto $(2Q, \tilde{\varphi}(2Q))$. This point lies in $B$ and in $\Omega_1$, as desired.

Note now that the projection of $\Phi$ onto the vertical axis has length $\pi$. Thus, since it contains the point $X$ lying below $D_1$, it cannot meet $D_2$. It follows by monotonicity that the portion of $\Phi$ lying to the right of $X$ can only leave the fundamental domain of the double cover by passing through the right-hand side of this domain at some point $Z$, say, as in Figure 14. Furthermore, the portion of $\Phi$ between $X$ and $Z$ lies below $D_2$. In particular, the portion of $\Phi$ between $Y$ and $Z$ does not meet the region of $\Omega_2$ lying between $C_2$ and $D_2$ and to the right of $B$ and therefore also does not meet $F$.

By a similar reasoning, $\Phi$ meets the left-hand side of the double cover at some point $W$, say, lying above $C_1$, as in Figure 14, and the portion of $\Phi$ lying between $W$ and $X$ does not meet $F$ either. It follows that the graph of $\varphi$ does not meet the graph of $f$ at any other point, and this completes the proof. □

'

We are now ready to prove Lemma 4.4.1.

*Proof of Lemma 4.4.1.* Let $f_+ := f$ be as in Lemma 4.8.1 so that the graph of $\varphi$ lies below that of $f_+$. By definition of $P$, $Q$, and $R$,

$$|f_+(x) - \varphi(x)| \leqslant 2\beta_\varphi(3I) \cdot \ell(\varphi(3I)).$$

Bearing in mind in addition Lemma (3.1.5),

$$\left| f'_+(x) - \frac{\ell(\varphi(3I))}{\ell(3I)} \right| \leqslant \left| f'_+(x) - \gamma(3I) \right| + \left| \gamma(3I) - \frac{\ell(\varphi(3I))}{\ell(3I)} \right| \leqslant 2\beta_\varphi(3I) \cdot \frac{\ell(\varphi(3I))}{\ell(3I)}.$$

Likewise, for $\beta_\varphi(3I) < 1/2$,

$$\left| f''_+(x) \right| \leqslant \gamma(3I) \cdot \ell(I)^{-1/4} + \frac{C \cdot \gamma(3I)}{\ell(I)} \cdot \sum_{\substack{I \subseteq J \\ \ell(J) < 24\delta}} \beta_\varphi(3J) \cdot \left( \frac{\ell(I)}{\ell(J)} \right)^{1-\eta}$$

$$\leqslant \left( 6\ell(I)^{3/4} + 6C \sum_{\substack{I \subseteq J \\ \ell(J) < 24\delta}} \beta_\varphi(3J) \cdot \left( \frac{\ell(I)}{\ell(J)} \right)^{1-\eta} \right) \cdot \frac{\ell(\varphi(3I))}{\ell(3I)^2}.$$

By symmetry, $\varphi$ also lies above the graph of a projective linear map $f_-$ with derivatives up to order 2 satisfying similar estimates. It follows by (4.2) that

$$\varepsilon_\varphi(I) \lesssim \ell(I)^{3/4} + \sum_{\substack{I \subseteq J \\ \ell(J) < 24\delta}} \beta_\varphi(3J) \cdot \left( \frac{\ell(I)}{\ell(J)} \right)^{1-\eta},$$

as desired. □

## 5. Complete maximal spacelike surfaces and renormalized area

Our third, and main, new characterization of Weil–Peterson homeomorphisms will be in terms of the renormalized area of maximal surfaces in $\mathbf{AdS}^{2,1}$. Using an estimate proven by Seppi in [Sep19] we show that if a homeomorphism $\varphi : \mathbf{RP}^1 \to \mathbf{RP}^1$ has finite epsilon sum, then it bounds a complete maximal spacelike surface of finite renormalized area.

### 5.1. Complete maximal spacelike surfaces.
Recall from Section 1.2 that the graph of a homeomorphism $\varphi : \mathbf{RP}^1 \to \mathbf{RP}^1$ defines a simple closed acausal curve $\Lambda_\varphi = \mathrm{graph}(\varphi)$ in $\mathbf{Ein}^{1,1}$ and that this curve bounds a unique entire maximal surface $\Sigma_\varphi$ in $\mathbf{AdS}^{2,1}$.

**Theorem 5.1.1.** *Let $\varphi : \mathbf{RP}^1 \to \mathbf{RP}^1$ be a homeomorphism. If*

$$\sum_{I \in \mathcal{D}} \varepsilon_\varphi(I)^2 < \infty \tag{5.1}$$

*for some dyadic decomposition $\mathcal{D}$ of $\mathbf{RP}^1$, then $\Sigma_\varphi$ has finite renormalized area, that is*

$$\int_{\Sigma_\varphi} \|\mathrm{I\!I}_\varphi\|^2 \mathrm{d\,Area} < +\infty \,, \tag{5.2}$$

*where $\mathrm{I\!I}_\varphi$ denotes the second fundamental form of $\Sigma_\varphi$.*

Theorem 5.1.1 is proven in Section 5.4 below. Let $\varphi : \mathbf{RP}^1 \to \mathbf{RP}^1$ be a homeomorphism and let $\mathcal{D}$ be a dyadic decomposition of $\mathbf{RP}^1$ for which (5.1) holds. Trivially, for all but finitely many dyadic intervals $I \in \mathcal{D}$,

$$\varepsilon_\varphi(I) < 1 \,. \tag{5.3}$$

In what follows, we will only be concerned with those dyadic intervals $I$ for which this property holds. For all such $I$, let $(x_I, f_{\pm,I})$ be a triplet realizing $\varepsilon_\varphi(I)$, as in Lemma 4.1.3.

### 5.2. Diamonds.
We first use diamonds to associate to every dyadic interval $I$ a compact subset $K_I$ of $\Sigma_\varphi$. We recall that, given two time-related points $p_1$ and $p_2$ in a globally hyperbolic lorentzian manifold $X$, their *diamond* $\mathrm{D}(p_1, p_2)$ is the set of all points lying in the past of $p_1$ and in the future of $p_2$, that is

$$\mathrm{D}(p_1, p_2) := \mathrm{C}^-(p_1) \cap \mathrm{C}^+(p_2) \,,$$

where $\mathrm{C}^-(p_1)$ and $\mathrm{C}^+(p_2)$ denote respectively the chronological past and future of $p_1$ and $p_2$.

In the present case, neither $\mathbf{AdS}^{2,1}$ nor $\mathbf{Ein}^{1,1}$ is globally hyperbolic, since both admit closed causal curves. The past and future of any point thus coincides with the entire space, and the above definition is of no practical use. This is however easily addressed by the following mild abuse of conventions. In what follows, we will work exclusively in the domain of some kleinian chart of $\mathbf{AdS}^{2,1}$, and in the domain of some (rotated) Penrose chart of $\mathbf{Ein}^{1,1}$. In both cases, the domain is globally hyperbolic, and diamonds are thus well-defined and non-trivial.

We first describe diamonds in $\mathbf{Ein}^{1,1}$. Consider two time-related points $p_1$ and $p_2$ in this space. We suppose that these points lie in the domain of some rotated Penrose chart which we identify with $]-\pi/2, \pi/2[\times]-\pi/2, \pi/2[$, so that

$$p_1 = (x_{1,1}, x_{1,2}) \quad \text{and} \quad p_2 = (x_{2,1}, x_{2,2}) \,.$$

We suppose furthermore that

$$x_{1,1} < x_{2,1} \quad \text{and} \quad x_{1,2} > x_{2,2} \,.$$

Recalling the causal structure of the matrix angle model 2.3.2, the diamond of $p_1$ and $p_2$ is then

$$\mathrm{D}^{\mathrm{ein}}(p_1, p_2) = [x_{1,1}, x_{2,1}] \times [x_{2,2}, x_{1,2}] \,, \tag{5.4}$$

as in Figure 15.

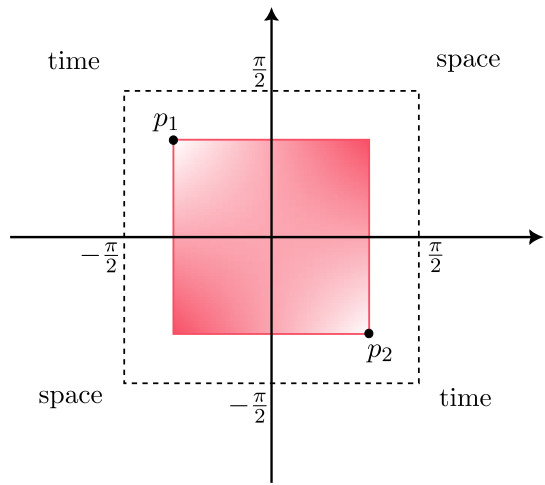

Figure 15. The diamond of $p_1$ and $p_2$.

We now describe diamonds in $\mathbf{AdS}^{2,1}$. We view the same points $p_1$ and $p_2$ as boundary points of $\mathbf{AdS}^{2,1}$, and, recalling the above mild abuse of convention, we denote by $\mathrm{D}^{\mathrm{ads}}(p_1, p_2)$ their diamond in $\mathbf{AdS}^{2,1}$. In particular,

$$\mathrm{D}^{\mathrm{ein}}(p_1, p_2) = \partial_\infty \mathrm{D}^{\mathrm{ads}}(p_1, p_2) . \tag{5.5}$$

We now work in the kleinian chart of $\mathbf{AdS}^{2,1}$ described in Section 2.1.4. Let $\tilde{\mathrm{D}}^{\mathrm{ads}}(p_1, p_2)$, $\tilde{p}_1$ and $\tilde{p}_2$ denote the respective images of $\mathrm{D}^{\mathrm{ads}}(p_1, p_2)$, $p_1$ and $p_2$ in this chart. Recall that, in this chart, lightlike geodesics identify with straight lines which are tangent to $\partial_\infty \mathbf{AdS}^{2,1} = \mathbf{Ein}^{1,1}$. It is then straightforward to show that $\tilde{\mathrm{D}}^{\mathrm{ein}}(p_1, p_2)$ is the set of all points of $\mathbf{AdS}^{2,1}$ lying below $\mathrm{T}_{\tilde{p}_1} \mathbf{Ein}^{1,1}$ and above $\mathrm{T}_{\tilde{p}_2} \mathbf{Ein}^{1,1}$ .

We also require the limiting case of diamonds of ideal boundary points in the kleinian chart. These are visualized as follows. We identify ideal points in this chart with causal rays in $\mathbf{R}^{2,2}$ leaving the origin, and we identify ideal boundary points with lightlike rays in $\mathbf{R}^{2,2}$ leaving the origin (note that two ideal points *in this chart* may identify with the same point in $\mathbf{AdS}^{2,1}$). When $p_1$ and $p_2$ correspond to the ideal points given by the lightlike rays $\tilde{R}_1$ and $\tilde{R}_2$, $\tilde{\mathrm{D}}^{\mathrm{ads}}(p_1, p_2)$ is the set of all points of $\mathbf{AdS}^{2,1}$ lying below the plane $\tilde{R}_1^\perp$ and above the plane $\tilde{R}_2^\perp$, as in Figure 16. Note that this is the limit of the case discussed in the preceding paragraph, as can be seen upon observing that, for all $p \in \mathbf{Ein}^{1,1}$, $\mathrm{T}_p \mathbf{Ein}^{1,1}$ is parallel to $\langle p \rangle^\perp$.

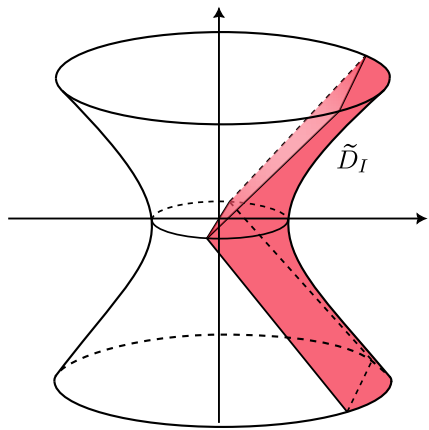

Figure 16. The diamond of two ideal points in the kleinian chart.

We now associate to every sufficiently small dyadic interval $I$ a compact subset $K_I$ of $\Sigma_\varphi$ as follows. Given a dyadic interval $I$ of depth $m \geqslant 2$, we define

$$D_I := D^{\text{ads}}\big(\big((3I)_-, \varphi(3I)_+\big), \big((3I)_+, \varphi(3I)_-\big)\big),$$

where $(3I)_\pm$ and $\varphi(3I)_\pm$ denote the respective end-points of $3I$ and $\varphi(3I)$. For every such dyadic interval $I$, we then define

$$(5.6) \qquad K_I := \Sigma_\varphi \cap \bigg(D_I \setminus \bigcup_{\substack{J \subseteq 3I \\ \ell(J) = \ell(I)/2}} \text{Int}(D_J)\bigg).$$

**Lemma 5.2.1.** *The family $(K_I)_{I \in \mathcal{D}_{\geqslant 2}}$ covers the complement of a compact subset of $\Sigma_\varphi$.*

*Proof.* Indeed, denoting

$$\Omega := \bigcup_{I \in \mathcal{D}_{\geqslant 2}} D_I,$$

we obtain

$$\Sigma_\varphi \cap \Omega \subseteq \bigcup_{I \in \mathcal{D}_{\geqslant 2}} K_I.$$

However, by (5.4) and (5.5), for all $I$,

$$\partial_\infty D_I = D^{\text{ein}}\big(\big((3I)_-, \varphi(3I)_+\big), \big((3I)_+, \varphi(3I)_-\big)\big) = 3I \times \varphi(3I).$$

In particular, the graph of the restriction of $\varphi$ to $3I$ is contained in this diamond. Since $I \in \mathcal{D}_{\geqslant 2}$ is arbitrary, $\Omega$ contains the graph of $\varphi$. Since the graph of $\varphi$ is the ideal boundary of $\Sigma_\varphi$, the complement of $\Omega$ in $\Sigma_\varphi$ is compact, and the result follows. $\square$

### 5.3. The canonical normalization.

We now describe, for every sufficiently small dyadic interval $I$, a canonical transformation $T_I \in \textbf{PSL}(2,\textbf{R}) \times \textbf{PSL}(2,\textbf{R})$ which sends the rectangle $3I \times \varphi(3I)$ to the square of side length $\pi/2$ centred on the origin. We continue to work in the matrix angle model of $\textbf{Ein}^{1,1}$.

Let $I \in \mathcal{D}$ be a dyadic interval. Let $E_1, E_2 \in \textbf{PSL}(2, \textbf{R})$ denote the unique rotations such that

$$E_1 \cdot x_I = E_2 \cdot \varphi(x_I) = 0.$$

**Lemma 5.3.1.** *There exist unique hyperbolic elements $H_1$ and $H_2$, each fixing $0$ and $\pi/2$, and unique parabolic elements $P_1$ and $P_2$, each fixing $\pi/2$, such that*

(1) *each of $(P_1 H_1 E_1) \cdot (3I)$ and $(P_2 H_2 E_2) \cdot \varphi(3I)$ has length $\pi/2$ ; and*

(2) *each of $(P_1 H_1 E_1) \cdot (3I)$ and $(P_2 H_2 E_2) \cdot \varphi(3I)$ is centered on $0$.*

*Proof.* Let $\widetilde{3I}$ denote the image of $E_1 \cdot (3I)$ in the affine chart of $\textbf{RP}^1$ about $0$, let $\lambda_1$ denote its length, and let $m_1$ denote its mid-point. Let $H_1$ and $P_1$ denote respectively the hyperbolic and parabolic elements given in this chart by

$$H_1 \cdot t := \frac{2t}{\lambda_1} \qquad \text{and} \qquad P_1 \cdot t := t - \frac{2m_1}{\lambda_1}.$$

Since the end-points of $(P_1 H_1) \cdot \widetilde{3I}$ are $-1$ and $1$ in the affine chart, the end-points of $(P_1 H_1 E_1) \cdot (3I)$ are $-\pi/4$ and $\pi/4$ in $\textbf{RP}^1$. This interval therefore has length $\pi/2$ and is centered on $0$, as desired. $H_2$ and $P_2$ are constructed in a similar manner, and this completes the proof. $\square$

For all $I$, let

$$(5.7) \qquad T_I := (T_1, T_2) := (P_1 H_1 E_1, P_2 H_2 E_2)$$

denote the canonical transformation constructed above, and define

$$\begin{aligned} y_I &:= T_1 \cdot x_I \,, \\ \psi_I &:= T_2 \circ \varphi \circ T_1^{-1} \,, \text{ and} \\ g_{\pm,I} &:= T_2 \circ f_{\pm,I} \circ T_1^{-1} \,. \end{aligned}$$

(5.8)

Let $d$ be a riemannian distance function over $\mathbf{PSL}(2,\mathbf{R})$. For example, we could use the distance function of some left-invariant metric. The main result of this subsection is the following estimate for $g_{\pm,I}$.

**Lemma 5.3.2.** *For all $I$,*

$$d(g_{\pm,I}, \mathrm{Id}) \lesssim \varepsilon_\varphi(I) \,.$$

To prove Lemma 5.3.2, we work in the rotated Penrose chart of $\mathbf{Ein}^{1,1}$ given in Section 2.3.3. Let $\tilde{y}_I$ denote the image of $y_I$ in this chart, and let $\tilde{\psi}$ and $\tilde{g}_{\pm,I}$ denote the respective images of $\psi$ and $g_{\pm,I}$ in this chart, that is

$$\tilde{\psi} := \tan \circ \psi \circ \arctan \qquad \text{and} \qquad \tilde{g}_{\pm,I} := \tan \circ g_{\pm,I} \circ \arctan \,.$$

We will require the following elementary result.

**Lemma 5.3.3.** *Choose $\delta < 1$. For $h \in \mathbf{PSL}(2,\mathbf{R})$, denoting by $\tilde{h}$ its realization in the affine chart of $\mathbf{RP}^1$, if*

(5.9) $$\left|\tilde{h}(0)\right| + \left|\tilde{h}'(0) - 1\right| + \left|\tilde{h}''(0)\right| < \delta \,,$$

*then*

$$d(h, \mathrm{Id}) \lesssim \left|\tilde{h}(0)\right| + \left|\tilde{h}'(0) - 1\right| + \left|\tilde{h}''(0)\right| \,.$$

*Proof.* Note that, for all $\delta < 1$, the set of all $h \in \mathbf{PSL}(2,\mathbf{R})$ satisfying (5.9) is compact. It is thus sufficient to prove the result in a neighbourhood of the identity. Consider now the function $\Phi$ defined near the identity in $\mathbf{PSL}(2,\mathbf{R})$ by

$$\Phi(h) := \left(\tilde{h}(0), \tilde{h}'(0), \tilde{h}''(0)\right) \,.$$

It suffices to show that $\Phi$ defines a local diffeomorphism near Id. By the inverse function theorem, this reduces to showing that $D\Phi(\mathrm{Id})$ is surjective. For this, it suffices in turn to construct a triplet $(\tilde{h}_{i,s})_{1 \leq i \leq 3, s \in ]-\delta,\delta]}$ of smooth families of fractional linear maps passing through Id whose derivatives with respect to $s$ at $s = 0$ are mapped by $D\Phi(\mathrm{Id})$ to a basis of $\mathbb{R}^3$. A straightforward calculation shows that the triplet

$$\tilde{h}_{1,s}(t) := t + s \,, \qquad \tilde{h}_{2,s}(t) := (1+s)t \qquad \text{and} \qquad \tilde{h}_{3,s}(t) = \frac{t}{1 - ts/2}$$

has the desired properties, and this completes the proof. □

*Proof of Lemma 5.3.2.* Indeed, let $I$ be a dyadic interval in $\mathcal{D}$. Since elliptic transformations are isometries, we may suppose that $E_1 = E_2 = \mathrm{Id}$. For brevity, we denote

$$\varepsilon := \varepsilon_\varphi(I) \,, \ f := f_{\pm,I} \,, \ \tilde{f} := \tilde{f}_{\pm,I} \,, \ \tilde{g} := \tilde{g}_{\pm,I} \,, \ L := \ell(3I) \,, \text{ and } L_\varphi := \ell(\varphi(3I)) \,.$$

Note that, by uniform continuity of $\varphi$, $L_\varphi$ tends to zero as $L$ tends to zero. We also denote $\tau := \tan$ and $\alpha := \arctan$, and we recall that the Taylor series of these functions are

$$\tau(x) = x \cdot \left(1 + \mathrm{O}(x^2)\right) \qquad \text{and} \qquad \alpha(t) = t \cdot \left(1 + \mathrm{O}(t^2)\right) \,.$$

By definition

$$\tilde{f} = \tau \circ f \circ \alpha \,.$$

Hence, by the chain rule,
$$\tilde{f}(0) = \tau(f(\alpha(0))),$$
$$\tilde{f}'(0) = \tau'(f(\alpha(0))) \cdot f'(\alpha(0)) \cdot \alpha'(0), \text{ and}$$
$$\tilde{f}''(0) = \tau''(f(\alpha(0))) \cdot f'(\alpha(0))^2 \cdot \alpha'(0)^2 + \tau'(f(\alpha(0))) \cdot f''(\alpha(0)) \cdot \alpha'(0)^2$$
$$+ \tau'(f(\alpha(0))) \cdot f'(\alpha(0)) \cdot \alpha''(0).$$

Applying the Taylor series of $\alpha$ simplifies these identities to
$$\tilde{f}(0) = \tau(f(0)),$$
$$\tilde{f}'(0) = \tau'(f(0)) \cdot f'(0), \text{ and}$$
$$\tilde{f}''(0) = \tau''(f(0)) \cdot f'(0)^2 + \tau'(f(0)) \cdot f''(0).$$

Applying the Taylor series of $\tau$ to the first identity and bearing in mind (4.2) yields

(5.10) $$\tilde{f}(0) = O(|f(0)|) = O(\varepsilon L_\varphi).$$

Likewise, the second and third identities yield

(5.11) $$\tilde{f}'(0) = \left(1 + O(|f(0)|^2)\right) \cdot f'(0) = \left(1 + O(\varepsilon L_\varphi)^2\right) \cdot \left(1 + O(\varepsilon)\right) \cdot \frac{L_\varphi}{L} = \left(1 + O(\varepsilon)\right) \cdot \frac{L_\varphi}{L},$$

and

(5.12) $$\tilde{f}''(0) = O\left(|f(0)| \cdot |f'(0)|^2 + |f''(0)|\right) = O\left((\varepsilon L_\varphi) \cdot \left(\frac{L_\varphi}{L}\right)^2 + \frac{\varepsilon L_\varphi}{L^2}\right) = O\left(\frac{\varepsilon L_\varphi}{L^2}\right).$$

Let $\widetilde{3I}$, $\widetilde{\varphi(3I)}$, $\lambda_1$, $\lambda_2$, $m_1$, and $m_2$ be as in the construction of $H_1$, $H_2$, $P_1$, and $P_2$ given in the proof of Lemma 5.3.1. By definition,
(5.13)
$$\lambda_1 = \ell(\widetilde{3I}), \qquad \lambda_2 = \ell(\widetilde{\varphi(3I)}), \qquad \tilde{y}_I = -\frac{2m_1}{\lambda_1}, \quad \text{and} \quad \tilde{g}(t) = \frac{2}{\lambda_2}\tilde{f}\left(\frac{\lambda_1 t}{2} + a_1\right) - \frac{2m_2}{\lambda_2}.$$

Applying the Taylor series of $\alpha$ a second time therefore yields
$$\lambda_1 = \ell(\widetilde{3I}) = \left(1 + O(L^2)\right) \cdot L \quad \text{and} \quad \lambda_2 = \ell(\widetilde{\varphi(3I)}) = \left(1 + O(L_\varphi^2)\right) \cdot L_\varphi,$$
so that
$$\frac{L}{\lambda_1}, \frac{L_\varphi}{\lambda_2} \simeq 1.$$

This reduces (5.10), (5.11), and (5.12) respectively to
$$\left|\tilde{f}(0)\right| \lesssim \varepsilon \lambda_2, \qquad \left|\tilde{f}'(0) - \frac{\lambda_2}{\lambda_1}\right| \lesssim \frac{\varepsilon \lambda_2}{\lambda_1}, \qquad \text{and} \qquad \left|\tilde{f}''(0)\right| \lesssim \frac{\varepsilon \lambda_2}{\lambda_1^2},$$

so that, by (5.13),
$$\left|\tilde{g}(\tilde{y}_I) + \frac{2a_2}{\lambda_2}\right| \lesssim \varepsilon, \qquad \left|\tilde{g}'(\tilde{y}_I) - 1\right| \lesssim \varepsilon, \qquad \text{and} \qquad \left|\tilde{g}''(\tilde{y}_I)\right| \lesssim \varepsilon.$$

To conclude, it remains only to show that

(5.14) $$\frac{2m_2}{\lambda_2} = -\tilde{y}_I + O(\varepsilon).$$

However, by Lemma 5.3.3, the fractional linear transformation
$$\tilde{g}' := \tilde{g} + \left(\frac{2a_2}{\lambda_2} + \tilde{y}_I\right)$$
lies in a $(C\varepsilon)$-neighborhood of the identity in $\mathbf{PSL}(2, \mathbf{R})$, for some suitable $C$. Hence,
$$\tilde{g}(1) + \left(\frac{2m_2}{\lambda_2} + \tilde{y}_I\right) = 1 + O(\varepsilon).$$

However, by definition,
$$\tilde{g}_{-,I}(1) \leq \tilde{\psi}_I(1) \leq \tilde{g}_{+,I}(1),$$
and
$$\tilde{\psi}_I(1) = 1.$$
Hence
$$\left(\frac{2m_2}{\lambda_2} + \tilde{y}_I\right) - O(\varepsilon) \leq \tilde{g}_{-,I}(1) - 1 \leq 0 \leq \tilde{g}_{+,I}(1) - 1 \leq \left(\frac{2m_2}{\lambda_2} + \tilde{y}_I\right) + O(\varepsilon),$$
and (5.14) follows. The result follows by Lemma 5.3.3 again, and this completes the proof. □

We will also require the following technical result. For any interval $I$, we denote respectively by $I_-$ and $I_+$ its lower and upper end-points.

**Lemma 5.3.4.** *Fix $\lambda > 0$, and let $I$ and $J$ be dyadic intervals such that*
$$J \subseteq \lambda I \quad \text{and} \quad \ell(J) = \frac{1}{2}\ell(I).$$
*Let $T := T_I := (T_1, T_2)$ denote the canonical transformation of $I$. There exists an integer $k \in [-\lambda/2, \lambda/2]$ such that*
$$(T_1 \cdot J)_\pm = \arctan\left(\frac{2k+1}{6} \pm \frac{1}{6}\right) + O\bigl(\ell(I)^2\bigr), \text{ and}$$
$$(T_2 \cdot \varphi(J))_\pm = \arctan\left(\frac{2k+1}{6} \pm \frac{1}{6}\right) + O\bigl(\ell(I)^2 + \varepsilon_\varphi(I)\bigr).$$

*Remark* 5.3.5. In the next section, we will be concerned with sets of the form
$$D_I := \bigl(3I \times \varphi(3I)\bigr) \setminus \bigcup_{\substack{J \subseteq 3I \\ \ell(J)=\ell(I)/2}} \bigl(3J \times \varphi(3J)\bigr),$$
where $I$ is a dyadic interval in $\mathcal{D}$. The purpose of Lemma 5.3.4 is to show that, as $\ell(I)$ tends to zero, up to a canonical transformation, the set $D_I$ Hausdorff converges to the set given in Figure 17.

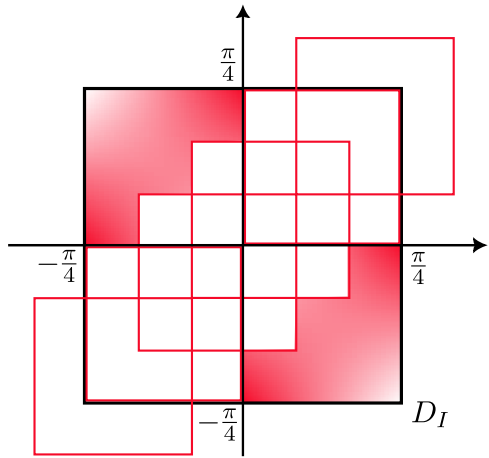

Figure 17. The limiting domain.

*Proof.* For $i \in \{1, 2\}$, let $P_i$, $H_i$, and $E_i$ denote respectively the parabolic, hyperbolic, and elliptic elements used in the construction of $T_i$ (Lemma 5.3.1). For each $i$, let $\tilde{P}_i$ and $\tilde{H}_i$ denote respectively the actions of $P_i$ and $H_i$ in the affine chart of $\mathbf{RP}^1$.

Let $\tilde{I}$ and $\tilde{J}$ denote respectively the images of $E_1 \cdot I$ and $E_2 \cdot J$ in the affine chart of $\mathbf{RP}^1$. Since the Taylor series of tan is

$$\tan(x) = x \cdot \left(1 + \mathrm{O}(x^2)\right),$$

we obtain

$$\ell(\tilde{I}) = \ell(I) \cdot \left(1 + \mathrm{O}(\ell(I)^2)\right),$$

$$\ell(\tilde{J}) = \frac{1}{2} \cdot \ell(I) \cdot \left(1 + \mathrm{O}(\ell(I)^2)\right), \text{ and}$$

$$\tilde{J}_- = \tilde{I}_- + \frac{k}{2} \cdot \ell(I) \cdot \left(1 + \mathrm{O}(\ell(I)^2)\right),$$

for some integer $k \in [-\lambda/2, \lambda/2]$. However, by definition,

$$(\tilde{P}_1 \tilde{H}_1) \cdot \tilde{I}_- = -\frac{1}{3} + \mathrm{O}\left(\ell(I)^2\right), \text{ and}$$

$$\ell\left((\tilde{P}_1 \tilde{H}_1) \cdot \tilde{I}\right) = \frac{2}{3} + \mathrm{O}\left(\ell(I)^2\right).$$

Since, in the affine chart, $\tilde{P}_1 \tilde{H}_1$ constitutes a rescaling followed by a translation, it follows that

$$(\tilde{P}_1 \tilde{H}_1) \cdot \tilde{J}_- = \frac{k-1}{3} + \mathrm{O}\left(\ell(I)^2\right), \text{ and}$$

$$\ell\left((\tilde{P}_1 \tilde{H}_1) \cdot \tilde{J}\right) = \frac{1}{3} + \mathrm{O}\left(\ell(I)^2\right),$$

and the estimates for $(T_1 \cdot J)_\pm$ follow.

Let $\widetilde{\varphi(J)}$ denote the image of $E_2 \cdot \varphi(J)$ in the affine chart of $\mathbf{RP}^1$. Let $\tilde{g}_\pm$ denote the image of $g_{\pm,I}$ in the affine chart, that is

$$\tilde{g}_I = \tan \circ g_{\pm,I} \circ \arctan.$$

By Lemma 5.3.2, for all $t$,

$$\tilde{g}_\pm(t) = t + \mathrm{O}\left(\varepsilon_\varphi(I)\right).$$

Since

$$\tilde{g}_-(\tilde{J}_-) \leqslant \widetilde{\varphi(J)}_- \leqslant \tilde{g}_+(\tilde{J}_+),$$

it follows that

$$\widetilde{\varphi(J)}_- = \frac{k-1}{3} + \mathrm{O}\left(\ell(I)^2 + \varepsilon_\varphi(I)\right).$$

A similar estimate holds for $\widetilde{\varphi(J)}_+$, and the estimates for $(T_2 \cdot \varphi(J))_\pm$ follow. This completes the proof. $\square$

**5.4. Area and curvature bounds.** We continue to use the notation of the previous sections. Theorem 5.1.1 will follow from the following estimates. For every dyadic interval $I$, let $K_I$ be as in (5.6).

**Lemma 5.4.1.** *For every dyadic interval $I \in \mathcal{D}_{\geqslant 2}$,*

$$\mathrm{Area}(K_I) \lesssim 1 \quad \text{and} \quad \sup_{x \in K_I} \|\mathrm{I\!I}_\varphi(x)\| \lesssim \varepsilon_\varphi(I).$$

We continue to work in the kleinian chart of $\mathbf{AdS}^{2,1}$. For all $I$, let $T_I$ be as (5.7), denote

$$S_I := T_I \cdot \Sigma_\varphi,$$

and let $\tilde{S}_I$ denote the image of $S_I$ in the kleinian chart. For $C > 0$, and for all $I$, consider the set $\tilde{X}_{I,C} \subseteq \mathbf{AdS}^{2,1}$ given by

$$\tilde{X}_{I,C} := \{(y_1, y_2, y_3) \mid y_1^2 + y_2^2 < 1 + y_3^2, \ |y_3| < C\varepsilon_\varphi(I)\},$$

as in Figure 18. Note that the two boundary components of $\tilde{X}_I$ in $\mathbf{AdS}^{2,1}$ are complete, totally geodesic planes. The constant $C$ is determined by the following result.

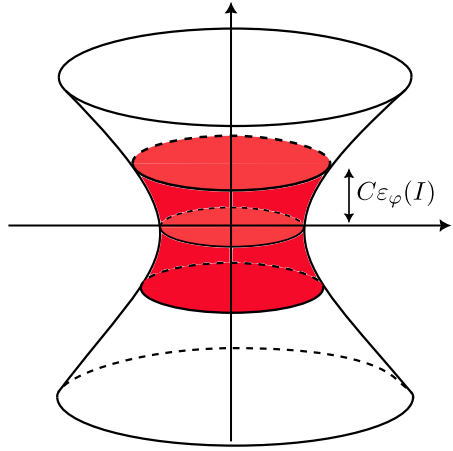

FIGURE 18. $\tilde{S}_I$ is pinched between two complete totally geodesic planes.

**Lemma 5.4.2.** *There exists $C > 0$ such that, for every dyadic interval $I$,*
$$\tilde{S}_I \subseteq \tilde{X}_{I,C} \,.$$

*Proof.* Indeed, denote
$$G_0 := \{(x,x) \mid x \in \mathbf{RP}^1\} \,,$$
and let $\tilde{G}_0$ denote its image in the kleinian chart. Thus
$$\tilde{G}_0 = \{(y_1, y_2, 0) \mid y_1^2 + y_2^2 = 1\} \,.$$

Let $G_\pm$ denote the graph of $g_{\pm,I}$ and let $\tilde{G}_\pm$ denote its image in the kleinian chart. Note that
$$G_\pm = \{(x, g_{\pm,I}(x)) \mid x \in \mathbf{RP}^1\} = M_\pm \cdot G_0 \,.$$
where $M_\pm \in \mathbf{PSL}(2,\mathbf{R}) \times \mathbf{PSL}(2,\mathbf{R})$ is given by
$$M_\pm := (\mathrm{Id}, g_{\pm,I}) \,.$$
Thus, denoting by $\tilde{M}$ the action of $M$ in the kleinian chart,
$$\tilde{G}_\pm = \tilde{M} \cdot \tilde{G}_0 \,.$$
By Lemma 5.3.2, there exists $C \geqslant 0$ such that
$$\sup_{(y_1,y_2,y_3) \in \tilde{G}_\pm} |y_3| \leqslant C\varepsilon_\varphi(I) \,.$$

Let $\Psi$ denote the graph of $\psi_I$, let $\tilde{\Psi}$ denote its image in the kleinian chart, and let $\mathrm{Conv}(\tilde{\Psi})$ denote its convex hull. Since $\tilde{\Psi}$ lies between $\tilde{G}_-$ and $\tilde{G}_+$,
$$\sup_{(y_1,y_2,y_3) \in \tilde{\Psi}} |y_3| \leqslant C\varepsilon_\varphi(I) \,,$$
so that
$$\tilde{\Psi} \subseteq \tilde{X}_{I,C} \,.$$
Since $\tilde{X}_{I,C}$ is the region bounded by two totally geodesic planes, it follows that
$$\mathrm{Conv}(\tilde{\Psi}) \subseteq \tilde{X}_{I,C} \,.$$
Since $S_I$ is the unique complete maximal spacelike surface bounded by $\psi_I$, by Proposition 4.6 of [BS10],
$$\tilde{S}_\psi \subseteq \mathrm{Conv}(\tilde{\Psi}) \subseteq \tilde{X}_{I,C} \,,$$
and this completes the proof. $\square$

For $r > 0$, let $B_r$ denote the open lorentzian ball of radius $r$ about $0$ in $\mathbf{R}^{2,1}$, that is
$$B_r := \{(y_1, y_2, y_3) \mid y_1^2 + y_2^2 < r^2 + y_3^2\} \,.$$

**Lemma 5.4.3.** *There exists $r < 1$ such that, for all $I$,*
$$T_I \cdot K_I \subseteq B_r \ .$$

*Proof.* Let $\tilde{D}_I$ denote the image of $T_I \cdot D_I$ in the kleinian chart. By definition,
$$T_I \cdot K_I \subseteq \tilde{D}_I \ .$$
By construction, $\partial_\infty(T_I \cdot D_I)$ is the square $[-\pi/2, \pi/2]^2$. In the kleinian chart, the points $p_1 := (-\pi/2, \pi/2)$ and $p_2 := (\pi/2, -\pi/2)$ are mapped to the ideal points corresponding to the lightlike rays
$$L_1 := \{(t, 0, t) \mid t \in \mathbf{R}\} \quad \text{and} \quad L_2 := \{(t, 0, -t) \mid t \in \mathbf{R}\} \ .$$
$\tilde{D}_I$ is thus the diamond of these two rays. Hence
$$\tilde{D}_I := \{(y_1, y_2, y_3) \mid y_1 \geqslant 0, \ -y_1 \leqslant y_3 \leqslant y_1, \ y_1^2 + y_2^2 < 1 + y_3^2\} \ .$$
as in Figure 16.

For every dyadic interval $J \subset I$, let $\tilde{D}_J$ denote the image of $T_I \cdot D_J$ in the kleinian chart. Consider now the set
$$\tilde{Y}_I := \tilde{X}_{I,C} \cap \left( \tilde{D}_I \setminus \bigcup_{\substack{J \subseteq 3I \\ \ell(J) = \ell(I)/2}} \tilde{D}_J \right) .$$

By Lemma 5.4.2,
$$T_I \cdot K_I \subseteq \tilde{Y}_I \ .$$

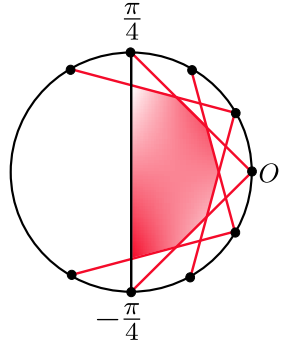

Figure 19. Another limiting domain.

We now determine the Hausdorff limit of $Y_I$ as $\ell(I)$ tends to zero. To this end, consider the six diamonds
$$D_{\infty,i} := D^{\text{ads}}\big((x_{i-1}, x_{i+2}), (x_{i+2}, x_{i-1})\big)$$
where, for each $j$, $x_j := \arctan(i/3)$, and $i$ varies over the index set $\mathcal{I} := \{-3, -2, \cdots, 2\}$. For all $i$, let $\tilde{D}_{\infty,i}$ denote the image of $D_{\infty,i}$ in the kleinian chart. By Lemma 5.3.4, $\tilde{D}_J$ subconverges towards $\tilde{D}_{\infty,i}$ for some $i$ in $\mathcal{I}$ as $\ell(I)$ tends to zero. Thus, as $\ell(I)$ tends to zero, $\tilde{Y}_I$ converges in the Hausdorff sense towards
$$\tilde{Y}_\infty := H_0 \cap \left( \tilde{D}_I \setminus \bigcup_{k \in S} \tilde{D}_{\infty,k} \right),$$
where
$$\tilde{H}_0 := \{(y_1, y_2, 0) \mid y_1^2 + y_2^2 < 1\} \ .$$
This set is illustrated in Figure 19. Since it is a compact subset of $\mathbf{AdS}^{2,1}$, there exists $r > 0$ such that
$$\tilde{Y}_\infty \subseteq B_r \ .$$

Thus, for all but finitely many $I$, $T_I \cdot K_I \subseteq B_r$, and the result follows upon increasing $r$ if necessary. □

We are now ready to prove Lemma 5.4.1.

*Proof of Lemma 5.4.1.* Choose a dyadic interval $I \in \mathcal{D}$. Let $H_0$ denote the totally geodesic surface $H_0$ given in this chart by
$$H_0 := \{(y_1, y_2, 0) \mid y_1^2 + y_2^2 < 1\}$$
Let $r < 1$ be as in Lemma 5.4.2, choose $r < r' < 1$, and let $B$ and $B'$ denote respectively the balls of radius $r$ and $r'$ about the origin in $H_0$. For $I$ sufficiently small, by Lemma 5.4.2, there is a smooth function $u_I : B' \to \mathbf{R}$ whose graph over this ball is a portion of $\tilde{S}_I$ and, furthermore,
$$\|u_I\|_{C^0(B')} \lesssim \epsilon_\varphi(I) \,.$$
It follows by Proposition 4.A of [Sep19] that, for all such $I$,
$$\|u_I\|_{C^2(B)} \lesssim \epsilon_\varphi(I) \,.$$
By Lemma 5.4.3, for all sufficiently small $I$, the graph of $u_I$ contains $K_I$. The area and curvature bounds of $K_I$ immediately follow, and this completes the proof. □

We are now ready to prove Theorem 5.1.1.

*Proof of Theorem 5.1.1.* Indeed,
$$\Sigma_\varphi = \Sigma_0 \cup \Sigma_\infty \,,$$
where $\Sigma_0$ is compact by Lemma 5.2.1, and
$$\Sigma_\infty := \bigcup_{I \in \mathcal{D}_{\geqslant 2}} K_I \,.$$
Hence
$$\int_{\Sigma_\varphi} \|\mathrm{I\!I}_\varphi\|^2 \, d\,\mathrm{Area} \leqslant \int_{\Sigma_0} \|\mathrm{I\!I}_\varphi\|^2 \, d\,\mathrm{Area} + \sum_{I \in \mathcal{D}_{\geqslant 2}} \int_{K_I} \|\mathrm{I\!I}_\varphi\|^2 \, d\,\mathrm{Area}$$
$$\leqslant \int_{\Sigma_0} \|\mathrm{I\!I}_\varphi\|^2 \, d\,\mathrm{Area} + \sum_{I \in \mathcal{D}_{\geqslant 2}} \mathrm{Area}(K_I) \cdot \sup_{x \in K_I} \|\mathrm{I\!I}_\varphi(x)\|^2 \,.$$
Thus, by Lemma 5.4.1,
$$\int_{\Sigma_\varphi} \|\mathrm{I\!I}_\varphi\|^2 \, d\,\mathrm{Area} \lesssim \int_{\Sigma_0} \|\mathrm{I\!I}_\varphi\|^2 \, d\,\mathrm{Area} + \sum_{I \in \mathcal{D}_{\geqslant 2}} \varepsilon_\varphi(I)^2 < \infty \,,$$
as desired. □

## 6. Minimal lagrangian diffeomorphisms

We now complete the chain of equivalences required to prove the main results of this paper. It remains only to show that a homeomorphism $\varphi : \mathbf{RP}^1 \to \mathbf{RP}^1$ bounds a complete maximal spacelike surface of finite renormalized area if and only if it bounds a minimal lagrangian diffeomorphism with square integrable Beltrami differential. It then follows by the characterization [She18] of Shen that in this case $\varphi$ is Weil–Petersson. In fact, we prove a result slightly stronger than that of [She18]. Indeed, Shen's result states that $\varphi$ is Weil–Petersson whenever its minimal lagrangian diffeomorphism has square integrable Beltrami differential *and* is quasiconformal. Using a recent result [Mor25] of Moriani, we will see that the hypothesis of quasiconformality is unnecessary. More precisely, for minimal lagrangian diffeomorphisms bounded by homeomorphisms of $\mathbf{RP}^1$, square integrability of the Beltrami differential implies quasiconformality.

6.1. **Minimal lagrangian diffeomorphisms.** Let $\varphi : \mathbf{RP}^1 \to \mathbf{RP}^1$ be a homeomorphism. Let $\Sigma_\varphi$ denote the unique complete maximal spacelike surface in $\mathbf{AdS}^{2,1}$ that it bounds, let $\Phi_\varphi : \mathbf{H}^2 \to \mathbf{H}^2$ denote the unique minimal lagrangian diffeomorphism that it bounds, and let $\mu_\varphi$ denote the Beltrami differential of $\Phi_\varphi$.

**Theorem 6.1.1.** *Let $\varphi : \mathbf{RP}^1 \to \mathbf{RP}^1$ be a homeomorphism. The surface $\Sigma_\varphi$ has finite renormalized area if and only if $\mu_\varphi \in L^2(\mathbf{H})$.*

6.2. **Algebraic properties of $\mathfrak{sl}(2,\mathbf{R})$.** Theorem 6.1.1 will be proven using the matrix model of $\mathbf{AdS}^{2,1}$, described in Section 2.1.2. We also refer the reader to Section 6.1 of [Bar18] for a deeper discussion of the formalism that will be used in the rest of this section (see also Section 6.3 of [BS20]).

Recall that $\mathsf{SL}(2,\mathbf{R})$ is the double cover of $\mathbf{AdS}^{2,1}$. Its tangent space at the identity is the Lie algebra $\mathfrak{sl}(2,\mathbf{R})$ of trace-free $2 \times 2$ real matrices. The restriction of $\langle \cdot, \cdot \rangle_{\mathrm{mat}}$ to this space is non-degenerate and of signature $(2,1)$. Let $H \subseteq \mathfrak{sl}(2,\mathbf{R})$ denote the quadric

$$H := \{\mathbf{M} \in \mathfrak{sl}(2,\mathbf{R}) \mid \|\mathbf{M}\|^2_{\mathrm{mat}} = -1\} = \mathfrak{sl}(2,\mathbf{R}) \cap \mathsf{SL}(2,\mathbf{R}) \ .$$

This quadric consists of two connected components, each of which is isometric to 2-dimensional hyperbolic space $\mathbf{H}^2$. An elementary calculation shows that $H$ consists of all elements $\mathbf{X} \in \mathfrak{sl}(2,\mathbf{R})$ such that

$$(6.1) \qquad \mathbf{X}^2 = -\mathrm{Id} \ .$$

In particular, the standard complex structure $\mathbf{J}$ given by (2.2) is an element of this quadric. We denote by $H^+$ the component of $H$ containing $\mathbf{J}$, and we identify this component with $\mathbf{H}^2$.

$H$ is furnished with a natural orientation and complex structure as follows. For all $\mathbf{X}$ in $H$, let $\mathbf{m} := \mathbf{m}_{\mathbf{X},l}$ denote the operation of multiplication on the left by $\mathbf{X}$. Note that the tangent space $T_\mathbf{X} H$ of $H$ at any such $\mathbf{X}$ is the orthogonal complement of the subspace $\langle \mathbf{X}, \mathrm{Id} \rangle$ with respect to $\langle \cdot, \cdot \rangle_{\mathrm{mat}}$. Since $\mathbf{m}$ is an isometry, and since it preserves $\langle \mathbf{X}, \mathrm{Id} \rangle$, it also preserves $T_\mathbf{X} H$. Furthermore, by (6.1), $\mathbf{m}^2 = -\mathrm{Id}$, so that $\mathbf{m}$ defines a complex structure over this space. The orientation of $H$ is now defined such that $\mathbf{m}_{\mathbf{X},l}$ rotates by $\pi/2$ in the *positive* direction for every $\mathbf{X}$. The reader may then verify that, for all $\mathbf{X}$ in $H$, multiplication on the right by $\mathbf{X}$ rotates by $\pi/2$ in the negative direction.

This construction yields an identification of the space of future-oriented timelike geodesics in $\mathbf{AdS}^{2,1}$ with $\mathbf{H}^2 \times \mathbf{H}^2$ as follows. Let $\mathsf{U}_t \mathsf{SL}(2,\mathbf{R})$ denote the bundle of unit timelike tangent vectors of $\mathsf{SL}(2,\mathbf{R})$, viewed as the submanifold

$$\mathsf{U}_t\mathsf{SL}(2,\mathbf{R}) := \{(\mathbf{M},\mathbf{N}) \in \mathcal{M}(2,\mathbf{R}) \times \mathcal{M}(2,\mathbf{R}) \mid \|\mathbf{M}\|^2_{\mathrm{mat}} = \|\mathbf{N}\|^2_{\mathrm{mat}} = -1 \ , \text{ and } \langle \mathbf{M}, \mathbf{N} \rangle_{\mathrm{mat}} = 0\} \ .$$

The *left* and *right* projections $\pi_l, \pi_r : \mathsf{U}_t\mathsf{SL}(2,\mathbf{R}) \to H^+ = \mathbf{H}^2$ are defined respectively by

$$\pi_l(\mathbf{M},\mathbf{N}) := \mathbf{M}^{-1} \cdot \mathbf{N} \qquad \text{and} \qquad \pi_r(\mathbf{M},\mathbf{N}) := \mathbf{N} \cdot \mathbf{M}^{-1} \ .$$

The function $\Pi := (\pi_l, \pi_r)$ is invariant under the geodesic flow (see, for example, [Bar18]), and thus projects to a diffeomorphism from the space of future-oriented timelike geodesics in $\mathbf{AdS}^{2,1}$ into $\mathbf{H}^2 \times \mathbf{H}^2$, as desired.

6.3. **The left and right Gauss maps.** Let $\mathbf{M} : \Sigma_\varphi \to \mathsf{SL}(2,\mathbf{R})$ be a lift of the canonical embedding. Let $\mathbf{N} : \Sigma_\varphi \to \mathcal{M}(2,\mathbf{R})$ denote the future-oriented unit normal vector field over $\mathbf{M}$. The *left* and *right Gauss maps* of $\Sigma_\varphi$ are defined respectively by

$$(6.2) \qquad G_{\varphi,l} := \pi_l(\mathbf{M},\mathbf{N}) = \mathbf{M}^{-1} \cdot \mathbf{N} \qquad \text{and} \qquad G_{\varphi,r} := \pi_r(\mathbf{M},\mathbf{N}) = \mathbf{N} \cdot \mathbf{M}^{-1} \ .$$

Each of these functions maps $\Sigma_\varphi$ into $H^+ = \mathbf{H}^2$. The following result is proven in [Tre24b] (see also Appendix A).

**Lemma 6.3.1.** *Let $\varphi : \mathbf{RP}^1 \to \mathbf{RP}^1$ be a homeomorphism. The functions $G_{\varphi,l}$ and $G_{\varphi,r}$ map $\Sigma_\varphi$ diffeomorphically onto the whole of $\mathbf{H}^2$. Furthermore, the unique minimal lagrangian diffeomorphism extending $\varphi$ is given by*

$$\Phi_\varphi := G_{\varphi,r} \circ G_{\varphi,l}^{-1}.$$

Let $A := A_\varphi$ denote the shape operator of $\Sigma_\varphi$, and let $\lambda := \lambda_\varphi$ denote its positive eigenvalue. As above, let $\mu := \mu_\varphi$ denote the Beltrami differential of $\Phi_\varphi$, and define the function $\tilde{\mu} := \tilde{\mu}_\varphi$ over $\Sigma_\varphi$ such that

$$\tilde{\mu}_\varphi^2 \, d\,\text{Area}_\Sigma = G_{\varphi,l}^*\left(\|\mu\|^2 d\,\text{Area}_{\mathbf{H}^2}\right),$$

where $d\,\text{Area}_\Sigma$ and $d\,\text{Area}_{\mathbf{H}^2}$ denote respectively the area forms of $\Sigma_\varphi$ and $\mathbf{H}^2$.

Let $g_{\varphi,l}$ and $g_{\varphi,r}$ denote the respective pull-backs through $G_{\varphi,l}$ and $G_{\varphi,r}$ of the hyperbolic metric of $\mathbf{H}^2$. We recall the following result (see [Bar18] and [BS20]).

**Lemma 6.3.2.** *The metrics $g_{\varphi,l}$ and $g_{\varphi,r}$ are given by*

$$g_{\varphi,l} = g\big((A_\varphi + J)\cdot,(A_\varphi + J)\cdot\big) \quad \text{and} \quad g_{\varphi,r} = g\big((A_\varphi - J)\cdot,(A_\varphi - J)\cdot\big),$$

*where $J$ denotes the complex structure of $\Sigma_\varphi$ compatible with the orientation.*

**Corollary 6.3.3.** *The Beltrami differential of $\Phi_\varphi$ is related to $\lambda_\varphi$ by*

$$(6.3) \quad \|\mu_\varphi \circ \Psi_l\|^2 = \frac{4\lambda_\varphi^2}{\left(1+\lambda_\varphi^2\right)^2} \quad \text{and} \quad \tilde{\mu}_\varphi^2 = \frac{4\lambda_\varphi^2 \cdot \left(1-\lambda_\varphi^2\right)}{\left(1+\lambda_\varphi^2\right)^2}.$$

*Proof of Lemma 6.3.2.* Indeed, for brevity, denote $G := G_{\varphi,l}$. Choose $p \in \Sigma_\varphi$, and let $\xi$ be a tangent vector of $\Sigma_\varphi$ at this point. Differentiating (6.2) yields

$$D G(p) \cdot \xi = -\left(\mathbf{M}^{-1} \cdot (D\mathbf{M}(p)\cdot\xi)\right)\cdot\left(\mathbf{M}^{-1}\cdot\mathbf{N}\right) + \mathbf{M}^{-1}\cdot(D\mathbf{N}(p)\cdot\xi).$$

Let $J'$ denote the complex structure of $H$ compatible with the orientation, as described in Section 6.2. By the discussion of the preceeding section, on the one hand, and by definition of $A$, on the other,

$$DG(p)\cdot\xi = J'\cdot\mathbf{M}^{-1}\cdot D\mathbf{M}(p)\cdot\xi + \mathbf{M}^{-1}\cdot D\mathbf{M}(p)\cdot A(p)\cdot\xi.$$

Since multiplication on the left by $\mathbf{M}^{-1}$ is an orientation-preserving isometry,

$$J'\cdot\mathbf{M}^{-1}\cdot D\mathbf{M}(p) = \mathbf{M}^{-1}\cdot D\mathbf{M}(p)\cdot J.$$

It follows that

$$DG(p)\cdot\xi = \mathbf{M}^{-1}\cdot D\mathbf{M}(p)\cdot(J(p)+A(p))\cdot\xi.$$

Since $\mathbf{M}^{-1}\cdot D\mathbf{M}(p)$ is an isometry, the formula for $g_{\varphi,l}$ follows. The formula for $g_{\varphi,r}$ follows in a similar manner, and this completes the proof. □

### 6.4. Asymptotic behaviour of the shape operator.
We will require the following result of Ishihara.

**Theorem 6.4.1** (Ishihara [Ish88]). *Let $\Sigma$ be a complete maximal spacelike surface in $\mathbf{AdS}^{2,1}$. If $\lambda$ denotes the positive eigenvalue of its shape operator, then*

$$(6.4) \qquad\qquad\qquad \lambda \leqslant 1.$$

*Furthermore, if equality is attained at one point, then it is attained at all points.*

An elementary compactness argument extends Ishihara's result as follows.

**Lemma 6.4.2.** *Let $\Sigma$ be a complete maximal spacelike surface in $\mathbf{AdS}^{2,1}$, let $A$ denote its shape operator, and let $\nabla$ denote its Levi–Civita covariant derivative.*

(1) *For all $p \in \Sigma$, and for all $r > 0$,*

$$\text{Area}(B_r(p)) \geqslant \pi r^2.$$

(2) *for all $k$, there exists $C_k > 0$, which does not depend on $\Sigma$, such that*

$$\|\nabla^k A\|^2 \leqslant C_k .$$

*Proof.* Indeed, by (6.4) and Gauss' equation, $\Sigma$ has non-positive curvature, and (1) follows by Rauch's comparison theorem. To prove (2), suppose the contrary holds for some $k \in \mathbf{N}$. There exists a sequence $(\Sigma_m)_{m \in \mathbf{N}}$ of complete maximal surfaces in $\mathbf{AdS}^{2,1}$, and a sequence $(p_m)_{m \in \mathbf{N}}$ in $\mathbf{AdS}^{2,1}$ such that $p_m \in \Sigma_m$ for all $m$, and

$$\lim_{m \to \infty} \|\nabla^k A_m(p_m)\| \to +\infty ,$$

where, for all $m$, $A_m$ denotes the shape operator of $\Sigma_m$. Let $q_0$ be some fixed point of $\mathbf{AdS}^{2,1}$, and let $n_0$ be a unit timelike vector tangent to $\mathbf{AdS}^{2,1}$ at this point. For all $m$, let $f_m$ be an isometry of $\mathbf{AdS}^{2,1}$ such that $f_m(p_m) = q_0$ and $f_m(\Sigma_m)$ is normal to $n_0$ at this point. For all $m$, denote $\Sigma'_m := f_m(\Sigma_m)$, and let $A'_m$ denote the shape operator of this surface. By Theorem 5.1 of [SST23], we may suppose that $(\Sigma_m)_{m \in \mathbf{N}}$ converges smoothly to a complete maximal spacelike surface $\Sigma'_\infty$. Denoting by $A'_\infty$ the shape operator of this surface, we obtain

$$\|\nabla^k A'_\infty(p_\infty)\| = \lim_{m \to \infty} \|\nabla^k A'_m(q_0)\| = \lim_{m \to \infty} \|\nabla^k A_m(p_m)\| = +\infty .$$

This is absurd, and the result follows. $\square$

**Lemma 6.4.3.** *Let $\varphi : \mathbf{RP}^1 \to \mathbf{RP}^1$ be a homeomorphism. If $\|\mu_\varphi\|_{L^2} < \infty$, then $\lambda_\varphi(p)$ tends either to $0$ or to $1$ as $p$ tends to infinity in $\Sigma_\varphi$.*

*Proof.* We show that the product $\lambda_\varphi^2(p)\bigl(1 - \lambda_\varphi^2(p)\bigr)$ tends to $0$ as $p$ tends to infinity in $\Sigma_\varphi$. Indeed, suppose the contrary. Note first that, by (6.4), this product is non-negative over $\Sigma_\varphi$. There therefore exists $c > 0$, and a diverging sequence $(p_m)_{m \in \mathbf{N}}$ of points in $\Sigma_\varphi$, such that

$$\lim_{m \to \infty} \lambda_\varphi^2(p_m)\bigl(1 - \lambda_\varphi^2(p_m)\bigr) \geqslant c > 0 .$$

By Part (2) of Lemma 6.4.2 with $k = 1$, there exists a constant $L$ such that $\lambda_\varphi^2(1 - \lambda_\varphi^2)$ is $L$-Lipschitz over $\Sigma_\varphi$. For all $m$, let $B_m$ denote the ball of radius $c/2L$ about $p_m$ in $\Sigma$, and note that, over this ball

$$\lambda_\varphi^2(1 - \lambda_\varphi^2) \geqslant \frac{c}{2} .$$

By Part (1) of Lemma 6.4.2, for all $m$,

$$\mathrm{Area}(B_m) \geqslant \frac{\pi c^2}{4L^2} .$$

Hence

$$\int_{B_m} \lambda_\varphi^2\bigl(1 - \lambda_\varphi^2\bigr) \mathrm{d}\,\mathrm{Area}_{\Sigma_\varphi} \geqslant \frac{\pi c^3}{8L^2} .$$

Upon extracting a subsequence, we may suppose that, for all $m \neq n$,

$$B_m \cap B_n = \emptyset .$$

Thus, by (6.3),

$$\|\mu_\varphi\|_{L^2}^2 = \int_{\mathbf{H}^2} \|\mu_\varphi\|^2 \mathrm{d}\,\mathrm{Area}_{\mathbf{H}^2} = \int_{\Sigma_\varphi} \tilde{\mu}_\varphi^2 \mathrm{d}\,\mathrm{Area}_{\Sigma_\varphi}$$

$$= \int_{\Sigma_\varphi} \frac{4\lambda_\varphi^2\bigl(1 - \lambda_\varphi^2\bigr)}{\bigl(1 + \lambda_\varphi^2\bigr)^2} \mathrm{d}\,\mathrm{Area}_{\Sigma_\varphi} \geqslant \sum_{m=1}^{\infty} \int_{B_m} \frac{4\lambda_\varphi^2\bigl(1 - \lambda_\varphi^2\bigr)}{\bigl(1 + \lambda_\varphi^2\bigr)^2} \mathrm{d}\,\mathrm{Area}_{\Sigma_\varphi} = \infty .$$

This is absurd, and the assertion follows. Finally, it follows by a connectedness argument that $\lambda_\varphi(p)$ converges either to $0$ or to $1$ as $p$ tends to infinity in $\Sigma_\varphi$, and this completes the proof. $\square$

We now show that, in our context, quasiconformality of the minimal lagrangian diffeomorphism follows from the square integrability of the Beltrami differential.

**Lemma 6.4.4.** *Let $\varphi : \mathbf{RP}^1 \to \mathbf{RP}^1$ be a homeomorphism, let $\Phi_\varphi$ denote the unique minimal lagrangian diffeomorphism that it bounds, and let $\mu_\varphi$ denote the Beltrami differential of $\Phi_\varphi$. If*

$$\|\mu_\varphi\|_{L^2} < \infty ,$$

*then $\Phi_\varphi$ is quasiconformal, that is*

$$\|\mu_\varphi\|_{L^\infty} < 1 .$$

*Proof.* Suppose first that $\lambda_\varphi(p)$ tends to $0$ as $p$ tends to infinity in $\Sigma_\varphi$. In particular $\lambda_\varphi$ attains its maximum at some point $p_0$, say. By Theorem 6.4.1,

$$\lambda_\varphi(p_0) \leqslant 1 .$$

Furthermore, since equality holds in 6.4 at one point if and only if it holds at all points, it follows that this inequality is strict. Consequently,

$$\|\lambda_\varphi\|_{L^\infty} = \lambda_\varphi(p_0) < 1 ,$$

and it follows by (6.3) that $\|\mu_\varphi\|_{L^\infty} < 1$, as desired.

We now show that $\lambda_\varphi(p)$ tends to $0$ as $p$ tends to infinity in $\Sigma_\varphi$. Indeed, suppose the contrary. Then, by Lemma 6.4.3, $\lambda_\varphi(p)$ converges to $1$ as $p$ tends to infinity in $\Sigma_\varphi$. It follows by (6.3) that, outside of some compact subset $\Sigma_0$ of $\Sigma_\varphi$,

$$\tilde{\mu}_\varphi^2 \geqslant 1 - \lambda_\varphi^2 .$$

Let $K_{\text{int}}$ denote the intrinsic curvature of $\Sigma_\varphi$. By Gauss' equation

$$1 - \lambda_\varphi^2 = -K_{\text{int}} = |K_{\text{int}}| .$$

Thus,

$$\int_{\Sigma_\varphi \setminus \Sigma_0} |K_{\text{int}}| \, d\operatorname{Area}_{\Sigma_\varphi} \leqslant \int_{\Sigma_\varphi \setminus \Sigma_0} G^*_{\varphi,l}\bigl(\|\mu_\varphi\|^2 d\operatorname{Area}_{\mathbf{H}^2}\bigr) \leqslant \int_{\mathbf{H}^2} \|\mu_\varphi\|^2 d\operatorname{Area}_{\mathbf{H}^2} = \|\mu_\varphi\|_{L^2}^2 < \infty ,$$

and $\Sigma_\varphi$ therefore has finite total curvature. By Theorem 1.2.1 of [Mor25], the ideal boundary of $\Sigma$ is a polygon consisting of lightlike segments. Since no such polygon is the graph of a homeomorphism, this is absurd, and it follows that the limit is $0$, as asserted. $\square$

We are now ready to prove Theorem 6.1.1.

*Proof of Theorem 6.1.1.* Suppose first that $\Sigma$ has finite renormalized area. By (6.3) and (6.4),

$$G^*_{\varphi,l}\bigl(\|\mu_\varphi\|^2 d\operatorname{Area}_{\mathbf{H}^2}\bigr) \lesssim \|A_\varphi\|^2 d\operatorname{Area}_{\Sigma_\varphi} .$$

Thus

$$\|\mu_\varphi\|_{L^2}^2 = \int_{\mathbf{H}^2} \|\mu_\varphi\|^2 d\operatorname{Area}_{\mathbf{H}^2} = \int_{\Sigma_\varphi} G^*_{\varphi,l}\bigl(\|\mu_\varphi\|^2 d\operatorname{Area}_{\mathbf{H}^2}\bigr) \lesssim \int_{\Sigma_\varphi} \|A_\varphi\|^2 d\operatorname{Area}_{\Sigma_\varphi} < \infty ,$$

so that $\Phi_\varphi$ has $L^2$ Beltrami differential, as desired.

Suppose now that $\Phi_\varphi$ has $L^2$ Beltrami differential. By Lemma 6.4.4, $\Phi_\varphi$ is quasiconformal. It follows by (6.3) that

$$\lambda_\varphi^2 < 1 ,$$

so that, by (6.3) again,

$$\|A_\varphi\|^2 d\operatorname{Area}_{\Sigma_\varphi} \lesssim G^*_{\varphi,l}\bigl(\|\mu_\varphi\|^2 d\operatorname{Area}_{\mathbf{H}^2}\bigr) .$$

Hence

$$\int_{\Sigma_\varphi \setminus \Sigma_0} \|A_\varphi\|^2 d\operatorname{Area}_{\Sigma_\varphi} \lesssim \int_{\Sigma_\varphi \setminus \Sigma_0} G^*_{\varphi,l}\bigl(\|\mu_\varphi\|^2 d\operatorname{Area}_{\mathbf{H}^2}\bigr) \leqslant \int_{\mathbf{H}^2} \|\mu_\varphi\|^2 d\operatorname{Area}_{\mathbf{H}^2} = \|\mu_\varphi\|_{L^2}^2 < \infty ,$$

and $\Sigma_\varphi$ thus has finite renormalized area, as desired. $\square$

## Appendix A. Minimal Lagrangian extension

For the reader's convenience, we include in this appendix a proof of the following fact:

An orientation preserving homeomorphism $\varphi : \mathbf{RP}^1 \to \mathbf{RP}^1$ extends to a unique minimal lagrangian diffeomorphism $\Phi_\varphi : \mathbf{H}^2 \to \mathbf{H}^2$. Furthermore, $\Phi_\varphi$ arises from the Gauss map of the unique maximal spacelike surface $\Sigma_\varphi$ bounded by $\varphi$ in $\mathbf{AdS}^{2,1}$.

In the case where $\varphi$ is quasisymmetric, this result has been proven by Bonsante–Seppi in [BS18, Corollary 1.6]. The general case, explained below, is part of the doctoral thesis of the fifth author (see Theorem L and Chapters 10 and 16 of [Tre24b]).

### A.1. Convex surfaces.
An embedded surface $S$ in $\mathbf{AdS}^{2,1}$ is *future-convex* if, in a suitable affine chart, the convex hull of $S$ is contained in the future of $S$. The notion of *past-convex* surface is analogous.

*Remark* A.1.1. If $S$ is $C^2$ and properly embedded, then it is future-convex or past convex if and only if its shape operator is non-negative or non-positive definite, respectively.

A.1.1. *Gauss map.* As explained in Section 6.2 (see also [BS20, Subsection 3.5]), the space of future oriented timelike geodesics in $\mathbf{AdS}^{2,1}$ identifies with $\mathbf{H}^2 \times \mathbf{H}^2$.

Let $S$ be a properly embedded spacelike surface. By the above identification, the Gauss map $G_S$ of $S$ writes as follows. For every point $x \in S$,

$$G_S(x) = (p,q) \text{ for } L_{p,q} := \exp_x(\mathbf{R} N_x S)$$

where $N_x S$ is the future directed unit normal to $S$ at $x$.

At a point $x \in S$, the data of the normal $N_x S$ is equivalent to the data of the tangent plane $T_x S$. Studying the Gauss map through the grassmanian is effective for convex surfaces, since tangent planes coincides with support planes.

A.1.2. *Normal flow.* Let $S$ be a properly embedded spacelike surface, and denote by $N$ its future-directed unitary normal vector field. The *normal flow* of $S$ is defined by

$$F : S \times \mathbf{R} \to \mathbf{AdS}^{2,1}$$
$$(x,t) \mapsto \exp_x(tN(x)) \,.$$

What precedes ensures that the Gauss map is invariant along the normal flow.

It is known (see for example [Sep19]) that the normal flow of $\Sigma_\varphi$ is a diffeomorphism onto when restricted to $\Sigma_\varphi \times [-\pi/4, \pi/4]$. Denote by $S_\varphi^+$ and $S_\varphi^-$ the images of the normal flow at time $\pi/4$ and $-\pi/4$, respectively. It turns out that $S_\varphi^+$ and $S_\varphi^-$ are respectively a past-convex and a future-convex smooth spacelike surface with constant sectional curvature $K = -2$.

A.1.3. *Achronal loops.* $K$–surfaces in $\mathbf{AdS}^{2,1}$ have been classified in [BS18, Theorem 1.3]. We present below a weaker version of this result, after introducing a technical definition.

**Definition A.1.2.** A *past-directed sawtooth* is a horizontal segment followed by an adjacent vertical segment in $\mathbf{RP}^1 \times \mathbf{RP}^1$. A *future-directed sawtooth* is a vertical segment followed by an adjacent horizontal segment in $\mathbf{RP}^1 \times \mathbf{RP}^1$. The common point of the two segments constituting a sawtooth is called the *vertex* of the sawtooth.

*Remark* A.1.3. A (past or future directed) sawtooth defines a triangle in $\mathbf{AdS}^{2,1} \cup \mathbf{Ein}^{1,1}$, that is part of a totally geodesic plane in $\mathbf{AdS}^{2,1}$ of degenerate signature $(1,0)$.

**Theorem A.1.4.** [BS18, Theorem 1.3] *Let $\Lambda$ be an achronal loop. There exist a properly embedded achronal past-convex surface $S_\Lambda^+$ and a properly embedded achronal future-convex surface $S_\Lambda^-$ in $\mathbf{AdS}^{2,1}$ such that:*
- *$\partial_\infty S_\Lambda^\pm = \Lambda$;*
- *the lightlike part of $S_\Lambda^+$ or $S_\Lambda^-$ is a union of lightlike triangles associated to past-directed sawteeth or future- directed sawteeth, respectively;*

- *the spacelike part of $S_\Lambda^+$ or $S_\Lambda^-$ is a smooth spacelike surface with constant sectional curvature $K = -2$ that is past-convex or future-convex, respectively.*

Applying this result to $\Lambda = \partial_\infty \Sigma$, where $\Sigma$ is a maximal surface, we obtain the existence of $S_\Lambda^+$ and $S_\Lambda^-$.

When the asymptotic loop $\Lambda$ is not acausal, the normal flow is a diffeomorphism onto only restricted to $\Sigma \times (-\pi/4, \pi/4)$. Nonetheless, $F(\Sigma, t)$ converges in the Haussdorf topology to $S_\Lambda^+$ and $S_\Lambda^-$ when $t$ tends respectively to $\pi/4$ and $-\pi/4$.

The Gauss map of $\Sigma$ still coincides with the Gauss map of $S_\Lambda^\pm$. For this reason, we will abusively denote $G_\Lambda$ both Gauss maps of $\Sigma$ and $S_\Lambda^\pm$.

A.2. **Extension.** We recall some important results.

**Lemma A.2.1.** [BS18, Lemma 4.2] *Let $S$ be a properly embedded convex spacelike surface in $\mathbf{AdS}^{2,1}$, then the projections $G_{S,l}, G_{S,r}$ are injective.*

A.2.1. *Duality.* There is a duality beetween totally geodesic spacelike planes in $\mathbf{AdS}^{2,1}$ and points in $\mathbf{AdS}^{2,1}$: indeed, a point in $\mathbf{AdS}^{2,1}$ is a negative definite line in $\mathbf{R}^{2,2}$, and its orthogonal complement is a signature $(2,1)$ subspace, which intersects $\mathbf{AdS}^{2,1}$ in a totally geodesic spacelike plane. Similarly, there is a duality between totally geodesic degenerate planes of signature $(1,0)$ in $\mathbf{AdS}^{2,1}$ and points in $\mathbf{Ein}^{1,1}$.

Let $(p, Q)$ be a pointed support plane in $S_\Lambda^+$. Denote by $(q^*, P^*)$ the dual pointed plane, namely $q^*$ is the dual of $Q$ and $P^*$ is the dual of $p$. One can prove that $(q^*, P^*)$ is a of pointed support plane for $S_\Lambda^-$ (see [Tre24b, Proposition 10.2.1].

**Lemma A.2.2.** [BS10, Lemma 3.18] *Let $S$ be a properly embedded convex achronal surface in $\mathbf{AdS}^{2,1}$. Let $(p_k)_{k \in \mathbb{N}}$ be a sequence of points in $S$ converging to $p_\infty = (\xi, \eta) \in \partial_\infty S$.*

*Let $((p_k, P_k))_{k \in \mathbb{N}}$ be a sequence of pointed support spacelike planes for $S$ converging to $(p_\infty, P_\infty)$. If $P_\infty$ is spacelike or is dual to $p_\infty$, then $(G_S((p_k, P_k)))_{k \in \mathbb{N}}$ converges to $(\xi, \eta)$.*

*Remark* A.2.3. The original statement is different and less technical, since a weaker result was sufficent for the goal of Bonsante–Schlenker. However, the proof given in [BS10] shows the stronger aforementioned result.

**Corollary A.2.4.** [Tre24b, Corollary P] *Let $\Lambda$ be an achronal loop. Then, $G_{\Lambda,l}$ and $G_{\Lambda,r}$ extend continuously to homeomorphisms of $\Lambda \to \mathbf{RP}^1$ if and only if $\Lambda$ is the graph of a homeomorphism $\varphi$.*

*Moreover, in this case $G_{\Lambda,r}|_\Lambda \circ G_{\Lambda,l}^{-1}|_{\mathbf{RP}^1} = \varphi$.*

We split the proof of Corollary A.2.4 in two lemmas.

**Lemma A.2.5.** *Let $\varphi$ be an orientation preserving homeomorphism, then $G_{\phi,l}$ and $G_{\phi,r}$ extend continuously to the boundary. Moreover, $G_{\phi,r}|_\Lambda \circ G_{\phi,l}^{-1}|_{\mathbf{RP}^1} = \varphi$.*

*Remark* A.2.6. The proof essentially follows [DS24, Proposition 6.1], which adresses the case where $S$ the boundary of the convex core, but can extended to this case.

*Proof.* Fix $p_\infty = (\xi, \varphi(\xi)) \in \operatorname{graph}(\varphi)$ and let $(p_k)_{k \in \mathbb{N}}$ be a sequence of points in $S_\varphi^+$ converging to $p_\infty$. Up to extracting a subsequence, the sequence $(T_{p_k} S_\varphi^+)_{k \in \mathbb{N}}$ of spacelike support planes converges to an acausal support plane $P_\infty$ for $S_\varphi$, containing $p_\infty$.

If $P_\infty$ satisfies the condition of Lemma A.2.2, we conclude that
$$\lim_{k \to \infty} G_\varphi(p_k) = (\xi, \varphi(\xi)).$$

The point $\xi \in \mathbf{RP}^1$ and sequence $(p_k)_{k \in \mathbb{N}}$ approximating $p_\infty$ being arbitrary, $G_{\varphi,l}$ and $G_{\phi,r}$ extend continuously to $\Lambda$ as homeomorphism, and
$$G_{\phi,r}|_\Lambda \circ G_{\phi,l}^{-1}|_{\mathbf{RP}^1} = \varphi.$$

By contradiction, assume that $P_\infty$ does not satisfies the condition of Lemma A.2.2. That is $P_\infty$ is the lightlike plane orthogonal to some point $q \in \partial_\infty \mathbf{AdS}^{2,1}$, with $q \neq p_\infty$.

The two lightlike lines through $q$ intersect the graph of $\varphi$ in two points $q_1, q_2$. If $q_1 \neq q_2$, $P_\infty$ is not a support plane for $S_\varphi^+$. Then, $q_1 = q_2 = q$, that is $q \in \partial_\infty S_\varphi^+$. It follows that $q$ and $p_\infty$ are contained in a lightlike line, hence the boundary of $S_\varphi^+$ is not acausal, contradicting the hypothesis and concluding the proof. □

**Lemma A.2.7.** *Let $\Lambda$ be an achronal loop containing a non-trivial lightlike segment. Then, there exist two sequence $((p_k, P_k))_{k \in \mathbb{N}}$, $((q_k, Q_k))_{k \in \mathbb{N}}$ such that either*

$$\lim_{k \to \infty} G_{\Lambda, l}((p_k, P_k)) = \lim_{k \to \infty} G_{\Lambda, l}((q_k, Q_k)) \quad \text{or} \quad \lim_{k \to \infty} G_{\Lambda, r}((p_k, P_k)) = \lim_{k \to \infty} G_{\Lambda, r}((q_k, Q_k)).$$

*In particular, $G_{\Lambda, r} \circ G_{\Lambda, l}^{-1}$ is not defined or non-injective.*

*Proof.* Let $L$ be a non-trivial maximal lightlike segment contained in $\Lambda$. And denote by $p_\infty$ be the past-endpoint of $L$, that is

(A.1) $$\mathcal{C}^-(p_\infty) \cap \Lambda = \{p_\infty\}.$$

We want to build a suitable sequence of pointed spacelike support planes satisfying the hypothesis of Lemma A.2.2, so that $G_\Lambda$ extends to $p_\infty$.

Let us distinguish two cases, namely if $p_\infty$ is the vertex of a past-directed sawtooth or not.

*First case: $p_\infty$ is a vertex.*
Then, $S_\Lambda^+$ contains the lightlike triangle which is the convex hull of the sawtooth by Theorem A.1.4. $S_\Lambda^+$ is smooth[3] around the unique edge $\gamma$ of the triangle contained in $\mathbf{AdS}^{2,1}$.

Fix an interior point $v \in \gamma$, take an auxiliary sequence $(v_k)_{k \in \mathbb{N}}$ of points in $S_\Lambda^+$ converging to $v$, and consider the sequence $((v_k, V_k))_{k \in \mathbb{N}}$, for $V_k := T_{v_k} S_\Lambda^+$. Up extracting a subsequence, $(V_k)_{k \in \mathbb{N}}$ converges to a support acausal plane for $S_\Lambda^+$ at $v$, which we denote by $V$.

Denote $(v^*, V^*)$ the pointed plane dual to $(v, V)$. By construction, $v$ is orthogonal to $p_\infty$. Hence, $V^*$ is a spacelike support plane for $S_\Lambda^-$ at $v^*$ which contains $p_\infty$.

If $V$ is spacelike, namely $v^* \in S_\Lambda^-$, then the half-geodesic $c$ connecting $v^*$ to $p_\infty$ is contained in $S_\Lambda^-$ by convexity. It follows that $((p_k, P_k))_{k \in \mathbb{N}} := ((c(k), V^*))_{k \in \mathbb{N}}$ is a sequence of pointed support plane for $S_\Lambda^-$ satisfying Lemma A.2.2.

Otherwise, $V$ is lightlike, that is $v^*$ belongs to the sawtooth. Since $V^*$ is spacelike, $v^* = p_\infty$: indeed, only vertex of the triangle admits spacelike support planes, and the other two vertices are not orthogonal to $v$ by construction. Then, the sequence $((p_k, P_k))_{k \in \mathbb{N}}$ dual to $((v_k, V_k))_{k \in \mathbb{N}}$ satisfies Lemma A.2.2.

*Second case: $p_\infty$ is not a vertex of a past-directed sawtooth.*
By Theorem A.1.4, $S_\Lambda^+$ is smooth around $p_\infty$. Let $(p_k)_{k \in \mathbb{N}}$ be a sequence converging to $p_\infty$, let $P_k := T_{p_k} S_\Lambda^+$, and denote by $((w_k, W_k))_{k \in \mathbb{N}}$ be the dual sequence to $((p_k, P_k))_{k \in \mathbb{N}}$. By [Tre24b, Proposition 10.2.1], $w_k$ lies in the past of $p_k$. By continuity and (A.1),

$$w_\infty \in \overline{\mathcal{C}^-(p_\infty)} \cap \Lambda = \{p_\infty\}.$$

By construction, $P_\infty$ is the orthogonal of $w_\infty = p_\infty$, hence we can apply Lemma A.2.2.

In the two cases, we can conclude that there exists a sequence $((p_k, P_k))_{k \in \mathbb{N}}$ such that

$$\lim_{k \to \infty} G_\Lambda((p_k, P_k)) = (\xi, \eta),$$

for $(\xi, \eta)$ be the coordinates of $p_\infty$ in $\mathbf{RP}^1 \times \mathbf{RP}^1$.

By an analogous argument, we can build a sequence $((q_k, Q_k))_{k \in \mathbb{N}}$ which satisfies the hypothesis of Lemma A.2.2 and such that $q_k$ converges to $q_\infty$. Since lightlike lines are either vertical or horizontal lines, the coordinates of $q_\infty$ in $\mathbf{RP}^1 \times \mathbf{RP}^1$ are either $(\xi, \theta)$ or $(\theta, \eta)$. In the first case,

$$\lim_{k \to \infty} G_\Lambda((q_k, Q_k)) = \xi = \lim_{k \to \infty} G_\Lambda((p_k, P_k)).$$

---

[3] Except for the case of the Barbot crown, where $\Lambda$ consists of exactly four lightlike lines. For this case, the Gauss map can be explicitly computed, and it is the product of two geodesic in $\mathbf{H}^2 \times \mathbf{H}^2$.

In the second,
$$\lim_{k\to\infty} G_\Lambda((q_k, Q_k)) = \eta = \lim_{k\to\infty} G_\Lambda((p_k, P_k)),$$
concluding the proof. □

A.2.2. *Existence.* We can conclude that
$$\Phi_\varphi = G_{\varphi,r} \circ G_{\varphi,l}^{-1}$$
is a minimal lagrangian diffeomorphism of $\mathbf{H}^2$ extending $\varphi$. Indeed, it is minimal lagrangian map and a local diffeomorphism by [BS10], it is injective by Lemma A.2.1, hence a diffeomorphism onto. By Corollary A.2.4, it extends to the boundary as the homeomorphism $\varphi$. It follows that $\Phi_\varphi$ is surjective, hence a global diffeomorphism, concluding the existence part of the proof.

A.2.3. *Uniqueness.* In [BS18] a more general representation formula for $K$–surface is presented. We recall more specific result, for the reader convenience.

**Lemma A.2.8.** [BS18, Corollary 5.13] *Let* $\Psi\colon \mathbf{H}^2 \to \mathbf{H}^2$ *be a minimal lagrangian diffeomorphism. Then, there exists a smooth embedding* $\sigma\colon \mathbf{H}^2 \to \mathbf{AdS}^{2,1}$ *such that*
  (1) *the image* $\Sigma$ *of* $\sigma$ *is a maximal surface;*
  (2) *the Gauss map of* $\Sigma$ *induces* $\Psi$.

Let $\Psi$ be a minimal lagrangian diffeomorphism extending $\varphi$. By the representation formula, $\Psi$ is induced by the Gauss map of a maximal surface $\Sigma$, which is properly embedded since $\Psi$ extends to the boundary, hence complete by [Tre24a]. Let $\Lambda$ be the boundary of $\Sigma$. By Corollary A.2.4, the Gauss map of $\Sigma_\Lambda$ extends to a homeomorphism if and only if $\Lambda$ is the graph of a homeomorphism, hence $\Lambda$ is the graph of $\varphi$ and $\Psi = \Phi_\varphi$, concluding the proof.

Farid Diaf: Institut de Recherche Mathématique Avancée, UMR 7501, Université de Strasbourg et CNRS, 7 rue René Descartes, 67000 Strasbourg, France.
*Email address*: f.diaf@unistra.fr

Alex Moriani: Dipartimento di Matematica *Giuseppe Peano*, Università degli Studi di Torino, Via Carlo Alberto 10, 10123 Torino, Italy
*Email address*: amoriani@proton.me

Rym Smaï: Laboratoire Jean Alexandre Dieudonné, Université Côte d'Azur, Parc Valrose, 28 avenue Valrose, 06000 Nice, France
*Email address*: rym.smai@univ-cotedazur.fr

Graham Andrew Smith: Departamento de Matemática, Pontifícia Universidade Católica do Rio de Janeiro (PUC-Rio), Rio de Janeiro, Brasil
*Email address*: grahamandrewsmith@gmail.com

Enrico Trebeschi: Laboratoire Jean Alexandre Dieudonné, Université Côte d'Azur, Parc Valrose, 28 avenue Valrose, 06000 Nice, France
*Email address*: enrico.trebeschi@univ-cotedazur.fr